\documentclass[a4paper, 11pt]{article}
\usepackage{float}
\usepackage [a4paper,left=2.5cm,bottom=2.5cm,right=2.5cm,top=2.5cm]{geometry}
\usepackage{bbm}
\usepackage[english]{babel}
\usepackage{amssymb}
\usepackage{amsmath,amsthm, dsfont}
\allowdisplaybreaks
\usepackage{subcaption}
\usepackage{comment}
\usepackage{tabularx,array}
\usepackage{graphicx, url}
\usepackage{enumitem} 
\usepackage{bm}
\usepackage{bbm}
\usepackage{xcolor}
\usepackage{graphicx}
\usepackage{todonotes}
\usepackage{todonotes}
\usepackage{hyperref}
\usepackage{stmaryrd}
\usepackage{booktabs}
\usepackage{cite}
\usepackage{float}
\hypersetup{colorlinks=true, linkcolor=blue,citecolor=gray, urlcolor=gray}

\theoremstyle{plain}
\theoremstyle{plain}
\newtheorem{assumption}{Assumption}
\newtheorem{theorem}{Theorem}[section]
\newtheorem{lemma}[theorem]{Lemma}
\newtheorem{corollary}[theorem]{Corollary}
\newtheorem{proposition}[theorem]{Proposition}

\newtheorem{remark}[theorem]{Remark}

\newcommand{\E}{\mathbb{E}}

\newcommand{\R}{\mathbb{R}}
\newcommand{\N}{\mathbb{N}}

\renewcommand{\P}{\mathbb{P}}

\title{Parameter estimation for discretely observed rough  
fractional interacting particle systems}

\author{Chiara Amorino\thanks{Universitat Pompeu Fabra and Barcelona School of Economics, Department of Economics and Business, Ram\'on Trias Fargas 25-27, 08005, Barcelona, Spain. CA gratefully acknowledges financial support of PID2022-138268NB-I00/AEI/10.13039/501100011033.}%
\and
Augustin Puel\thanks{Universit\'e Côte d’Azur, Laboratoire J.A. Dieudonné, UMR CNRS 7351, Nice, 06108, France, and Centrale Méditerranée.}
\and
Radomyra Shevchenko\footnotemark[2]%
}

\date{\today}
\begin{document}
\maketitle
\begin{abstract}
We study parametric inference for interacting particle systems driven by independent fractional Brownian motions, based on discrete observations of all particles over a fixed time horizon. Our goal is to estimate an unknown drift parameter in the joint limit as the number of particles increases and the observation mesh vanishes. For \(H>1/2\), we construct contrast-based estimators by extending the discrete-time least-squares approach of Neuenkirch and Tindel to the interacting setting. For the rough regime \(H\in(1/4,1/2)\), we introduce a modified contrast and use the stochastic sewing lemma to handle the first-order approximations available from discrete observations. Our analysis combines the propagation of chaos for the associated fractional McKean–Vlasov equation with estimates controlling the temporal discretization error. When the Hurst parameter and the diffusion coefficient are known, we establish strong consistency and asymptotic normality of the proposed estimators. The asymptotic results reveal an interplay between the number of particles and the temporal discretization: the convergence rate with respect to the observation mesh changes at \(H=3/4\), while the asymptotic normality holds in both regimes considered. We further develop a joint estimation procedure for the drift parameter, the Hurst parameter, and the diffusion coefficient, extending quadratic-variation methods for fractional Brownian motion to the interacting particle framework. Numerical experiments illustrate the performance of the proposed estimators.

\end{abstract}

\tableofcontents

\section{Introduction}
Diffusion phenomena, often modeled as solutions to stochastic differential equations (SDEs), constitute a fundamental mathematical framework for describing dynamical systems subject to random fluctuations and have become an indispensable tool in probability theory and its various applications. While classical SDEs effectively describe the evolution of isolated stochastic systems, many phenomena of practical interest involve large populations of interacting agents whose dynamics are heavily influenced by collective behavior, naturally leading to systems of interacting stochastic differential equations, commonly referred to as interacting particle systems (IPS). Originating from the pioneering works of Boltzmann~\cite{boltzmann1970weitere}, Kac~\cite{kac56}, and McKean~\cite{mckean1967propagation,mck66}, IPS have become a central object of modern probability theory, leading to a rich theoretical framework that includes fundamental results on the propagation of chaos, McKean 
Vlasov limits, well-posedness, and stability~\cite{car18,szn91}. Beyond their mathematical interest, IPS serve as a vital modeling paradigm in numerous scientific disciplines, finding applications in finance~\cite{dje22,gie20}, statistical physics~\cite{dem86,mar13}, biology and self-organization~\cite{mun14}, neuroscience~\cite{bal12, tou14}, and social sciences~\cite{gom12}. More recently, IPS and mean-field models have emerged as an important framework in modern machine learning, particularly in the analysis of neural networks~\cite{rot22,chizat2018global}, optimization and Sequential Monte Carlo methods~\cite{carrillo2021consensus,grassi2021particle}, diffusion generative models~\cite{par24}, and transformer-based language models~\cite{ges25,rig25}, where they naturally describe the collective dynamics induced by self-attention mechanisms. 

Despite these extensive applications, the classical SDE framework often falls short in capturing the inherent memory and long-range dependence observed in real-world phenomena, prompting the replacement of standard Brownian motion with fractional Brownian motion (fBm). This fractional noise has proven to be highly effective in modeling processes driven by self-similar noise in fields such as finance~\cite{Com12, Ros09}, geophysics~\cite{Earthquakes, Hurst}, traffic patterns~\cite{Cag04, Will}, and medical research~\cite{Lai00}. Fractional Brownian motion provides a natural alternative when memory and long-range dependence need to be incorporated into the driving noise. For example, IPS driven by fBm are particularly well-suited to model opinion dynamics~\cite{chazelle2017well, goddard2022noisy}, such as changes in public opinion on COVID-19 vaccinations~\cite{ancona2022model, muller2022echo}, where the Hurst parameter $H \in (0, 1)$ governs the extent of memory and captures how interactions disseminate information with varying memory effects. Similarly, in the geophysical context, empirical observations suggest that non-Markovian driving noise offers a more accurate representation of turbulent motion and turbulent kinetic energy~\cite{Bossy2022, Franzke2015}. A further compelling application arises in a recent study~\cite{House2025}, which examines the growth of $N$ axon cells in vertebrate brains; while typically modeled using independent fBm paths, the authors argue that a self-repulsion effect observed in these cells necessitates an interacting particle framework. 

The numerous applications of such complex interacting models naturally require the development of a solid statistical framework. Although the statistical analysis of standard IPS has recently become an active area of research encompassing parametric, semi-parametric, and non-parametric inference~\cite{amo25, bel23,gen21,del23,sharrock2021parameter}, adapting these rigorous estimation procedures and asymptotic results for IPS driven by fractional noise remains a critical frontier. This methodological gap naturally motivates the core statistical investigations of this work. Specifically, we focus on estimating a drift parameter for the following interacting particle system driven by fractional Brownian motion:
$$
\begin{cases}
     \mathrm{d}X_t^{\theta, i,N}=b(\theta, X_t^{\theta, i,N},\mu_t^{\theta,N})\mathrm{d}t +\sigma \mathrm{d}B_t^{i,H}\\
    \mathcal{L}(X^{\theta,1,N}_0, \dots, X^{\theta,N,N}_0)=\mu_0^{\otimes N},
\end{cases}
$$
where the true parameter $\theta_0$ belongs to the interior of a compact set $\Theta$ of $\R$. The processes $(B^{i,H}_t)_{t\geq 0}$ are independent fractional Brownian motions with Hurst parameter $H \in (0, 1)$, independent of the initial value of the system, $\sigma >0$ is a volatility parameter and $\mu_t^{\theta, N}$ is the empirical measure of the system at time $t$, given by
$$
\mu_t^{\theta, N} := \frac{1}{N} \sum_{i = 1}^N \delta_{X_t^{\theta, i, N}}.
$$
For $l \ge 1$, let $\mathcal{P}_l$ denote the set of probability measures on $\R$ with a finite $l$-th moment, endowed with the Wasserstein $l$-metric
\begin{equation}
W_l(\mu, \nu) := \Big( \inf_{m \in \Gamma (\mu, \nu)} \int_{\R^2} |x - y|^l m(dx, dy) \Big)^{\frac 1 l},
\label{eq:wass}
\end{equation}
where $\Gamma(\mu, \nu)$ denotes the set of probability measures on $\R^2$ with marginals $\mu$ and $\nu$. We assume that we have access to discrete observations of all particles over a fixed time horizon $T$, meaning the observations considered in this work are given by 
\begin{equation}{\label{eq:observations}}
\big(X_{t_k}^{\theta,i,N}\big)_{k=0, ... , n}^{i=1, \ldots, N},
\end{equation}  
where the step size $\alpha_n := \frac{T}{n}$ tends to $0$ and the number of particles $N$ tends to infinity.

To the best of our knowledge, statistical inference for interacting particle systems driven by fractional Brownian motion remains largely unexplored. In the interacting setting, an initial contribution was made by~\cite{amorino2025fractional}, who investigated drift parameter estimation for models where $b$ is linear in the parameter, under continuous observation of all particles when $H > 1/2$. Subsequently, non-parametric estimation of the drift coefficient in IPS driven by fBm was addressed in~\cite{Shen}. Parallel to these developments in interacting systems, a related line of research has focused on the framework of $N$ independent and identically distributed processes driven by fractional noise. Within this setting,~\cite{marie2025computable} analyzed parameter estimation for the model considered in~\cite{amorino2025fractional} by introducing a fixed-point estimator. This approach was later extended in~\cite{amorino2025fixed} to account for multiplicative noise and the rough regime, while non-parametric estimation for $N$ i.i.d. copies driven by fBm was developed in~\cite{comte2021nonparametric}. The Brownian case, corresponding to $H=\frac12$, was studied in
\cite{amorino2023parameter}, where consistency and asymptotic
normality were established for contrast-based estimation from
discrete observations of interacting particle systems. The present
work extends this framework to the fractional setting, covering
both the Young regime $H\in(\frac12,1)$ and the rough regime
$H\in(\frac14,\frac12)$.

Classical estimation methods for diffusions and standard IPS rely heavily on the Markov and semimartingale properties of
Brownian motion. Replacing Brownian motion with fractional
noise removes both properties, making statistical inference
substantially more challenging and helping explain the
relatively limited literature in this area. However, a well-documented parallelism exists between the large-population limit ($T$ fixed, $N \rightarrow \infty$), which is the focus of this work, and the long-time limit ($N = 1$, $T \rightarrow \infty$). Leveraging this analogy allows us to draw valuable insights from the extensive literature on single-particle fractional inference. In that context, the Maximum Likelihood Estimator (MLE) and the Least Squares Estimator (LSE) are the most prevalent techniques. Although these two estimators are identical under standard Brownian noise, this equivalence breaks down in the fractional setting~\cite{HuNua10, HuNuaZho19, KleLeB02, TudVie07}. For the MLE based on continuous observations, consistency within the fractional Ornstein-Uhlenbeck (fOU) model for $H > \frac{1}{2}$ was initially established by Kleptsyna and Le Breton~\cite{KleLeB02}. Subsequent central limit theorems for the MLE in this regime were provided by Brouste and Kleptsyna, and independently by Bercu et al.~\cite{Bercu}. The theoretical scope was later broadened by Tudor and Viens~\cite{TudVie07}, who utilized Malliavin calculus to prove the consistency of both continuous and discrete MLEs across the entire range $H \in (0, 1)$ for models with linear drift. Conversely, the LSE was introduced for the fOU process in~\cite{HuNua10}, where the authors proved consistency and asymptotic normality for $H \in [\frac{1}{2}, \frac{3}{4})$. These results were successfully generalized to all $H \in (0, 1)$ in~\cite{hu2019parameter} using Malliavin calculus; however, establishing the asymptotic distribution in this framework relies strictly on the fOU structure, as the estimator belongs to the second Wiener chaos, permitting the use of limit theorems for multiple stochastic integrals. Moving beyond the one-dimensional fOU model, multidimensional systems with linear drifts were analyzed in~\cite{HuNuaZho19}, establishing LSE consistency for $H > 1/4$. A significant drawback of this continuous LSE is its reliance on Skorokhod integrals, which hinders practical computation. To overcome this limitation, Neuenkirch and Tindel~\cite{NeuTin} proposed a computable discrete-time version of the LSE for $H > \frac{1}{2}$, proving strong consistency by applying Young's inequality from rough path theory. 

The work of Neuenkirch and Tindel~\cite{NeuTin} is particularly relevant to our approach, given the aforementioned parallelism between the $N \rightarrow \infty$ and $T \rightarrow \infty$ asymptotics. Specifically, the estimator we propose is defined as the minimizer of a contrast function, which serves as a natural adaptation of their discrete-time LSE to the interacting particle framework for $H > \frac{1}{2}$. Furthermore, to handle the rougher noise regime where $H \in (\frac{1}{4}, \frac{1}{2})$, we introduce a necessary and nontrivial modification to this contrast function (see Equation~\eqref{eq:contrast rough}).

In this paper, we aim to estimate the drift parameter $\theta_0$ from the discrete observations specified in~\eqref{eq:observations}. Since the core of our strategy relies on the minimization of a contrast function, we require asymptotic results on the particle system as the number of particles $N \rightarrow \infty$. In the standard Brownian motion case, it is well known that under appropriate assumptions on the drift coefficient, the interacting particle system exhibits a propagation of chaos phenomenon, yielding the McKean-Vlasov equation as the mean-field limit (see~\cite{szn91}). This phenomenon extends to the case of fractional noise; we provide a concise proof of this extension in Section~\ref{prop chaos}. Consequently, we define the limit system via the following McKean-Vlasov equation:
\begin{align}\label{McKean-Vlasov}
\mathrm{d}\bar X_t^i=b(\theta_0,\bar X_t^i,\bar \mu_t)\mathrm{d}t+\sigma \mathrm{d}B^{i,H}_t
\end{align}
where $\bar \mu_t$ is the law of $\bar X_t^i$ and $\bar X_0^i=X^i_0$ is the initial condition. In this equation, the drift coefficient depends not only on the spatial position but also on the current marginal law of the process. The well-posedness of this equation is established in Lemma 19 of~\cite{pathwiseMcKeanvlasov}.

To construct our estimators, we must also handle the time discretization of the process. To do so, we rely on tools from rough path theory~\cite{friz2014course}. In the regime $H>\frac{1}{2}$, the pathwise stochastic integration is well understood and can be handled using the standard Sewing Lemma. However, when $H\in(\frac{1}{4},\frac{1}{2})$, defining pathwise stochastic integrals becomes significantly more complex, as rough path theory typically requires iterated integrals (second-order information) to construct them. Since this additional information is unavailable from our discrete data, we are constrained to use a first-order approximation. To overcome this theoretical hurdle, we leverage the Stochastic Sewing Lemma introduced by Lê~\cite{StochasticSewingLe}, and more specifically its improved version by Perkowski and Matsuda~\cite{stochasticsewingupgrade}. This powerful tool exploits stochastic cancellations to gain regularity, allowing for a construction of Stratonovich-type stochastic integrals of the fBm via a well-chosen Riemann approximation. We generalize this construction to measure-dependent integrals, obtaining a result of independent interest for fractional IPS. Specifically, when $H \in (1/4, 1/3]$, the measure dependence can be handled due to the extra regularity of the McKean–Vlasov law; indeed, we show that it exhibits $2H$-regularity even though the particles only have $H$-regularity. More generally, our strategy exploits the improved regularity
in expectation of centered fBm-driven processes: although their
sample paths have limited H\"older regularity, their zero-mean
property yields sharper estimates.

Building on this mathematical framework, our paper establishes several key theoretical guarantees for the estimation of $\theta_0$. First, assuming the Hurst and diffusion parameters are known, we prove the strong consistency of our contrast-based estimators in both regimes ($H>\frac{1}{2}$ and $H\in(\frac{1}{4},\frac{1}{2})$) as $N \rightarrow \infty$ and the time step goes to zero. Second, we establish the asymptotic normality of these estimators. We find that in $N$, the estimators converge with a rate $\sqrt{N}$, typical for the IPS literature, and also established in~\cite{amorino2025fractional} for the fractional setup, while the rate in $n$ changes behavior at the threshold $H=\frac{3}{4}$, similarly to the drift estimator for the fOU in the $T\to\infty$ asymptotics~\cite{hu2019parameter}. In contrast to the long-time asymptotics of the fOU model,
our estimators are asymptotically normal throughout the
fractional regimes considered in this paper.

Finally, we demonstrate that our approach can be successfully extended to perform the joint estimation of all model parameters by extending Coeurjolly's results on weighted quadratic variations~\cite{estimationparameterfBm} for a discretely observed fBm to the interacting particle framework, allowing us to achieve a faster convergence rate. Here, again, the rough case requires a separate treatment and is handled through a regularity gain for a particular choice of weights. We conclude with numerical experiments illustrating the performance of our estimators. 

The remainder of the paper is organized as follows. Section~\ref{section model assumptions} introduces the contrast functions for both regimes, $H>\frac{1}{2}$ and $H\in(\frac{1}{4},\frac{1}{2})$, the associated estimators, and the main assumptions of our model. Section~\ref{main results} presents our theoretical results: it establishes the propagation of chaos (Section~\ref{prop chaos}), followed by the consistency (Section~\ref{section consistency}) and the asymptotic normality (Section~\ref{section CLT}) of the estimators. Section~\ref{section joint estimation} extends this analysis by detailing the joint estimation of all parameters. Numerical experiments are provided in Section~\ref{section numerics}. Section~\ref{section proofs} contains the proofs of the key theorems, while auxiliary results are deferred to Section~\ref{section preliminary results}. Finally, Sections~\ref{app:malliavin_fbm} and~\ref{app: Lions derivative} provide the necessary technical background on Malliavin calculus and differential calculus on probability measures, respectively.

\subsection*{Notations}

Throughout the paper, $c$ denotes a generic positive constant whose
value may change from one occurrence to another. We denote by
$\partial_x f$ the spatial derivative of a function $f$ and by
$\partial_\mu f$ its Lions derivative with respect to the measure
variable (see Appendix~\ref{app: Lions derivative} for details).\\
For a process $Y:[0,T]\to\R$, we use the notation
\[
Y_{st}:=Y_t-Y_s.
\]
For a function
$f:\Theta\times\R\times\mathcal{P}_2(\R)\to\R$, we define
\[
\delta_{\theta,\tilde\theta} f(\cdot,x,\mu):=f(\theta,x,\mu)-f(\tilde\theta,x,\mu).
\]

\section{Model and Assumptions}\label{section model assumptions}


Our aim is to estimate the unknown parameter $\theta_0$ from discrete observations of particle trajectories in the asymptotic regime $N,n\to\infty$. We first consider the case $H>\frac12$, which provides a simpler setting in which to introduce our strategy. We then turn to the rougher regime $H\in(\frac14,\frac12)$. For $H\in(\frac12,1)$, we consider the contrast function
\begin{equation}\label{eq:contrast Young}
Q_n^N(\theta)=\frac{1}{N\alpha_n}\sum_{i=1}^N\sum_{k=0}^{n-1}\left (\left(X^{i,N}_{t_kt_{k+1}}-\alpha_n b(\theta,X^{i,N}_{t_k},\mu^N_{t_k})\right)^2-\sigma^2\alpha_n^{2H}\right).
\end{equation}
This contrast is inspired by the approach of Neuenkirch and Tindel~\cite{NeuTin}. It can be viewed as a least squares type error with a correction term of order $\alpha_n^{2H}$, which removes the contribution of the second-order moments of the fractional Brownian motion. For $H\in(\frac14,\frac12)$, we instead consider 
\begin{equation}\label{eq:contrast rough} 
Q_n^N(\theta)=\frac{1}{N\alpha_n} \sum_{i=1}^N\sum_{k=0}^{n-1}\left(\left(X^{i,N}_{t_kt_{k+1}}-\frac{\alpha_n}{2}
\left(b(\theta,X^{i,N}_{t_k},\mu^N_{t_k})+b(\theta,X^{i,N}_{t_{k+1}},\mu^N_{t_{k+1}})\right)\right)^2-\sigma^2\alpha_n^{2H}\right).
\end{equation}

The two contrast functions differ only in the approximation used for the drift integral over each time interval. When $H>\frac12$, the choice of approximation is not essential for our analysis. In the rough regime ($H\leq\frac12$), however, the choice becomes relevant. In particular, midpoint-type approximations can benefit from stochastic cancellations, as shown, for instance, in Proposition 3.5 of \cite{stochasticsewingupgrade}.
For both contrast functions, we define the estimator
\[
\theta_n^N \in \operatorname*{arg\,min}_{\theta\in\Theta} \left|Q_n^N(\theta)\right|^2.
\]
Following the idea of Neuenkirch and Tindel~\cite{NeuTin}, who consider the absolute value of the contrast, we minimize its square. This modification does not affect the consistency argument and leads to a more convenient formulation of the proof of the central limit theorem.

We now introduce the assumptions required for our asymptotic results. The following two assumptions are standard and guarantee existence and uniqueness of the IPS solution.

\begin{assumption}[Bounded moments]\label{Bdd moment}
For every $k\geq 1$, there exists a constant $C_k>0$ such that
\[
\int_{\R}|x|^k\,\mu_0(\mathrm{d}x)\leq C_k.
\]
\end{assumption}

\begin{assumption}[Lipschitz conditions]\label{Lip cond}
The drift coefficient $b$ is Lipschitz continuous in $(x,\mu)$, uniformly in $\theta\in\Theta$. More precisely, for every $\theta\in\Theta$, there exists $C>0$ such that, for all $(x,\mu),(y,\nu)\in\R\times\mathcal{P}_2(\R)$,
\[
|b(\theta,x,\mu)-b(\theta,y,\nu)| \leq C\bigl(|x-y|+W_2(\mu,\nu)\bigr).
\]

Moreover, $b$ is Lipschitz continuous in the parameter $\theta$, with polynomial growth in $(x,\mu)$. Namely, there exist constants $C>0$ and $L,K\in\N$ such that, for all $\theta,\tilde\theta\in\Theta$ and $(x,\mu)\in\R\times\mathcal{P}_2(\R)$,
\[
|b(\theta,x,\mu)-b(\tilde\theta,x,\mu)| \leq C|\theta-\tilde\theta| \left( 1+|x|^K+W_2^L(\delta_0,\mu)\right).
\]
\end{assumption}

The Lipschitz assumption may potentially be relaxed. For instance, one could consider locally Lipschitz drift coefficients with polynomial growth under  additional structural assumptions on the interaction. Although our propagation-of-chaos argument does not directly cover such settings, related results for singular drift coefficients can be found in~\cite{singulardrift}; see also the discussion and examples in Section 2.2 of~\cite{amorino2023parameter}.

\begin{assumption}[Identifiability condition]\label{Identif cond}
For $\theta\in\Theta$, we assume that
\[
\mathbb{E}\left[\int_0^T |\delta_{\theta,\theta_0}b(\cdot,\bar X_s,\bar\mu_s)|^2 \,\mathrm{d}s-2\sigma \int_0^T \delta_{\theta,\theta_0}b(\cdot,\bar  X_s,\bar\mu_s)\circ\mathrm{d}B_s^H\right]=0,
\]
implies $\theta=\theta_0$. The stochastic integral is either a Young integral when $H\in(\frac12,1)$ or a Stratonovich integral defined in Lemma~\ref{lemma Strat integral} for $H\in (\frac14,\frac12)$. 
\end{assumption}

     Although the above assumption is usually difficult to check by hand, it holds for the fOU process, as shown in Section~\ref{section numerics}. In some models, the identifiability condition can be verified
locally by restricting the parameter space to a neighborhood
of $\theta_0$, provided that $\theta_0$ is an isolated zero
of the corresponding identifiability function.

\begin{assumption}[Non-degeneracy]\label{ass invertibility}
    We assume that:
    \[
    \E[\int_0^T \partial_\theta b(\theta_0,\bar X_s,\bar\mu_s)\circ \mathrm{d}B_s^H]\neq0 .
    \]
\end{assumption}
This assumption is required to obtain a non-zero variance in the CLT. It is not a restrictive condition, meaning that, except for a very specific structure of the drift $b$, it is often verified.
\begin{assumption}[Rate of convergence]\label{CV rate step size 1}
We assume that
\[
\frac{n^{3-4H}\log^2(n)}{N} \xrightarrow[N,n\to\infty]{}0.
\]
\end{assumption}

This condition is restrictive only when $H\leq\frac34$, reflecting the change in the behavior of the second-order moments of fractional Brownian motion. This phenomenon is illustrated in Proposition~\ref{L2_Estimate_fBm}.

\begin{assumption}[Rate of convergence for the Young CLT]
\label{CV rate step size 2}
We assume that
\[
N\alpha_n \xrightarrow[N,n\to\infty]{}0.
\]
\end{assumption}

This condition is also present in the Brownian case; see~\cite{amorino2023parameter}. However, this condition is compatible with Assumption~\ref{CV rate step size 1} only when $H>\frac12$. Therefore, for $H\in(\frac14,\frac12)$, we require different conditions.

\begin{assumption}[Rate of convergence for the rough CLT ($H\in(1/3,1/2)$)]
\label{CV rate step size 3}
We assume that
\[
N\alpha_n^{1+2H} \xrightarrow[N,n\to\infty]{}0.
\]
\end{assumption}

\begin{assumption}[{Rate of convergence for the rough CLT ($H\in(1/4,1/3]$)}] \label{CV rate step size 4}
We assume that
\[
N\alpha_n^{1+4H} \xrightarrow[N,n\to\infty]{}0.
\]
\end{assumption}

Assumption \ref{CV rate step size 1}  is imposed throughout, including for the
consistency results, whereas Assumptions~\ref{ass invertibility},~\ref{CV rate step size 2},~\ref{CV rate step size 3} and \ref{CV rate step size 4} are required only for the asymptotic normality results.

\begin{assumption}\label{assumptions derivative}
The function $f:\Theta\times\R\times\mathcal{P}_2(\R)\to\R$ is twice continuously differentiable with respect to $\theta$, uniformly in the remaining variables. The function f and its derivatives in $\theta$ are Lipschitz continuous in $(x,\mu)$, uniformly in $\theta$, and Lipschitz continuous in $\theta$ with polynomial growth in $(x,\mu)$.
\end{assumption}

\begin{assumption}\label{strong assumption on derivative}
The function $f:\Theta\times\R\times\mathcal{P}_2(\R)\to\R$ is three times continuously differentiable with respect to $x$, with bounded derivatives uniformly in $(\theta,\mu)$ and it admits a Lions derivative in $\mu$ uniformly in $(\theta,x)$ and the derivative is bounded.
\end{assumption}

\begin{assumption}\label{very strong assumption on derivative}
The function $f:\Theta\times\R\times\mathcal{P}_2(\R)\to\R$ is five times continuously differentiable with respect to $x$, with bounded derivatives uniformly in $(\theta,\mu)$ and is twice Lions differentiable in $\mu$ uniformly in $(\theta,x)$  with $\partial_\mu f(\theta,x,\nu)(z)$ differentiable in $z$ and all the derivatives are bounded. Finally, assume that $\partial_x f$ also admits a Lions derivative with $\partial_\mu\partial_x f(\theta,x,\mu)(z)$ differentiable in $z$ and bounded.
\end{assumption}

The last three assumptions will be applied primarily to the drift coefficient $b$, either in the consistency analysis or in the central limit theorem, in particular in the rough regime $H\in(\frac14,\frac12)$, where more precise estimates are required.

Furthermore, in practical applications, the discretization step size cannot always be freely chosen. In scenarios where Assumption~\ref{CV rate step size 1} fails to hold, provided the drift is sufficiently smooth to satisfy Assumption~\ref{very strong assumption on derivative}, we can alternatively rely on Assumption~\ref{CV rate step size 4} for all $H > \frac{1}{4}$.

\section{Main results}\label{main results}

\subsection{Propagation of chaos}\label{prop chaos}

The propagation of chaos is a classical and well-studied phenomenon in interacting particle systems. A wide range of results are available under various  assumptions on the coefficients and noise. In our setting, the strong moment and regularity assumptions imposed on the drift allow us to obtain a  relatively short and direct proof. In particular, we derive an almost sure rate of convergence arbitrarily close to $N^{-1/2}$, which is sufficient for our purposes. The result holds for all $H\in(0,1)$.

To establish propagation of chaos, we use the coupling argument introduced by Sznitman~\cite{szn91}. More precisely, for each fractional Brownian motion and each initial condition, we associate the solution of the corresponding McKean-Vlasov equation:
\begin{align}\label{McKean-Vlasov equation}
\mathrm{d}\bar X_t^i=b(\theta_0,\bar X_t^i,\bar\mu_t)\,\mathrm{d}t+\sigma\,\mathrm{d}B_t^{i,H},\qquad\bar X_0^i=X_0^i.
\end{align}
Since the initial conditions and the fractional Brownian motions are independent, the processes $(\bar X^i)_{i\geq1}$ are independent and identically distributed. We define their empirical measure by
\[
\bar\mu_t^N=\frac{1}{N}\sum_{i=1}^N\delta_{\bar X_t^i}.
\]

Since the particles are i.i.d., we can use a large deviation estimate for their empirical measure. Using the concentration inequality given in Theorem 2 of~\cite{CVmeasure}, we obtain the following result.

\begin{theorem}\label{a.s limit of empirical measure}
Under Assumptions~\ref{Bdd moment} and~\ref{Lip cond}, for every $\gamma\in(0,\frac12)$, there exists a random variable $\eta_\gamma$ such that
\[
\sup_{t\in[0,T]}W_2(\bar\mu_t^N,\bar\mu_t) \leq \frac{\eta_\gamma}{N^\gamma}
\]
almost surely. Moreover, $\eta_\gamma\in L^p(\Omega)$ for every $p\geq1$.
\end{theorem}

The proof is postponed to Section~\ref{section proofs prop chaos}.

\begin{remark}
The same argument applies when the state space is $\R^d$. In that case, one can work with the $p$-Wasserstein distance for $p>d/2$ in order to apply Theorem 2 of~\cite{CVmeasure}. For $p\geq2$, the $2$-Wasserstein distance is controlled by the $p$-Wasserstein distance.
\end{remark}

The previous estimate yields the propagation of chaos for the interacting particle system. Indeed, since the noise is additive and the drift coefficient is Lipschitz, the distance $W_2(\mu^N,\bar\mu^N)$ can be used to control the difference between the particle system and its McKean-Vlasov approximation. We therefore obtain the following proposition, whose proof is given in Section~\ref{section proofs prop chaos}.

\begin{proposition}\label{prop of chaos}
Under Assumptions~\ref{Bdd moment} and~\ref{Lip cond}, for every $\gamma\in(0,\frac12)$, there exists a random variable $\tilde\eta_\gamma>0$, independent of $N$, such that
\[
\|X^{i,N}-\bar X^i\|_{\infty,[0,T]} \leq \frac{\tilde\eta_\gamma}{N^\gamma}, \qquad \|W_2(\mu^N,\bar\mu)\|_{\infty,[0,T]} \leq \frac{\tilde\eta_\gamma}{N^\gamma},
\]
almost surely.
\end{proposition}

\begin{remark}
The Gaussian nature of the driving noise is not essential for this argument. The proof only requires the corresponding empirical-measure estimate together with sufficient moment and Hölder-regularity properties of the driving noise.
\end{remark}

The previous proposition can be combined with an interpolation inequality to obtain convergence in Hölder norms.

\begin{corollary}\label{prop chaos Holder Norms}
Under Assumptions~\ref{Bdd moment} and~\ref{Lip cond}, for every $H\in(0,1)$, $\varepsilon>0$, and $\gamma\in(0,\frac12)$, there exists a random variable $\tilde\eta_{\gamma,\varepsilon}>0$, independent of $N$, such that
\[
\|X^{i,N}-\bar X^i\|_{H-\varepsilon,[0,T]} \leq \frac{\tilde\eta_{\gamma,\varepsilon}}{N^\gamma}, \qquad \|W_2(\mu^N,\bar\mu)\|_{H-\varepsilon,[0,T]} \leq \frac{\tilde\eta_{\gamma,\varepsilon}}{N^\gamma},
\]
almost surely.
\end{corollary}

The proof is postponed to Section~\ref{section proofs prop chaos}.

\subsection{Consistency of the estimator}\label{section consistency}

We now establish the consistency of the estimator in both regimes $H>\frac12$ and $H\in(\frac14,\frac12)$. Although the two cases rely on different technical tools, the overall strategy is inspired by Neuenkirch and Tindel~\cite{NeuTin}, who considered the long-time asymptotics of a single particle.

\subsubsection{The Young case}

\label{Young case}

We first consider the case $H>\frac12$ and the contrast function defined in~\eqref{eq:contrast Young}. The consistency proof is based on a standard result concerning minimizers of random functions, see e.g.~\cite{kasonga1988}.

\begin{proposition}\label{CV minimizers}
Assume that a family of random variables $L_n(\theta)$, $n\in\N$, $\theta\in\Theta$, satisfies the following conditions:
\begin{enumerate}[label=\arabic*)]
    \item Almost surely, $L_n(\theta)\to L(\theta)$ uniformly in $\theta\in\Theta$ as $n\to\infty$.
    \item The limit $L$ is deterministic and satisfies $L(\theta_0)\leq L(\theta)$ for every $\theta\in\Theta$.
    \item $L(\theta)=L(\theta_0)$ if and only if $\theta=\theta_0$.
\end{enumerate}
Then
\[
\hat\theta_n\longrightarrow\theta_0 \qquad\text{almost surely},
\]
where
\[
L_n(\hat\theta_n)=\min_{\theta\in\Theta}L_n(\theta).
\]
\end{proposition}

To apply this result, we need to characterize the limiting behavior of the contrast function. We first study the contrast evaluated at the true parameter.

\begin{proposition}\label{CV contrast in theta 0}
For $H>\frac12$, under Assumptions~\ref{Bdd moment},~\ref{Lip cond}, and~\ref{CV rate step size 1},
\[
Q_n^N(\theta_0)\xrightarrow[N,n\to\infty]{}0
\]
almost surely.
\end{proposition}

Note that the additional rate condition in Assumption \ref{CV rate step size 1} is not required in the Brownian case~\cite{amorino2023parameter}. Next, we study the difference between the contrast evaluated at $\theta$ and the true parameter $\theta_0$.

\begin{theorem}\label{Diff CV}
For $H>\frac12$, under Assumptions~\ref{Bdd moment} and \ref{Lip cond},
\begin{align*}
\sup_{\theta\in\Theta}\Bigg|Q_n^N(\theta)-Q_n^N(\theta_0)-\mathbb{E}\Bigg[&\int_0^T\left|\delta_{\theta,\theta_0}b(\cdot,\bar X_s,\bar\mu_s)\right|^2\,\mathrm{d}s -2\sigma\int_0^T\delta_{\theta,\theta_0}b(\cdot,\bar X_s,\bar\mu_s)\circ\mathrm{d}B_s^H\Bigg]\Bigg|\xrightarrow[N,n\to\infty]{}0
\end{align*}
almost surely.
\end{theorem}

Combining the previous two results yields the consistency of the estimator.

\begin{theorem}\label{Consistency Young case}
For $H>\frac12$, assume that \ref{Bdd moment},~\ref{Lip cond},~\ref{Identif cond}, and \ref{CV rate step size 1} hold. Then, for every $\theta_0\in\Theta$,
\[
\theta_n^N\xrightarrow[N,n\to\infty]{}\theta_0
\]
almost surely.
\end{theorem}

\begin{proof}
Set $L_n^N(\theta)=|Q_n^N(\theta)|^2$
and
\[
L(\theta)=\left|\mathbb{E}\left[\int_0^T\left|\delta_{\theta,\theta_0}b(\cdot,\bar X_s,\bar\mu_s)\right|^2\,\mathrm{d}s-2\sigma\int_0^T\delta_{\theta,\theta_0} b(\cdot,\bar X_s,\bar\mu_s) \circ\mathrm{d}B_s^H \right]\right|^2.
\]
The first condition of Proposition~\ref{CV minimizers} follows from Theorem~\ref{Diff CV}. The function $L$ is deterministic since it isdefined as an expectation. Moreover, $L(\theta_0)=0$ by Proposition~\ref{CV contrast in theta 0}, so that $\theta_0$ is a minimizer of $L$. Finally, the identifiability condition in Assumption~\ref{Identif cond} ensures that this minimizer is unique.
\end{proof}

\subsubsection{The rough case}\label{rough section}

We now consider the rough regime $H\in(\frac14,\frac12)$. In this case, the contrast function used in the Young regime is no longer appropriate. More precisely, Proposition~\ref{CV contrast in theta 0} does not hold in this range. To understand the reason, we use the following expansion of the
particle trajectories:
\begin{align}\label{Taylor particule}
X^{i,N}_{t_k t_{k+1}}&=\alpha_n b(\theta_0,X^{i,N}_{t_k},\mu^N_{t_k})+\sigma B^{i,H}_{t_k t_{k+1}}+r_k^{i,N},
\end{align}
where
\begin{equation}\label{eq: remainder young}
r_k^{i,N}:=\int_{t_k}^{t_{k+1}}\left(b(\theta_0,X_s^{i,N},\mu_s^N)-b(\theta_0,X_{t_k}^{i,N},\mu_{t_k}^N)\right)\mathrm{d}s.
\end{equation}

Substituting this expansion into the contrast evaluated at the true parameter gives
\begin{equation}\label{expansion QnN(theta0)}
Q_n^N(\theta_0)=\frac{\sigma^2}{N\alpha_n}\sum_{i=1}^N\sum_{k=0}^{n-1}\left(|B^{i,H}_{t_k t_{k+1}}|^2-\sigma^2\alpha_n^{2H}\right)+\frac{2\sigma}{N\alpha_n}\sum_{i=1}^N\sum_{k=0}^{n-1}r_k^{i,N}B^{i,H}_{t_k t_{k+1}}+\frac{1}{N\alpha_n}\sum_{i=1}^N\sum_{k=0}^{n-1}|r_k^{i,N}|^2.
\end{equation}

The first and last terms can be treated as in the Young case. However, the second term involving $r_k^{i,N}B^{i,H}_{t_k t_{k+1}}$, has a different behavior. Even in the simple case $b(\theta,x,\mu)=\theta x$, direct computations using the regularity of the fractional Brownian motion and of the particle paths show that this term is of order $\alpha_n^{2H-1}$, that no longer tends to zero for $H\leq \frac12$.

A way to overcome this difficulty is to use an average approximation instead of a left-point one. This is motivated by the following lemma.

\begin{lemma}\label{lemma Strat integral}
Let $H\in(\frac14,\frac12)$ and assume that $f$ is $C^1$ in space with a Lipschitz derivative. When $H\in(\frac14,\frac13]$, assume further that $f$ and $\partial_x f$ admit a Lions derivative $\partial_\mu f$ and $\partial_\mu \partial_xf$, finally we assume that $\partial_\mu f$ and $\partial_\mu\partial_x f(\theta,x,\nu)(z)$ are differentiable in $z$ with a bounded derivatives. Then, for every $T>0$ and every partition $\pi$ of $[0,T]$, one can define the "Stratonovich integral",
\[
\int_0^T f(\theta,\bar X_s^i,\bar\mu_s) \circ\mathrm{d}B_s^{i,H}
\]
through the limit
\begin{equation}\label{def Strato rough}
\lim_{|\pi|\to0} \sum_{[s,t]\in\pi} \frac{f(\theta,\bar X_t^i,\bar\mu_t)+f(\theta,\bar X_s^i,\bar\mu_s)}{2}B^{i,H}_{st}=\int_0^Tf(\theta,\bar X_s^i,\bar\mu_s)\circ\mathrm{d}B_s^{i,H}.
\end{equation}
The convergence holds in $L^p(\Omega)$ for every $p\geq2$. Moreover,
there exists $\varepsilon>0$ such that
\[
\left\|\sum_{[s,t]\in\pi}\frac{f(\theta,\bar X_t^i,\bar\mu_t)+f(\theta,\bar X_s^i,\bar\mu_s)}{2}B^{i,H}_{st}-\int_0^T f(\theta,\bar X_s^i,\bar\mu_s)\circ\mathrm{d}B_s^{i,H}\right\|_{L^p(\Omega)}\lesssim |\pi|^\varepsilon,
\]
where the implicit constant is independent of $\theta$.

Consequently, if the discretization step is sufficiently small, the convergence in~\eqref{def Strato rough} also holds almost surely by the Borel--Cantelli lemma.
\end{lemma}

The proof relies on the Stochastic Sewing Lemma, introduced by Lê~\cite{StochasticSewingLe}. More precisely, our argument is based on the improved version developed by Perkowski and Matsuda~\cite{stochasticsewingupgrade}. The present result can be viewed as a generalization of Proposition 3.6 in~\cite{stochasticsewingupgrade}, where a function of fractional Brownian motion is integrated against the same fractional Brownian motion. Moreover, the integral obtained through this approximation coincides with the corresponding integral defined in the rough path framework using the canonical rough path of the fBm; see, for instance,~\cite{SelectedAspectsofFractionalBM}, or also the Russo-Valois integral defined in~\cite{russo1993forward}. 

We recall the definition of the contrast function for $H\in(\frac14,\frac12)$:
\begin{equation*}
\begin{aligned}
Q_n^N(\theta)= \frac{1}{N\alpha_n} \sum_{i=1}^N\sum_{k=0}^{n-1} \Bigg(\left( X^{i,N}_{t_k t_{k+1}}-\frac{\alpha_n}{2} \Big[b(\theta,X^{i,N}_{t_{k+1}},\mu^N_{t_{k+1}})+b(\theta,X^{i,N}_{t_k},\mu^N_{t_k})\Big]\right)^2-\sigma^2\alpha_n^{2H}\Bigg).
\end{aligned}
\end{equation*}

As in the Young case, the estimator is defined by
\[
\theta_n^N=\operatorname*{arg\,min}_{\theta\in\Theta}|Q_n^N(\theta)|^2.
\]

We can now establish the convergence of the contrast function.

\begin{theorem}\label{CV contrast rough case}
Let $H\in(\frac14,\frac12)$.  Under Assumptions~\ref{Bdd moment},~\ref{Lip cond},~\ref{CV rate step size 1}, and additionally~\ref{strong assumption on derivative} \textup{(}resp.~\ref{very strong assumption on derivative}\textup{)} on $b$ for $H\in(\frac13,\frac12)$ \textup{(}resp.\ $H\in(\frac14,\frac13]$\textup{)}, we have, for $N,n \rightarrow \infty$:

\begin{enumerate}
    \item $\sup_{\theta\in\Theta}
    \Bigg|Q_n^N(\theta)-Q_n^N(\theta_0)-\mathbb{E}\Bigg[\int_0^T\left|\delta_{\theta,\theta_0}b(\cdot,\bar X_s,\bar\mu_s)\right|^2\,\mathrm{d}s-2\sigma     \int_0^T \delta_{\theta,\theta_0} b(\cdot,\bar X_s,\bar\mu_s) \circ\mathrm{d}B_s^H \Bigg]\Bigg|\rightarrow 0$
    almost surely;
    \item
    $Q_n^N(\theta_0)\rightarrow 0$
    almost surely.
\end{enumerate}
\end{theorem}

Although the first statement has the same form as the corresponding result in the Young regime, the stochastic integral that appears in the limit is defined as in Lemma~\ref{lemma Strat integral}. Moreover, the second statement requires a stronger regularity assumption on the drift since sharper estimates on the remaining terms are needed to prove convergence to zero.

We can now deduce the consistency of the estimator.

\begin{theorem}\label{consistency rough case}
Let $H\in(\frac14,\frac12)$.  Under Assumptions~\ref{Bdd moment}, \ref{Lip cond}, \ref{Identif cond}, \ref{CV rate step size 1} and additionally~\ref{strong assumption on derivative} \textup{(}resp.~\ref{very strong assumption on derivative}\textup{)} on $b$ for $H\in(\frac13,\frac12)$ \textup{(}resp.\ $H\in(\frac14,\frac13]$\textup{)}, we have:
\[
\theta_n^N \xrightarrow[N,n\to\infty]{} \theta_0
\]
almost surely, for every $\theta_0\in\Theta$.
\end{theorem}

\begin{proof}
The proof is identical to that of Theorem \ref{Consistency Young case}, with Theorem \ref{CV contrast rough case} replacing the corresponding convergence results in the Young regime.
\end{proof}

\subsection{A Central Limit Theorem}\label{section CLT}

We now establish the asymptotic normality of the estimator in both regimes, $H>\frac12$ and $H\in(\frac13,\frac12)$. The overall strategy is the same in both cases, although the assumptions and the normalization differ. Unless otherwise specified, $Q_n^N$ denotes the appropriate contrast function for the regime under consideration. The proofs are given in Section~\ref{section proofs CLT}.

Since the estimator is defined implicitly as a minimizer of a contrast function, in general there is no explicit expression for $\theta_n^N$ available. We therefore  use the fact that the estimator is a critical point of the squared contrast and apply a Taylor expansion. Following the approach used, for example, in~\cite{amorino2023parameter}, we obtain 
\begin{align}\label{Taylor constrast function}
(\theta_{n}^N-\theta_{0}) \int_0^1 \partial_{\theta}^2 (Q_n^N)^2 \bigl(s\theta_0+(1-s)\theta_n^N\bigr)\,\mathrm{d}s=-2Q_n^N(\theta_0)\partial_{\theta}Q_n^N(\theta_0).
\end{align}

We therefore need to study the asymptotic behavior of the two terms on the right-hand side. The contrast $Q_n^N(\theta_0)$ converges to zero by the consistency analysis, while we establish a central limit theorem for its suitably normalized version. We also prove the convergence of the first derivative. The joint convergence is then obtained using the Cram\'er--Wold theorem (see Theorem 29.4 in~\cite{billingsley1995probability}).

The analysis of $Q_n^N(\theta_0)$ is based on the expansion~\eqref{expansion QnN(theta0)}. The leading contribution is given by a sum of independent random variables, while the remaining terms arise from the Taylor remainder and are of smaller order. Since both $N$ and $n$ tend to infinity, we prove the following triangular-array central limit theorem based on a Lindeberg condition.

\begin{theorem}\label{Linderberg CLT}
Let $(Y_n^i)_{i\geq1,n\geq1}$ be a triangular array of real-valued random variables such that, for every $n$, the random variables $(Y_n^i)_{i\geq1}$ are i.i.d. and centered, with variance
\[
\operatorname{Var}(Y_n^1)=\lambda_n^2 \xrightarrow[n\to\infty]{} \lambda^2. 
\]
Assume moreover that, for every $\delta>0$, $\mathbb{E}\left[|Y_n^1|^2\mathbf{1}_{\{|Y_n^1|\geq\sqrt{N}\delta\}}\right]\xrightarrow[N,n\to\infty]{} 0$.
Then
\[
\frac{1}{\sqrt{N}} \sum_{i=1}^N Y_n^i \xrightarrow[N,n\to\infty]{\mathcal{D}}\mathcal{N}(0,\lambda^2).
\]
\end{theorem}
The proof can be found in Section~\ref{section proofs CLT}.\\
We can now state the asymptotic behavior of the contrast evaluated at the true parameter.

\begin{theorem}\label{CLt contrast function}
For $H>\frac12$, under Assumptions \ref{Bdd moment},~\ref{Lip cond},~\ref{CV rate step size 1},~\ref{CV rate step size 2}, and \ref{assumptions derivative} applied to $b$, we have
\[
\sqrt{N\gamma_n^2}\,Q_n^N(\theta_0)\xrightarrow[N,n\to\infty]{\mathcal{D}}\mathcal{N}(0,\sigma^2c_H^2),
\]
where
\[
\gamma_n^2=
\begin{cases}
n^{4H-3},& \frac12<H<\frac34,\\[3mm]
\log(n)^{-1},& H=\frac34,\\[3mm]
1,& H>\frac34,
\end{cases}
\]
and $c_H$ is given in Proposition~\ref{L2_Estimate_fBm}.

Similarly, for $H\in(\frac14,\frac12)$, under Assumptions \ref{Bdd moment},~\ref{Lip cond}, \ref{CV rate step size 1},~\ref{CV rate step size 3} \textup{(}resp.~\ref{CV rate step size 4}\textup{)},~\ref{assumptions derivative} applied to b, and additionally~\ref{strong assumption on derivative} \textup{(}resp.~\ref{very strong assumption on derivative}\textup{)} on $b$ for $H\in(\frac13,\frac12)$ \textup{(}resp.\ $H\in(\frac14,\frac13]$\textup{)}, we have
\[
\sqrt{Nn^{4H-3}}\,Q_n^N(\theta_0)\xrightarrow[N,n\to\infty]{\mathcal{D}}\mathcal{N}(0,\sigma^2c_H^2).
\]

\end{theorem}

It remains to study the asymptotic behavior of the first and second derivatives of the contrast. Most of the required estimates were already established in the consistency analysis.

\begin{proposition}\label{assymptotic first and second order derivative}
Under Assumptions~\ref{Bdd moment},~\ref{Lip cond},~\ref{CV rate step size 1}, \ref{assumptions derivative}, as well as Assumption~\ref{strong assumption on derivative} \textup{(}resp.~\ref{very strong assumption on derivative}\textup{)} on $b$ for $H\in(\frac13,\frac12)$ \textup{(}resp.\ $H\in(\frac14,\frac13]$\textup{)} and $\partial_{\theta}b,\partial^2_{\theta}b$ satisfies the hypotheses of Lemma~\ref{lemma Strat integral}, we have:
\begin{align*}
&\sup_{\theta\in\Theta} \Bigg|\partial_{\theta}Q_n^N(\theta)-2\mathbb{E}\Bigg[ \int_0^T \delta_{\theta,\theta_0} b(\cdot,\bar X_s,\bar\mu_s) \partial_{\theta}b(\theta,\bar X_s,\bar\mu_s)\,\mathrm{d}s- \sigma\int_0^T\partial_{\theta}b(\theta,\bar X_s,\bar\mu_s)\circ\mathrm{d}B_s^H\Bigg]\Bigg| \xrightarrow[N,n\to\infty]{}0,
\end{align*}
almost surely and
\begin{align*}
&\sup_{\theta\in\Theta}\Bigg|\partial_{\theta\theta}^2Q_n^N(\theta)- 2\mathbb{E}\Bigg[\int_0^T\left|\partial_{\theta}b(\theta,\bar X_s,\bar\mu_s)\right|^2\,\mathrm{d}s+\int_0^T\delta_{\theta,\theta_0}b(\cdot,\bar X_s,\bar\mu_s)\partial_{\theta}^2b(\theta,\bar X_s,\bar\mu_s)\,\mathrm{d}s\\
&\hspace{2.8cm} -\sigma\int_0^T\partial_{\theta}^2b(\theta,\bar X_s,\bar\mu_s)\circ\mathrm{d}B_s^H\Bigg]\Bigg|\xrightarrow[N,n\to\infty]{}0 \qquad\text{a.s.}.
\end{align*}

As a consequence,
\[
\int_0^1\partial_{\theta\theta}^2(Q_n^N)^2\bigl(s\theta_0+(1-s)\theta_n^N\bigr)\,\mathrm{d}s\xrightarrow[N,n\to\infty]{}8\sigma^2\left(\mathbb{E}\left[\int_0^T \partial_{\theta}b(\theta_0,\bar X_s,\bar\mu_s) \circ\mathrm{d}B_s^H\right]\right)^2
\]
almost surely.
\end{proposition}

Combining the previous results yields the central limit theorem for the estimator.

\begin{theorem}\label{CLT}
For $H>\frac12$, under Assumptions \ref{Bdd moment},~\ref{Lip cond},~\ref{Identif cond},~\ref{ass invertibility},~\ref{CV rate step size 1},~\ref{CV rate step size 2}, and \ref{assumptions derivative} applied to $b$, we have
\[
\sqrt{N\gamma_n^2}\, (\theta_n^N-\theta_0) \xrightarrow[N,n\to\infty]{\mathcal{D}} \mathcal{N}(0,\Sigma_H^2),
\]
where
\[
\gamma_n^2=
\begin{cases}
n^{4H-3},& \frac12<H<\frac34,\\[3mm]
\log(n)^{-1},& H=\frac34,\\[3mm]
1,& H>\frac34,
\end{cases}
\]
and
\[
\Sigma_H^2=\frac{c_H^2}{4}\left(\mathbb{E}\left[\int_0^T\partial_{\theta}b(\theta_0,\bar X_s,\bar\mu_s)\circ\mathrm{d}B_s^H\right]\right)^{-2},
\]
with $c_H$ given in Proposition~\ref{L2_Estimate_fBm}.

For $H\in(\frac14,\frac12)$, under Assumptions~\ref{Bdd moment},~\ref{Lip cond},~\ref{Identif cond},~\ref{ass invertibility},~\ref{CV rate step size 1},~\ref{CV rate step size 3} (resp.~\ref{CV rate step size 4}),~\ref{assumptions derivative} applied to $b$, and additionally~\ref{strong assumption on derivative} \textup{(}resp.~\ref{very strong assumption on derivative}\textup{)} on $b$ for $H\in(\frac13,\frac12)$ \textup{(}resp.\ $H\in(\frac14,\frac13]$\textup{)} and $\partial_{\theta}b,\partial^2_{\theta \theta}b$ satisfying the hypotheses of Lemma~\ref{lemma Strat integral}, we have
\[
\sqrt{Nn^{4H-3}}\,(\theta_n^N-\theta_0)\xrightarrow[N,n\to\infty]{\mathcal{D}}\mathcal{N}(0,\Sigma_H^2).
\]
\end{theorem}

The rate condition imposed on the discretization step in the Young regime coincides with the corresponding condition in the Brownian setting; see, for example,~\cite{amorino2023parameter}.
The convergence rates we obtain for our CLTs are novel; indeed, very few Central Limit Theorems have been established in the fractional case. Furthermore, for $H \leq \frac{3}{4}$, the rate of convergence depends on the discretization step. This dependence arises because the second-order moments change behavior at this threshold.

\subsection{Joint estimation of the Hurst parameter and diffusion coefficient}
\label{section joint estimation}

In practical applications, neither the Hurst parameter $H$ nor the diffusion coefficient $\sigma$ should be assumed to be known. Both parameters can be estimated directly from the discretely observed trajectories, and our estimation procedure remains valid when their estimators are substituted into the contrast function. Our approach is strongly inspired by the work of Coeurjolly~\cite{estimationparameterfBm}.

    The estimators are based on empirical quadratic variations, obtained from weighted increments of the observed trajectories. We consider the filter $a=(1,-3,3,-1)$, which is of order $3$ in the sense of~\cite{estimationparameterfBm}. The motivation for this choice is twofold: First, any filter of order 2 or higher is sufficient to ensure that the convergence rate of an estimator based on quadratic variations does not worsen with $H$. For $H>\frac{1}{2}$ we will show, moreover, that for any such filter we can construct a plug-in estimator for the drift coefficient having the same asymptotic properties as the initial estimator $\hat\theta_n^N$. This is no longer true for $H<\frac{1}{2}$ under our assumptions, however, an antisymmetric filter, i.e. a filter with $a_q=-a_{3-q}$, induces stochastic cancellations that allow us to handle the rough case as well. This provides the second piece of motivation for our filter $a$. For this filter, we define the corresponding filtered time series by
\[
V^aX^{i,N}\left(\frac{kT}{n}\right)=\sum_{j=0}^3a_jX^{i,N}_{\frac{(k-j)T}{n}},\qquad k=3,\ldots,n-1.
\]

To estimate $H$ and $\sigma$ jointly, we introduce the dilated filters $a^m$, $m=1,2$, defined by
\[
a_i^m=
\begin{cases}
a_j,& \text{if }i=jm,\\
0,& \text{otherwise},
\end{cases}
\qquad i=0,\ldots,3m.
\]

We then define the empirical second moment of the corresponding discrete variations by
\[
S_nX^{i,N}(a^m)= \frac{1}{n-3m} \sum_{k=3m}^{n-1} \left| V^{a^m}X^{i,N} \left(\frac{kT}{n}\right)\right|^2.
\]

The estimators of $H$ and $\sigma$ are first defined separately for each trajectory. For a trajectory $X^{i,N}$, we denote them by $H_n^i$ and $\sigma_n^i$, respectively. They are defined as the estimators $\tilde H^{\mathrm{ols}}$ and $\tilde C^{\mathrm{ols}}$ introduced in Section 3.2 of~\cite{estimationparameterfBm}, with the interval $[0,1]$ replaced by $[0,T]$, the fractional Brownian motion replaced by the observed trajectory $X^{i,N}$, and the filter $a=(1,-3,3,-1)$ together with $m=1,2$. The results in~\cite{estimationparameterfBm} are based on the asymptotic behavior of the functionals $S_nB^{i,H}$. In the present setting, it is therefore sufficient to show that the corresponding functionals of the interacting trajectories differ negligibly from those of the driving fractional Brownian motions. Since the filters $a^m$ are of order three, the first-order
term vanishes in the corresponding Taylor expansion.

We can obtain similar results as Proposition 4 of~\cite{estimationparameterfBm} for the estimators based on each particle.   

\begin{proposition}\label{estimator one particle sigma and H}
Let $a$ be a filter of order $p\geq1$, and let $H_n^i$ and $\sigma_n^i$ be the estimators of $H$ and $\sigma$ defined from the corresponding filter, as described above. Under Assumptions~\ref{Bdd moment} and~\ref{Lip cond}, for every $H\in(0,1)$,
\begin{enumerate}[label=\roman*)]
    \item for every $i\geq1$,
    \[
    (H_n^i,\sigma_n^i) \xrightarrow[n\to\infty]{}(H,\sigma)  \qquad\text{almost surely};
    \]
    \item if $p>\frac14+H$, then for every $i\geq1$,
    \[
    \sqrt{n}\,(H_n^i-H) \xrightarrow[n\to\infty]{\mathcal{D}} \mathcal{N}(0,\Sigma^2),
    \]
    and
    \[
    \frac{\sqrt{n}}{\log(n)} \frac{\sigma_n^i-\sigma}{\sigma}\xrightarrow[n\to\infty]{\mathcal{D}}\mathcal{N}(0,\Sigma^2).
    \]
\end{enumerate}
Here $\Sigma$ is the asymptotic variance specified in Proposition 4 (with $k=2$) of~\cite{estimationparameterfBm}.
\end{proposition}

The proof is given in Section~\ref{section proofs Unknown H and sigma}. The previous proposition concerns a single trajectory. Since in our setting we observe $N$ interacting trajectories, we can exploit the large number of particles by averaging the individual estimators. We denote these averages by $H_n^N$ and $\sigma_n^N$, respectively. Although the estimators are not independent because of the interaction, the leading term in their asymptotic expansions is determined by the fractional Brownian noises. Since these noises are independent across particles, this leads to the following result.

\begin{theorem}\label{average estimator of sigma and H}
Under Assumptions~\ref{Bdd moment},~\ref{Lip cond}, and~\ref{CV rate step size 1}, and additionally
Assumption~\ref{CV rate step size 2} for $H\in(1/2,1)$, Assumption~\ref{CV rate step size 3} and Assumption~\ref{strong assumption on derivative}
for $H\in(1/3,1/2)$, and Assumption~\ref{CV rate step size 4} and Assumption~\ref{very strong assumption on derivative} for
$H\in(1/4,1/3]$, we have
\begin{enumerate}[label=\roman*)]
    \item $\sqrt{nN}\,(H_n^N-H)\xrightarrow[N,n\to\infty]{\mathcal{D}}\mathcal{N}(0,\Sigma^2);$
    \item $\frac{\sqrt{nN}}{\log(n)}\frac{\sigma_n^N-\sigma}{\sigma}\xrightarrow[N,n\to\infty]{\mathcal{D}}\mathcal{N}(0,\Sigma^2).$
\end{enumerate}
\end{theorem}
The proof is postponed to Section~\ref{section proofs Unknown H and sigma}.\\
We can now construct an estimator of $\theta_0$ by substituting $H_n^N$ and $\sigma_n^N$ into the contrast function. For $H>\frac12$, we define
\[
\widetilde Q_n^N(\theta)=\frac{1}{N\alpha_n}\sum_{i=1}^N\sum_{k=0}^{n-1}\left(\left(X^{i,N}_{t_k t_{k+1}}-\alpha_n b(\theta,X^{i,N}_{t_k},\mu^N_{t_k})\right)^2-(\sigma_n^N)^2\alpha_n^{2H_n^N}\right),
\]
while for $H\in(\frac14,\frac12)$ we set
\[
\widetilde Q_n^N(\theta)=\frac{1}{N\alpha_n}\sum_{i=1}^N\sum_{k=0}^{n-1}\left(\left(X^{i,N}_{t_k t_{k+1}}-\frac{\alpha_n}{2}\left[b(\theta,X^{i,N}_{t_k},\mu^N_{t_k})+b(\theta,X^{i,N}_{t_{k+1}},\mu^N_{t_{k+1}})\right]\right)^2-(\sigma_n^N)^2\alpha_n^{2H_n^N}\right).
\]

The corresponding estimator is
\[
\widetilde\theta_n^N=\operatorname*{arg\,min}_{\theta\in\Theta}\left|\widetilde Q_n^N(\theta)\right|^2.
\]

We obtain the following consistency result whose proof is given in Section~\ref{section proofs Unknown H and sigma}.

\begin{theorem}\label{consistency joint estimation}
For $H>\frac12$, under the assumptions of Theorem~\ref{Consistency Young case}, and for $H\in(\frac14,\frac12)$ under the assumptions of Theorem~\ref{consistency rough case}, we have
\[
\widetilde\theta_n^N \xrightarrow[N,n\to\infty]{} \theta_0
\]
almost surely.
\end{theorem}

The asymptotic behavior of the estimator changes according to the value of $H$. Indeed, the normalization in the central limit theorem for the drift parameter depends on $H$, whereas the estimators of $H$ and $\sigma$ have different rates. Consequently, the joint asymptotic behavior is determined by the slower of the relevant rates.

\begin{theorem}\label{CLT joint estimation}
Under the assumptions of Theorem~\ref{CLT}, we have:
\begin{align*}
    \sqrt{N\gamma_n^2}\, (\widetilde\theta_n^N-\theta_0) \xrightarrow[N,n\to\infty]{\mathcal{D}}\mathcal{N}(0,\widetilde\Sigma_H).
\end{align*}
Here,
\begin{center}
$
\widetilde\Sigma_H:=
\begin{cases}
\displaystyle
\frac{\Lambda_H^2}{4}\left(\mathbb{E}\left[\int_0^T\partial_{\theta}b(\theta_0,\bar X_s,\bar\mu_s) \circ\mathrm{d}B_s^H\right]\right)^{-2},H<\frac34,\\[5mm]
\displaystyle 
\frac{c_H^2}{4}\left(\mathbb{E}\left[\int_0^T\partial_{\theta}b(\theta_0,\bar X_s,\bar\mu_s)\circ\mathrm{d}B_s^H\right]\right)^{-2},H\geq\frac34,
\end{cases}
$
\end{center}
where $c_H^2$ is given in Proposition~\ref{L2_Estimate_fBm}, and $\Lambda_H^2$ is given in the proof.
\end{theorem}

Substituting the estimators for $\sigma$ and $H$ does not degrade the overall rate of convergence. This phenomenon may seem surprising, given that the estimation of both the diffusion coefficient and Hurst parameter typically exhibits a slower convergence rate. However, this slower rate does not propagate to the final result due to structural cancellations that occur when evaluating the difference between $\tilde Q_n^N$ and $Q_n^N$.

\section{Numerical results}\label{section numerics}

In this section, we illustrate the performance of our estimators through numerical experiments. Throughout this section, we consider the time horizon $T=1$. We first simulate the fractional Brownian motion paths using the Davies--Harte method~\cite{davies1987tests} and then generate the particle trajectories using an Euler scheme. For each set of parameters, we perform $100$ Monte Carlo repetitions and report the resulting bias and root mean squared error (RMSE).

We first consider the linear model
\begin{equation}\label{linear model}
\mathrm{d}X_t^{i,N}=-\theta_0\left(X_t^{i,N}-\frac{1}{N}\sum_{j=1}^N X_t^{j,N}\right)\mathrm{d}t+\mathrm{d}B_t^{i,H}, \quad X_0^i=0.
\end{equation}
This model satisfies all of our required assumptions. In particular, the identifiability condition detailed in Assumption~\ref{Identif cond} can be readily verified through direct computation. We consider the mean-field equation and write
\[
b(\theta,x,\mu)=-\theta\left(x-\int y\,\mu(dy)\right).
\]
Let $m_0:=\mathbb E[\bar X_0]$. Taking expectations in the
McKean-Vlasov equation shows that
$\mathbb E[\bar X_t]=m_0$ for every $t\geq0$. Setting
$Y_t:=\bar X_t-m_0$, we obtain
\[
\mathrm{d}Y_t=-\theta_0Y_t\,\mathrm{d}t+\,\mathrm{d}B_t^H.
\]
In particular,
\[
\delta_{\theta,\theta_0}b(\cdot,\bar X_t,\bar\mu_t)
=-(\theta-\theta_0)Y_t.
\]
The left-hand side of the identifiability condition in
Assumption~\ref{Identif cond} can therefore be written as
\[
F(\theta)=
(\theta-\theta_0)^2
\mathbb E\left[\int_0^T Y_s^2\,ds\right]
+ 2(\theta-\theta_0)
\mathbb E\left[\int_0^T Y_s\circ \mathrm{d}B_s^H\right].
\]
Using the fact that
\[
\mathrm{d}B_t^H=\mathrm{d}Y_t+\theta_0 Y_t\mathrm{d}t,
\]
we obtain: 
\[
\E[\int_0^T Y_s\circ \mathrm{d}B_s^H]=\E[\int_0^T Y_s\circ \mathrm{d}Y_s+\theta_0\int_0^T Y_s^2\mathrm{d}s]=\E[\frac12Y_T^2+\theta_0\int_0^TY_s^2\mathrm{d}s]
\]
where we used the chain rule for the last equality since we have either a Young integral when $H>\frac{1}{2}$ or a geometric rough integration when $H\in (\frac14,\frac12)$.  Thus, we have:
\[
F(\theta)=(\theta-\theta_0)\left((\theta+\theta_0)\E[\int_0^TY_s^2\mathrm{d}s]+\E[Y_T^2]\right).
\]
Then the roots of $F$ are 
\[
\theta=\theta_0 \qquad and \quad\theta=-\theta_0-\frac{\E[Y_T^2]}{\E[\int_0^T Y_s^2ds]}<0 .
\]
Hence, for any parameter space such that $\Theta\subset\R_+$, the identifiability condition is verified.

\begin{table}[H]
    \centering
    \begin{tabular}{cccc}
        \toprule
        & $H=0.4$ & $H=0.6$ & $H=0.8$ \\
        \midrule
        & $N=30,n=1000$ & $N=30,n=1000$ & $N=30,n=1000$ \\
        \midrule
        $\theta_0=4$
        & $0.34\,(0.002)$
        & $0.119\,(0.021)$
        & $0.045\,(0.041)$ \\
        $\theta_0=10$
        & $0.003\,(-0.002)$
        & $0.003\,(-0.003)$
        & $0.003\,(-0.003)$ \\
        \bottomrule
    \end{tabular}
    \caption{Simulation results for the model~\eqref{linear model}, reported in the format RMSE (Bias).}
\end{table}
We take $\theta_0=4$ and $\theta_0=10$, and parameter space $\Theta=[3,13]$. The numerical results are consistent with the qualitative behavior predicted by the theory: as \(H\) increases, the RMSE decreases, while the bias remains relatively small in the examples considered. The results for $\theta_0=4$ exhibit a larger bias than those for $\theta_0=10$, which may be due to its proximity to the boundary of the parameter space. We also observe that a relatively small number of particles, $N=30$, already leads to satisfactory numerical performance. 

\vspace{0.3 cm}

We next consider a nonlinear model arising in opinion dynamics:
\begin{equation}
\mathrm{d}X_t^{i,N}=-\frac{1}{N}\sum_{j=1}^N\varphi_{\theta_0}\left(|X_t^{i,N}-X_t^{j,N}|\right)\left(X_t^{i,N}-X_t^{j,N}\right)\mathrm{d}t+\,\mathrm{d}B_t^{i,H}, \quad X_0^i=0,
\label{non linear model}
\end{equation}
where 
\[
\varphi_\theta(x)=\exp\left(-\frac{0.01}{1-(x-\theta)^2}\right)\mathbf{1}_{[\theta-1,\theta+1]}(x).
\]
The interaction kernel $\varphi_\theta$ is a smooth, compactly
supported function whose support is centered at $\theta$. It is, in particular, a smooth approximation of the indicator function $\mathbf{1}_{[\theta-1,\theta+1]}$.
Since its argument is the nonnegative distance between particles,
the effective interaction range is
$[\theta-1,\theta+1]\cap[0,\infty)$.
Thus, the parameter $\theta$ determines the range of distances
over which the interaction is active. Because the interaction term is attractive, particles tend to
move towards one another whenever their distance lies in the
effective interaction range. In the opinion-dynamics
interpretation, this mechanism promotes the formation of
clusters of agents with similar opinions. For further details on this opinion dynamics model, we refer to \cite{sharrock2021parameter} and the references therein. Furthermore, the drift satisfies all the regularity assumptions since it is linear in the measure and the function $\varphi_\theta$ is $C^\infty$.  

For this model, we obtain the following results.

\begin{table}[H]
    \centering
    \begin{tabular}{cccc}
        \toprule
        & $H=0.4$ & $H=0.6$ & $H=0.8$ \\
        \midrule
        & $N=30,n=1000$ & $N=30,n=1000$ & $N=30,n=1000$ \\
        \midrule
        $\theta_0=-0.8$
        & $0.41\,(0.003)$
        & $0.226\,(0.012)$
        & $0.052\,(-0.021)$ \\
        $\theta_0=0.5$
        & $0.003\,(-0.002)$
        & $0.005\,(-0.003)$
        & $0.003\,(0.003)$ \\
        \bottomrule
    \end{tabular}
    \caption{Simulation results for the model~\eqref{non linear model},
    reported in the format RMSE (Bias).}
\end{table}
We take $\theta_0=-0.8$ and $\theta_0=0.5$, and parameter space $\Theta=[-1,1]$.
The numerical behavior of the estimator in this nonlinear model is similar to that observed in the linear setting. For $\theta_0=-0.8$, the RMSE decreases as $H$ increases.
For $\theta_0=0.5$, it remains small across the three values
of the Hurst parameter considered.

\section{Proof of the main results}\label{section proofs}

\subsection{Proof of propagation of chaos}
\label{section proofs prop chaos}

\begin{proof}[Proof of Theorem~\ref{a.s limit of empirical measure}]
We use the concentration inequality from Theorem 2 of \cite{CVmeasure}. More precisely, for every $t\in[0,T]$, $x>0$, $q>4$, and $\varepsilon\in(0,q)$,
\begin{equation}\label{concentration inequality emprirical measure}
\P\left(W_2(\bar\mu_t^N,\bar\mu_t)>x\right)\leq e^{-cNx^2}\mathbf{1}_{\{x\leq1\}}+N(Nx)^{-\frac{q-\varepsilon}{2}}.
\end{equation}

We first establish a Hölder estimate in time for the Wasserstein distance. Using the $L^p$ estimates and the Hölder regularity of the processes $\bar X^i$, we obtain
\begin{align}
&\E\left[\left|W_2(\bar\mu_t^N,\bar\mu_t)-W_2(\bar\mu_s^N,\bar\mu_s)\right|^p\right]\nonumber\\
&\qquad\leq c\,\mathbb{E}\left[W_2(\bar\mu_t^N,\bar\mu_t)-W_2(\bar\mu_t,\bar\mu_s^N)\right]^p+c\,\mathbb{E}\left[\left|W_2(\bar\mu_t,\bar\mu_s^N)-W_2(\bar\mu_s^N,\bar\mu_s)\right|^p\right]\nonumber\\
&\qquad\leq c|t-s|^{pH}.\label{regularity distance muN mu barre}
\end{align}

We can therefore obtain a bound which is uniform in time. We consider the mesh $t_k^N=kT/N^{\frac{1}{2H}}$, $k=0,\ldots,\lfloor N^{\frac{1}{2H}}\rfloor $, and set $\Delta_N:=\frac{T}{N^{\frac{1}{2H}}}.$
Then, for $q>4$, $r\geq1$, and $\varepsilon\in(0,q)$,
\begin{align*}
&\P\left(\sup_{t\in[0,T]}W_2(\bar\mu_t^N,\bar\mu_t)\geq2x \right)\\
&\quad\quad\leq \P\left(\max_{k=0,\ldots,N} W_2(\bar\mu_{t_k^N}^N,\bar\mu_{t_k^N})\geq x \right)+\P\left(\sup_{\substack{s,t\in[0,T]\\|t-s|\leq\Delta_N}}\left|W_2(\bar\mu_t^N,\bar\mu_t)-W_2(\bar\mu_s^N,\bar\mu_s)\right|\geq x\right)\\
&\quad\quad\leq N^{\frac{1}{2H}}e^{-cNx^2}+N^{1+\frac{1}{2H}}(Nx)^{-\frac{q-\varepsilon}{2}}+c\Delta_N^{rH}x^{-r}.
\end{align*}
Here we used~\eqref{concentration inequality emprirical measure} and~\eqref{regularity distance muN mu barre}. The almost sure convergence then follows from the Borel--Cantelli lemma. To obtain the stated rate, define $\eta_\gamma:=\sup_{N\in\N}N^\gamma\|W_2(\bar\mu^N,\bar\mu)\|_{\infty,[0,T]}.$
It remains to show that $\eta_\gamma$ has finite moments of every order. For $p\geq1$,
\begin{align*}
\E[|\eta_\gamma|^p]&\leq c\left(1+\E\left[|\eta_\gamma|^p\mathbf{1}_{\{\eta_\gamma\geq1\}}\right]\right)\\
&\leq c\left(1+\sum_{N\geq1}\E\left[N^{\gamma p}\|W_2(\bar\mu^N,\bar\mu)\|_\infty^p\mathbf{1}_{\{N^\gamma\|W_2(\bar\mu^N,\bar\mu)\|_\infty\geq1\}}\right]\right)\\
&\leq c\left(1+\sum_{N\geq1}\int_1^\infty px^{p-1}\P\left(N^\gamma\|W_2(\bar\mu^N,\bar\mu)\|_\infty\geq x\right)\mathrm{d}x\right)\\
&\leq c\Bigg(1+\sum_{N\geq1}\int_1^\infty px^{p-1}\Big[N^{\frac{1}{2H}}e^{-cN^{1-2\gamma}x^2}+N^{1+\frac{1}{2H}}(N^{1-\gamma}x)^{-\frac{q-\varepsilon}{2}}+N^{r\gamma}\frac{\Delta_N^{Hr}}{x^r}\Big]\mathrm{d}x\Bigg).
\end{align*}
By choosing $r$ and $q$ sufficiently large, the series is finite. This proves the result.
\end{proof}

\begin{proof}[Proof of Proposition~\ref{prop of chaos}]
We use a Grönwall argument to control the distance between the interacting particle system and its McKean-Vlasov limit. Let $p\geq2$. Using the equations satisfied by $X^{i,N}$ and $\bar X^i$, together with Jensen's inequality, we obtain
\begin{align*}
\E\left[\sup_{t\in[0,T]}|X_t^{i,N}-\bar X_t^i|^p\right]&\leq cT^{p-1}\int_0^T\E\left[\sup_{s\in[0,t]}|b(\theta_0,X_s^{i,N},\mu_s^N)-b(\theta_0,\bar X_s^i,\bar\mu_s)|^p\right]\mathrm{d}t\\
&\leq cT^{p-1}\int_0^T\E\left[\sup_{s\in[0,t]}|X_s^{i,N}-\bar X_s^i|^p+\sup_{s\in[0,t]}W_2^p(\mu_s^N,\bar\mu_s)\right]\mathrm{d}t,
\end{align*}
where the last inequality follows from the Lipschitz continuity of $b$. We next control the Wasserstein distance. By the triangle inequality,
\[\E\left[\sup_{s\in[0,t]}W_2^p(\mu_s^N,\bar\mu_s)\right]\leq c\left(\E\left[\sup_{s\in[0,t]}W_2^p(\mu_s^N,\bar\mu_s^N)\right]+\E\left[\sup_{s\in[0,t]}W_2^p(\bar\mu_s^N,\bar\mu_s)\right]\right).
\]
For the first term, using the definition of the Wasserstein distance and Jensen's inequality,

\begin{align*}
\E\left[\sup_{s\in[0,t]}W_2^p(\mu_s^N,\bar\mu_s^N)\right]&\leq \E\left[\sup_{s\in[0,t]}\left( \frac1N\sum_{j=1}^N|X_s^{j,N}-\bar X_s^j|^2\right)^{p/2}\right]\\
&\leq\E\left[\frac1N\sum_{j=1}^N\sup_{s\in[0,t]}|X_s^{j,N}-\bar X_s^j|^p\right] =\E\left[\sup_{s\in[0,t]}|X_s^{i,N}-\bar X_s^i|^p\right],
\end{align*}
where the last equality follows because the particles are identically distributed.\\
For the second term, Theorem~\ref{a.s limit of empirical measure} gives
\[
\E\left[\|W_2(\bar\mu^N,\bar\mu)\|_\infty^p\right]\leq\frac{\E[|\eta_\gamma|^p]}{N^{\gamma p}}.
\]
Hence, Grönwall's lemma yields
\[
\E\left[\|X^{i,N}-\bar X^i\|_\infty^p\right]\leq  c\E\left[\|W_2(\bar\mu^N,\bar\mu)\|_\infty^p\right]\leq c\frac{\E[|\eta_\gamma|^p]}{N^{\gamma p}}.
\]
The desired almost sure estimate then follows from Lemma~\ref{fast Lp implies as}. The estimate for $W_2(\mu^N,\bar\mu)$ follows directly from the previous bounds.
\end{proof}

\begin{proof}[Proof of Corollary~\ref{prop chaos Holder Norms}]
   The bound on the Hölder norms follows by direct computation, utilizing the estimates provided in Proposition~\ref{prop of chaos}. Let $\varepsilon>0$ and $\gamma\in(0,\frac 12)$:
    \begin{align*}
        \|X^{i,N}-\bar X^i\|_{H-\varepsilon,[0,T]}&=\sup_{s\neq t} \frac{|X^{i,N}_t-\bar X^i_t-(X^{i,N}_s-\bar X^i_s)|}{|t-s|^{H-\varepsilon}}\\
        &=\sup_{s\neq t} \frac{|\int_s^t b(\theta_0,X^{i,N}_r,\mu_r^N)-b(\theta_0,\bar X^{i}_r,\bar \mu_r)dr|}{|t-s|^{H-\varepsilon}}\\
        &\leq c \sup_{s\neq t}\frac{|\int_s^t|X^{i,N}_r-\bar X^i_r|+W_2(\mu^N_r,\bar\mu_r)dr|}{|t-s|^{H-\varepsilon}}\\
        &\leq cT^{1-(H-\varepsilon)}\frac{\eta_\gamma}{N^\gamma},
    \end{align*}
    where $\eta_\gamma$ is the random variable defined in Proposition~\ref{prop of chaos}. Finally, for the rate on the laws we use the coupling to get the same rate.
\end{proof}

\subsection{Proof of the Young case}\label{section proofs Young case}

We first establish a decomposition of the difference between the contrast evaluated at an arbitrary parameter $\theta$ and at the true parameter $\theta_0$.

\begin{lemma}\label{lemma diff contrast fct Young case}
Under Assumptions~\ref{Bdd moment},~\ref{Lip cond}, for every $N,n$ and $\theta\in\Theta$,
\[
Q_n^N(\theta)-Q_n^N(\theta_0)=\widehat Q_n^N(\theta)-2\sigma\widetilde Q_n^N(\theta)+R_n^N(\theta),
\]
where
\[
\sup_{\theta\in\Theta}|R_n^N(\theta)|\xrightarrow[N,n\to\infty]{}0
\]
almost surely, and
\begin{align*}
\widehat Q_n^N(\theta)&=\frac{\alpha_n}{N}\sum_{i=1}^N\sum_{k=0}^{n-1}\left|\delta_{\theta,\theta_0}b(\cdot,X_{t_k}^{i,N},\mu_{t_k}^N)\right|^2, \qquad \widetilde Q_n^N(\theta)&=\frac1N\sum_{i=1}^N\sum_{k=0}^{n-1}\delta_{\theta,\theta_0}b(\cdot,X_{t_k}^{i,N},\mu_{t_k}^N)B_{t_kt_{k+1}}^{i,H}.
\end{align*}
\end{lemma}

\begin{proof}
 The first result is just a computation using the Taylor expansion of the trajectories. Indeed, for all $i,N$ we have:
     \[
     X_{t_{k+1},t_k}^{i,N}=b(\theta_0,X^{i,N}_{t_k},\mu^N_{t_k})\alpha_n+\sigma B^{H,i}_{t_kt_{k+1}}+r^{i,N}_k,
     \]
     where $r^{i,N}_k$ is the remainder defined in the equation~\eqref{eq: remainder young}. We compute the difference:
     \begin{align*}
         & N\alpha_n(Q_n^N(\theta)-Q_n^N(\theta_0)) \\
         &=\sum_{i=1}^N\sum_{k=0}^{n-1}\left(| X^{i,N}_{t_{k+1},t_k}-\alpha_nb(\theta,X^{i,N}_{t_k},\mu^N_{t_k})|^2-\sigma\alpha^{2H}_n \right) -\sum_{i=1}^N\sum_{k=0}^{n-1}\left(|\ X^{i,N}_{t_{k+1},t_k}-\alpha_nb(\theta_0,X^{i,N}_{t_k},\mu^N_{t_k})|^2-\sigma\alpha^{2H}_n \right)\\
         &=\sum_{i=1}^N\sum_{k=0}^{n-1}\left( |b(\theta,X^{i,N}_{t_k},\mu_{t_k}^N)|^2-|b(\theta_0,X^{i,N}_{t_k},\mu_{t_k}^N)|^2\right)\alpha_n^2 -2\sum_{i=1}^N\sum_{k=0}^{n-1} X^{i,N}_{t_{k+1},t_k} \cdot \delta_{\theta,\theta_0}b(.,X^{i,N}_{t_k},\mu^N_{t_k}) \alpha_n.
        \end{align*}
        Then we use the Taylor expansion of $X^{i,N}$, so that the object above equals 
        \begin{align*}
  &\sum_{i=1}^N\sum_{k=0}^{n-1}\left( |b(\theta,X^{i,N}_{t_k},\mu_{t_k}^N)|^2+|b(\theta_0,X^{i,N}_{t_k},\mu_{t_k}^N)|^2\right)\alpha_n^2-2\sum_{i=1}^N\sum_{k=0}^{n-1} b(\theta,X^{i,N}_{t_k},\mu_{t_k}^N) b(\theta_0,X^{i,N}_{t_k},\mu_{t_k}^N) \alpha_n^2\\
         &-2\sigma \sum_{i=1}^N\sum_{k=0}^{n-1}\delta_{\theta,\theta_0}b(\cdot,X_{t_k}^{i,N},\mu_{t_k}^N) B^{H,i}_{t_{k+1},t_k} \alpha_n\underbrace{-2\sum_{i=1}^N\sum_{k=0}^{n-1} \delta_{\theta,\theta_0}b(\cdot,X^{i,N}_{t_k},\mu^N_{t_k})r_k^{i,N} \alpha_n}_{:=N\alpha_n R_n^N(\theta)}\\
         &=N\alpha_n(\tilde Q_n^N(\theta) -2\sigma \hat Q_n^N(\theta)+R_n^N(\theta)).
         \end{align*}
     We have the desired identity; now we prove the convergence of $R_n^N$. It will be in two steps; firstly we find an estimate of the $L^p$ moments of $r_k^{i,N}$, then we plug it into the sum: 
     \begin{align}
\E[|r_k^{i,N}|^p]&=\E[\left|\int_{t_k}^{t_{k+1}}(b(\theta_0,X^{i,N}_u,\mu_u^N)-b(\theta_0,X_{t_k}^{i,N},\mu^N_{t_k}))du \right|^p] \nonumber\\
         &\leq \alpha_n^{p-1}\E[\int_{t_k}^{t_{k+1}}\left|b(\theta_0,X^{i,N}_u,\mu_u^N)-b(\theta_0,X_{t_k}^{i,N},\mu^N_{t_k})\right|^pdu]\nonumber\\
         &\leq c \alpha_n^{p-1}\int_{t_k}^{t_{k+1}} \E[|X_u^{i,N}-X^{i,N}_{t_k}|^p+W_2^p(\mu^N_u,\mu^N_{t_k})] \leq c\alpha_n^{p(1+H)}. \label{estimate remainder}
     \end{align}
     To obtain this estimate, we used Jensen's inequality, the Lipschitz assumptions on $b$, and finally the regularity of the fractional Brownian motion, which is proven in Proposition~\ref{prop solution}. With this estimate, we can now bound $R_n^N(\theta)$:
     \begin{align*}
         \E[\underset{\theta \in \Theta}{\sup} |R_n^N(\theta)|^p]&=\E[\underset{\theta \in \Theta}{\sup}\left|\frac{2}{N\alpha_n}\sum_{i=1}^N\sum_{k=0}^{n-1}\delta_{\theta,\theta_0}b(\cdot,X^{i,N}_{t_k},\mu^N_{t_k}) r_k^{i,N} \alpha_n\right|^p ]\\
         &\leq c\E[\underset{\theta \in \Theta}{\sup}\left|\frac{1}{N}\sum_{i=1}^N\sum_{k=0}^{n-1}|\theta-\theta_0|(1+|X_{t_k}^{i,N}|^K+W_2^L(\delta_0,\mu_{t_k}^N))|r_k^{i,N}|\right|^p]\\
         &\leq cn^{p-1}\sum_{k=0}^{n-1}\E[|1+|X_{t_k}^{1,N}|^K+W_2^L(\delta_0,\mu^N_{t_k})||r_k^{1,N}|^p]\\
         &\leq c n^{p-1}\sum_{k=0}^{n-1}\E[|1+|X_{t_k}^{1,N}|^K+W_2^L(\delta_0,\mu^N_{t_k})|^{2p}]^\frac{1}{2}\mathbb{E[}|r_k^{1,N}|^{2p}]^\frac{1}{2}\leq c\alpha_n^{pH}.
     \end{align*}
     To obtain the last inequality, we used the uniform boundedness of moments of $X^{i,N}$ as in Proposition~\ref{prop solution}. Thus, we can conclude by using  Lemma~\ref{fast Lp implies as}.
 \end{proof}

\begin{remark}
The preceding argument illustrates a strategy that will be used repeatedly below: first derive suitable $L^p$ estimates for the remainder terms and then deduce almost sure convergence using Lemma~\ref{fast Lp implies as}.
\end{remark}

We next prove the convergence of the contrast evaluated at the true parameter.

\begin{proof}[Proof of Proposition~\ref{CV contrast in theta 0}]\label{proof Q_n^N(theta_0) H>1/2}
    Firstly we are going to rewrite $Q_n^N(\theta_0)$ using the Taylor expansion of $X^{i,N}$:
    \begin{align*}
        Q_n^N(\theta_0)&=\frac{1}{N\alpha_n}\sum_{i=1}^N\sum_{k=0}^{n-1}|X_{t_{k+1},t_k}^{i,N}-\alpha_nb(\theta_0,X_{t_k}^{i,N},\mu^N_{t_k})|^2-\sigma^2\alpha_n^{2H}\\
        &=\frac{1}{N\alpha_n}\sum_{i=1}^N\sum_{k=0}^{n-1}|\sigma B^{i,H}_{t_{k+1},t_k}+r_k^{i,N}|^2-\sigma^2\alpha_n^{2H}\\
        &=\frac{\sigma^2}{N\alpha_n}\sum_{i=1}^N\sum_{k=0}^{n-1}|B^{i,H}_{t_{k+1},t_k|}|^2-\alpha_n^{2H}+\frac{1}{N\alpha_n}\sum_{i=1}^N\sum_{k=0}^{n-1}|r_k^{i,N}|^2+\frac{2\sigma}{N\alpha_n}\sum_{i=1}^N\sum_{k=0}^{n-1}B_{t_{k+1},t_k}^{i,H} r^{i,N}_k.
    \end{align*}
    It is quite straightforward to show that the last two terms go to zero. Indeed, we already have an estimate for $r^{i,N}_k$ and using the regularity estimates of the fractional brownian motion we have:
     \begin{align*}
         & \E[|\frac{1}{N\alpha_n}\sum_{i=1}^N\sum_{k=0}^{n-1}B_{t_{k+1},t_k}\cdot r^{i,N}_k|^p]\leq \frac{1}{\alpha_n^p}n^{p-1}\sum_{k=0}^{n-1}\E[|B_{t_{k+1},t_k}^{1,H}|^p|r_k^{1,N}|^p]\\
        &\leq \frac{1}{\alpha_n^p}n^{p-1}\sum_{k=0}^{n-1}\E[|B_{t_kt_{k+1}}|^{2p}]^{1/2}\E[|r^{i,N}_k|^{2p}]^{1/2} \leq c\frac{1}{\alpha_n^p}n^p\alpha_n^{pH}\alpha_n^{p(1+H)} \leq c\alpha_n^{p(2H-1)}
    \end{align*}
    we used the fact that $n\alpha_n=T$ which is constant. Then we compute the last term:
    \begin{align*}
    \E [\left|\frac{1}{N\alpha_n}\sum_{i=1}^N\sum_{k=0}^{n-1}|r_{k}^{i,N}|^2\right|^p]
        \leq \frac{n^{p-1}}{\alpha_n^p}\sum_{k=0}^{n-1}\E[|r_k^{i,N}|^{2p}] \leq c\frac{n^p}{\alpha_n^p}\alpha_n^{2p(1+H)}=c\alpha_n^{2pH}.
    \end{align*}
    Hence we conclude using Lemma~\ref{fast Lp implies as}. For the first term we are going to use the equivalence of moments in a Gaussian chaos as discussed in~\cite{MalliavinNualart}. And for the second moment we have the Lemma~\ref{L2_Estimate_fBm} which gives us:
   \begin{align*}
\E[\left|\frac{1}{N\alpha_n}\sum_{i=1}^N\sum_{k=0}^{n-1}|B_{t_kt_{k+1}}|^2-\alpha_n^{2H}\right|^2]\underset{N,n\rightarrow \infty}{\sim}
        \begin{cases}
    c_H\frac{1}{N} \quad \text{if} \quad H\in (3/4,1),\\
    c_H\frac{\log(n)}{N} \quad if \quad H=3/4,\\
    c_H\frac{\alpha_n^{4H-3}}{N} \quad \text{if} \quad H\in(1/2,3/4).
        \end{cases}
   \end{align*}
   Thus, it goes to zero in all cases because of Assumption~\ref{CV rate step size 1}. Therefore by the equivalence of moments we have the same bounds for large $L^p$ norms. Finally, we use Lemma~\ref{fast Lp implies as} to conclude the almost sure convergence towards zero of this term as well. 
    
\end{proof}

We next study the two terms $\widehat Q_n^N$ and $\widetilde Q_n^N$. The first one converges to the deterministic quadratic term appearing in the limiting contrast. The second one converges to a pathwise stochastic integral. Unlike in the Brownian case, this integral is generally not centered. In both cases, we first remove the time discretization, then use propagation of chaos, and finally apply a law of large numbers.

\begin{lemma}\label{CV somme Riemann classique}
Let $f,g:\Theta\times\R\times\mathcal{P}_2(\R)\to\R$ satisfy, for some $c>0$ and $K,L\in\N$,
\[
|f(\theta,x,\mu)-f(\theta,y,\nu)| \leq c\bigl(|x-y|+W_2(\mu,\nu)\bigr),
\]
\[
|f(\theta,x,\mu)-f(\tilde\theta,x,\mu)| \leq c|\theta-\tilde\theta| \left(1+|x|^K+W_2^L(\mu,\delta_0)\right),
\]
and assume that the same bounds hold for $g$. Then, for every $H\in(0,1)$,
\begin{align*}
\sup_{\theta\in\Theta}\Bigg|\frac{\alpha_n}{N}\sum_{i=1}^N\sum_{k=0}^{n-1}f(\theta,X_{t_k}^{i,N},\mu_{t_k}^N)g(\theta,X_{t_k}^{i,N},\mu_{t_k}^N)-\E\left[\int_0^Tf(\theta,\bar X_s^1,\bar\mu_s) g(\theta,\bar X_s^1,\bar\mu_s)\,\mathrm{d}s\right] \Bigg|\xrightarrow[N,n\to\infty]{}0
\end{align*}
almost surely. In particular, under Assumptions~\ref{Bdd moment} and
\ref{Lip cond},
\[
\sup_{\theta\in\Theta}\left|\widehat Q_n^N(\theta)-\E\left[ \int_0^T|\delta_{\theta,\theta_0}b(\cdot,\bar X_s,\bar\mu_s)|^2\,\mathrm{d}s\right]\right|\xrightarrow[N,n\to\infty]{}0 \qquad \mbox{almost surely.}
\]

\end{lemma}

\begin{proof}
We denote 
\[
L^{(1)}(\theta):=\E\left[\int_0^Tf(\theta,\bar X_s^1,\bar\mu_s)g(\theta,\bar X_s^1,\bar\mu_s)\,\mathrm{d}s\right]
\]
     and we decompose the limit in three parts:
     \begin{align}
         &\underset{\theta \in \Theta}{\sup}| \frac{\alpha_n}{N}\sum_{i=1}^N\sum_{k=0}^{n-1}f(\theta,X^{i,N}_{t_k},\mu^N_{t_k})g(\theta,X^{i,N}_{t_k},\mu^N_{t_k})-L^{(1)}(\theta)|\nonumber\\
         &\leq \underset{\theta \in \Theta}{\sup}|\frac{\alpha_n}{N}\sum_{i=1}^N\sum_{k=0}^{n-1}f(\theta,X^{i,N}_{t_k},\mu^N_{t_k})g(\theta,X^{i,N}_{t_k},\mu^N_{t_k})-\frac{1}{N}\sum_{i=1}^N\int_0^Tf(\theta,X^{i,N}_{s},\mu^N_{s})g(\theta,X^{i,N}_{s},\mu^N_{s})\mathrm{d}s|\nonumber\\
         &+ \underset{\theta \in \Theta}{\sup}|\frac{1}{N}\sum_{i=1}^N\int_0^Tf(\theta,X^{i,N}_{s},\mu^N_{s})g(\theta,X^{i,N}_{s},\mu^N_{s})\mathrm{d}s-\frac{1}{N}\sum_{i=1}^N\int_0^Tf(\theta,\bar X^{i}_{s},\bar \mu_{s})g(\theta,\bar X^{i}_{s},\bar \mu_{s})\mathrm{d}s|\nonumber\\
         &+\underset{\theta \in \Theta}{\sup}|\frac{1}{N}\sum_{i=1}^N\int_0^Tf(\theta,\bar X^{i}_{s},\bar \mu_{s})g(\theta,\bar X^{i}_{s},\bar \mu_{s})\mathrm{d}s-L^{(1)}(\theta)|.\label{eq decomposition Young case  riemman sums}
     \end{align}
    For the first term, we are just approximating some Riemann sums, so using regularity estimates on $f$ and $g$ we have:
    \begin{align*}
       &\E[ \sup_{\theta \in \Theta}|\frac{\alpha_n}{N}\sum_{i=1}^N\sum_{k=0}^{n-1}f(\theta,X^{i,N}_{t_k},\mu^N_{t_k})g(\theta,X^{i,N}_{t_k},\mu^N_{t_k})-\frac{1}{N}\sum_{i=1}^N\int_0^Tf(\theta,X^{i,N}_s,\mu_s^N)g(\theta,X^{i,N}_s,\mu_s^N)\mathrm{d}s|^p]\\
        &\leq \E[\sup_{\theta \in \Theta}\sum_{k=0}^{n-1}\int_{t_k}^{t_{k+1}}\left|f(\theta,X^{i,N}_{t_k},\mu^N_{t_k})g(\theta,X^{i,N}_{t_k},\mu^N_{t_k})-f(\theta,X^{i,N}_{s},\mu^N_{s})g(\theta,X^{i,N}_{s},\mu^N_{s})\right|^p\mathrm{d}s]\\
        &\leq c\E[\sup_{\theta \in \Theta} \sum_{k=0}^{n-1}\|g(\theta,X^{i,N},\mu^N)\|_\infty^p\int_{t_k}^{t_{k+1}}|f(\theta,X^{i,N}_{t_k},\mu_{t_k}^N)-f(\theta,X^{i,N}_s,\mu_s^N)|^p\mathrm{d}s]\\
        &+\E[\sup_{\theta \in \Theta} \sum_{k=0}^{n-1}\|f(\theta,X^{i,N},\mu^N)\|_\infty^p\int_{t_k}^{t_{k+1}}|g(\theta,X^{i,N}_{t_k},\mu_{t_k}^N)-g(\theta,X^{i,N}_s,\mu_s^N)|^p\mathrm{d}s]\\
        &\leq c \E[\|f(\theta,X^{1,N},\mu^N)\|_\infty^{2p}+\|g(\theta,X^{1,N},\mu^N)\|_\infty^{2p})^{\frac{1}{2}} \sum_{k=0}^{n-1}\int_{t_k}^{t_{k+1}} \E[ |X^{1,N}_s-X^{1,N}_{t_k}|^{2p}+W_2^{2p}(\mu_{t_k}^N,\mu^N_s)\mathrm{d}s]^{\frac{1}{2}}\\
        &\leq c \alpha_n^{pH}\E[(1+\|X^{1,N}\|_\infty^{2p}+\|W_2^p(\mu^N,\delta_0)\|^{2p}_\infty)]^{\frac{1}{2}},
    \end{align*}
    where we have used that the Hölder norm in $L^p$ does not depend on $N$ and that the expectation of the infinite norm is finite independently of $N$ by the Proposition~\ref{prop solution}. So, using Lemma~\ref{fast Lp implies as}, we have the almost sure convergence. 
     Now for the second term we are going to use our propagation of chaos:
     \begin{align*}
         \sup_{\theta \in \Theta} |\frac{1}{N}\sum_{i=1}^N \int_0^T &f(\theta,X^{i,N}_s,\mu^N_s)g(\theta,X^{i,N}_s,\mu^N_s)-f(\theta,\bar{X}^i_s,\bar{\mu}_s)g(\theta,\bar X^{i}_s,\bar \mu_s)\mathrm{d}s|\\&\leq \sup_{\theta \in \Theta} T\frac{1}{N}\sum_{i=1}^N \|f(\theta,X^{i,N},\mu^N)g(\theta,X^{i,N},\mu^N)-f(\theta,\bar{X}^i,\bar{\mu})g(\theta,\bar X^{i},\bar \mu)\|_\infty\\
         &\leq \sup_{\theta\in\theta}T\frac{1}{N}\sum_{i=1}^N\|g(\theta,\bar{X}^i,\bar{\mu})\|_\infty\|f(\theta,X^{i,N},\mu^N)-f(\theta,\bar{X}^i,\bar{\mu})\|_\infty\\&
         +\|f(\theta,\bar{X}^i,\bar{\mu})\|_\infty\|g(\theta,X^{i,N},\mu^N)-g(\theta,\bar{X}^i,\bar{\mu})\|_\infty\\
         &\leq \sup_{\theta\in\theta} \frac{c}{N}\sum_{i=1}^N(\|X^{i,N}-\bar{X}^i\|_\infty+\|W_2(\mu^N,\bar{\mu})\|_\infty)(\|g(\theta,\bar{X}^i,\bar{\mu})\|_\infty+\|f(\theta,\bar{X}^i,\bar{\mu})\|_\infty)\\
         &\leq c\tilde{\eta}_\gamma\frac{1}{N^\gamma}\frac{1}{N}\sum_{i=1}^N (\|\bar X^i\|_\infty+\|W_2(\delta_0,\bar \mu)\|_\infty),
     \end{align*}
    where the last sum is finite due to the law of large numbers.\\
     For the last part of equation~\eqref{eq decomposition Young case  riemman sums}, it is a bit easier since we have that the $\bar{X}^i$ are i.i.d . Then, for a fixed $\theta$, we can use the law of large numbers. To make this convergence uniform, we use Ascoli-Arzela. In order to use it, we need to check that:
     \[Y^N(\theta)=\frac{1}{N}\sum_{i=1}^N\int_0^Tf(\theta,\bar{X}_s^i,\bar{\mu}_s)g(\theta,\bar{X}_s^i,\bar{\mu}_s)\mathrm{d}s\]
     is a.s. equicontinuous. So let us take $\theta,\tilde{\theta} \in \Theta$. Then, for all $N$,
     \begin{align*}
         |Y^N(\theta)-Y^N(\tilde{\theta})|
         &\leq \frac{1}{N}\sum_{i=1}^N \int_0^T \left|f(\theta,\bar{X}^i_s,\bar{\mu}_s)g(\theta,\bar{X}_s^i,\bar{\mu}_s)-f(\tilde{\theta},\bar{X}^i_s,\bar{\mu}_s)g(\theta,\bar{X}_s^i,\bar{\mu}_s)\right|\mathrm{d}s\\
         &\leq \frac{1}{N} \sum_{i=1}^N\int_0^T |f(\theta,\bar{X}^i_s,\bar{\mu}_s)-f(\tilde{\theta},\bar{X}^i_s,\bar{\mu}_s)||g(\theta,\bar{X}^i_s,\bar{\mu}_s)|\\
         &+|g(\theta,\bar{X}^i_s,\bar{\mu}_s)-g(\tilde{\theta},\bar{X}^i_s,\bar{\mu}_s)||f(\theta,\bar{X}^i_s,\bar{\mu}_s)|\mathrm{d}s\\
         &\leq c|\theta-\tilde{\theta}|\frac{1}{N} \sum_{i=1}^N\int_0^T (1+|\bar{X}^i_s|^K+W_2^L(\delta_0,\bar{\mu}_s))\mathrm{d}s.
     \end{align*}
    We define $A:=\{\omega:\frac{1}{N}\sum_{i=1}^N\int_0^T|(1+|\bar{X}^i_s|^K+W_2^L(\delta_0,\bar{\mu}_s))\mathrm{d}s \, \text{converges}\}$ this set is of probability one by the law of large numbers and on this set this sum is bounded uniformly in N. Thus, we can apply Ascoli-Arzela to conclude.
\end{proof}

To control $\widetilde Q_n^N$, we use the following result.

\begin{lemma}\label{CV somme Riemann fBm, smooth case}
Let $H>\frac12$ and
$f:\Theta\times\R\times\mathcal{P}_2(\R)\to\R$
satisfy, for some $K, L \in \mathbb{N}$,
\[
|f(\theta,x,\mu)-f(\theta,y,\nu)|\leq c\bigl(|x-y|+W_2(\mu,\nu)\bigr),
\]
and
\[
|f(\theta,x,\mu)-f(\tilde\theta,x,\mu)| \leq c|\theta-\tilde\theta| \left( 1+|x|^K+W_2^L(\mu,\delta_0) \right).
\]
Then
\[
\sup_{\theta\in\Theta} \left| \frac1N \sum_{i=1}^N \sum_{k=0}^{n-1} f(\theta,X_{t_k}^{i,N},\mu_{t_k}^N) B_{t_kt_{k+1}}^{i,H}-\E\left[\int_0^T f(\theta,\bar X_s^i,\bar\mu_s) \circ\mathrm{d}B_s^{i,H}\right] \right|\xrightarrow[N,n\to\infty]{}0
\]
almost surely, where the stochastic integral is understood in the Young sense. Consequently, under Assumptions~\ref{Bdd moment} and~\ref{Lip cond}:
\[
\sup_{\theta\in\Theta}\left|\widetilde Q_n^N(\theta)-\E\left[\int_0^T\delta_{\theta,\theta_0}b(\cdot,\bar X_s^1,\bar\mu_s)\circ\mathrm{d}B_s^{1,H}\right]\right|\xrightarrow[N,n\to\infty]{}0
\]
almost surely.
\end{lemma}

\begin{proof}
The proof is going to follow the same approach as Lemma~\ref{CV somme Riemann classique}. So we denote
\[
L^{(2)}(\theta):=\E[\int_0^Tf(\theta,\bar{X}^i_s,\bar{\mu}_s)\circ\mathrm{d}B_s^H]
\]
and we split the limit:
    \begin{align*}
    \underset{\theta\in\Theta}{\sup}&|\frac{1}{N}\sum_{i=1}^N\sum_{k=0}^{n-1} f(\theta,X^{i,N}_{t_k},\mu^N_{t_k})B^{i,H}_{t_kt_{k+1}}-L^{(2)}(\theta)|\\
    &\leq \underset{\theta\in\Theta}{\sup} |\frac{1}{N}\sum_{i=1}^N\sum_{k=0}^{n-1} f(\theta,X^{i,N}_{t_k},\mu^N_{t_k})B^{i,H}_{t_kt_{k+1}}-\frac{1}{N}\sum_{i=1}^N \int_0^Tf(\theta,X^{i,N}_s,\mu^N_s)\circ \mathrm{d}B^{i,H}_s|\\
    &+\sup_{\theta \in \Theta}|\frac{1}{N}\sum_{i=1}^N \int_0^Tf(\theta,X^{i,N}_s,\mu^N_s)\circ \mathrm{d}B^{i,H}_s-\frac{1}{N}\sum_{i=1}^N\int_0^Tf(\theta,\bar{X}^{i}_s,\bar{\mu}_s^N) \circ \mathrm{d}B^{H,i}_s|\\
    &+\underset{\theta\in\Theta}{\sup}|\frac{1}{N}\sum_{i=1}^N\int_0^Tf(\theta,\bar{X}^i_s,\bar{\mu}_s) \circ \mathrm{d}B^{H,i}_s-L^{(2)}(\theta)|\ .
    \end{align*}
     We are going to prove that each part goes to 0. \\
     Regarding the first term, we use the same strategy as in the proof of Lemma~\ref{CV somme Riemann classique}. We start with the Riemann sums, we take $\varepsilon>0$ such that $H-\varepsilon>\frac{1}{2}$ and note that:
     \begin{align*}
         &\E[\sup_{\theta \in \Theta}|\frac{1}{N}\sum_{i=1}^N\sum_{k=0}^{n-1}\int_{t_k}^{t_{k+1}} f(\theta,X^{i,N}_{t_k},\mu^N_{t_k})-f(\theta,X^{i,N}_s,\mu^N_s)\circ \mathrm{d}B^{i,H}_s|^p \\
         &\leq c\E[\sup_{\theta \in \Theta}(\sum_{k=0}^{n-1}\|f(\theta,X^{i,N},\mu^N)\|_{H-\varepsilon,[t_kt_{k+1}]}\|B^{i,H}\|_{H-\varepsilon,[t_kt_{k+1}]}\alpha_n^{2(H-\varepsilon)})^p] \leq c \alpha_n^{p(2(H-\varepsilon)-1)},
     \end{align*}
     where we still used the fact that $f$ is Lipschitz uniformly in $\theta$ and that the Hölder norm of $X^{i,N}$ is uniform in $N$ in the $L^p$ moments. Finally, Lemma~\ref{fast Lp implies as} gives us the almost sure convergence. For the second term we use the propagation of chaos for the Hölder norm:
     \begin{align*}
         &\sup_{\theta\in\Theta}|\frac{1}{N}\sum_{i=1}^N\int_0^T f(\theta,X^{i,N}_s,\mu_s^N)-f(\theta,\bar{X}^i_s,\bar{\mu}_s)\mathrm{d}B^{i,H}_s|\\&\leq c\frac{1}{N} \sum_{i=1}^N \|f(\theta,X^{i,N},\mu^N)-f(\theta,\bar{X}^i,\bar{\mu})\|_{H-\varepsilon,[0,T]}\|B^{i,H}\|_{H-\varepsilon,[0,T]}\\
         &\leq c\frac{1}{N} \sum_{i=1}^N \|X^{i,N}-\bar{X}^i\|_{H-\varepsilon,[0,T]}+\|W_2(\mu^N,\bar{\mu})\|_{H-\varepsilon,[0,T]} \leq c \tilde{\eta}_\gamma \frac{1}{N^\gamma}
     \end{align*}
    for some $\frac{1}{2}>\gamma >0$.\\
    For the final step, we apply the Arzelà–Ascoli theorem to establish uniform convergence in $\theta$. Having already obtained pointwise convergence via the law of large numbers, it remains to verify the equicontinuity of 
\[
    Y_n^N(\theta) = \frac{1}{N}\sum_{i=1}^N \int_0^T f(\theta, \bar{X}^i_s, \bar{\mu}_s) \circ \mathrm{d}B^{H,i}_s.
\]
    To do it, we are going to use the Lipschitz continuity in $\theta$ of $b$ and the independence of the particles in the limit. Let $\theta,\tilde{\theta}\in\Theta$:
    \begin{align*}
       & | \frac{1}{N}\sum_{i=1}^N\int_0^Tf(\theta,\bar{X}^i_s,\bar{\mu}_s)\circ \mathrm{d}B^{H,i}_s-\frac{1}{N}\sum_{i=1}^N\int_0^Tf(\tilde\theta,\bar{X}^i_s,\bar{\mu}_s)\circ \mathrm{d}B^{H,i}_s |\\
        \qquad &\leq \frac{1}{N}\sum_{i=1}^N|\int_0^T\delta_{\theta,\tilde{\theta}}f(\cdot,\bar{X}^i_s,\bar{\mu}_s)\circ \mathrm{d}B^{H,i}_s|\\
        &\leq \frac{1}{N}\sum_{i=1}^N|\int_0^T\delta_{\theta,\tilde{\theta}}f(\cdot,\bar{X}^i_s,\bar{\mu}_s)-\delta_{\theta,\tilde\theta}f(\cdot,X_0^i,\bar\mu_0)\circ \mathrm{d}B^{H,i}_s|+\frac{1}{N}\sum_{i=1}^N|\delta_{\theta,\tilde\theta}f(\cdot,X_0^i,\bar\mu_0)B_T^{i,H}|\\
        &\leq \frac{c}{N}\sum_{i=1}^N\|f(\theta,X^i,\bar\mu)-f(\tilde \theta,X^i,\bar\mu)\|_{H-\varepsilon,[0,T]}\|B^{i,H}\|_{H-\varepsilon,[0,T]} +\frac{1}{N}\sum_{i=1}^N|\delta_{\theta,\tilde\theta}f(\cdot,X_0^i,\bar\mu_0)B_T^{i,H}|\\
        &\leq c|\theta-\tilde{\theta}| \frac{1}{N}\sum_{i=1}^N [(\|1+|\bar{X}^i|^K+W_2^L(\delta_0,\bar{\mu})\|_{H-\varepsilon,[0,T]})\|B^{H,i}\|_{H-\varepsilon,[0,T]} +(1+|X_0^i|^L+W_2^K(\delta_0,\bar\mu_0))|B_T^{i,H}|].
    \end{align*}
    By the law of large numbers, both sums are finite almost surely. Thus, almost surely $Y_n^N(\theta)$ is equicontinuous, which concludes the proof.
\end{proof}

\subsection{Proof for the rough case}\label{section proofs rough case}

We begin with the proof of the main technical tool used in the rough regime, namely Lemma~\ref{lemma Strat integral}.

\begin{proof}[Proof of Lemma~\ref{lemma Strat integral}]
The proof is inspired by Theorem 3.5 of Perkowski and Matsuda \cite{stochasticsewingupgrade}. More  precisely, we use the extended version of the stochastic sewing lemma stated in Theorem 5.5 of~\cite{stochasticsewingupgrade}.

For a two-parameter function $A:[0,T]^2_{\leq}\to\R$, we use the standard notation
\[
\delta A_{sut}=A_{st}-A_{su}-A_{ut}, \qquad 0\leq s\leq u\leq t\leq T.
\]
When $A_{st}=(f_t+f_s)g_{st},$
we have
\begin{align}
\delta A_{sut}=(f_t-f_u)(g_u-g_s)-(f_u-f_s)(g_t-g_u).
\label{delta of average}
\end{align}
We consider the germ $A_{st}=\left(f(\theta,\bar X_t^i,\bar\mu_t)+f(\theta,\bar X_s^i,\bar\mu_s)\right)B_{st}^{i,H}.$
To apply the stochastic sewing result, we need to establish estimates for all $0\leq v<s<u<t\leq T$ of the form
\begin{equation}\label{sto sewing cond Lp}
\|\delta A_{sut}\|_{L^p(\Omega)}\lesssim|t-s|^{\beta_3},
\end{equation}
and, whenever $t-s\leq s-v$,
\begin{equation}\label{sto sewing cond Lp conditional}
\left\|\E[\delta A_{sut}\mid\mathcal F_v]\right\|_{L^p(\Omega)}\lesssim |s-v|^{-\alpha}|t-s|^{\beta_1}+|t-v|^\gamma|t-s|^{\beta_2},
\end{equation}
where the nonnegative exponents satisfy
\begin{align}
\min\{\beta_1,2\beta_3\}>1,\qquad \gamma+\beta_2>1, \qquad 1+\alpha-\beta_1 <\alpha\left(\frac{\gamma+\beta_2-1}{\gamma}\right),
\label{cond sto sewing}
\end{align}
although Theorem 5.5 in Matsuda and Perkowski~\cite{stochasticsewingupgrade} is not stated in our general setting, their Remark 5.6 confirms that the result applies to our framework. Note that they initially impose the condition $1+\alpha-\beta_1<\alpha(2\beta_3-1)$ (Assumption 1.12) however, they demonstrate in Appendix A that this assumption can be dropped. 
We first estimate $\delta A$ in $L^p$. By
\eqref{delta of average},
\[
\delta A_{sut}=f(\theta,\bar X,\bar\mu)_{ut}B_{su}^H-f(\theta,\bar X,\bar\mu)_{su}B_{ut}^H.
\]
Using the Lipschitz continuity of $f$ and the Hölder regularity of $\bar X$ and $\bar\mu$, we obtain
\[
\|\delta A_{sut}\|_{L^p(\Omega)}\lesssim|t-s|^{2H},
\]
which means that~\eqref{sto sewing cond Lp} is satisfied with $\beta_3=2H$. The second part of the proof relies on estimating the conditional expectation $\E[\delta A_{sut}|\mathcal{F}_v]$. We take $t-s\leq s-v$, and we decompose $X$ into its drift and noise components. Since
the drift component is sufficiently regular, we treat its
contribution separately and focus on the noise component. We denote for $v\leq u$
\[
\psi_u^v:=X_0+\int_0^vb(\theta_0,\bar{X}_r,\bar{\mu}_r)\mathrm{d}r+(u-v)b(\theta_0,\bar X_v,\mu_v).
\]
 This means we can write:
\begin{align*}
    A_{st}&=(f(\theta,\bar{X}_t,\bar{\mu}_t)+f(\theta,\bar{X}_s,\bar{\mu}_s))B_{st}^H\\
    &=\underbrace{(f(\theta, \bar X_t,\bar{\mu}_t)+f(\theta, \bar X_s,\bar{\mu}_s)-f(\theta, \psi_t^v+\sigma B_t^H,\bar{\mu}_t)-f(\theta, \psi_s^v+\sigma B_s^H,\bar{\mu}_s))B_{st}^H}_{:=R_{st}}\\
    &+\underbrace{(f(\theta, \psi_t^v+\sigma B_t^H,\bar{\mu}_t)+f(\theta, \psi_s^v+\sigma B_s^H,\bar{\mu}_t))B_{st}^H.}_{:=\hat{A}_{st}}
\end{align*}
Using the fact that the $\delta$ operator is linear, we start by looking at $\delta R$ and using~\eqref{delta of average} we have:
\begin{align}
    \|\E[\delta R_{sut}|\mathcal{F}_v]\|_{L^p}&\leq\| \delta R_{sut}\|_{L^p}\nonumber\\
    &\leq \| (f(\theta,\bar X_t,\bar{\mu}_t)-f(\theta,\psi_t^v+\sigma B_t^H,\bar{\mu}_t)-\nonumber\\
    &\quad \quad(f(\theta,\bar X_u,\bar{\mu}_u)-f(\theta,\psi_u^v+\sigma B_u^H,\bar{\mu}_u))B_{su}^H+(f(\theta,\bar X_u,\bar{\mu}_u)\nonumber\\
    &\qquad-f(\theta,\psi_u^v+\sigma B_u^H,\bar{\mu}_u)-(f(\theta,\bar X_s,\bar{\mu}_s)-f(\theta,\psi_s^v+\sigma B_s^H,\bar{\mu}_s))B_{ut}^H\|_{L^p}\nonumber\\
    &\leq c  \left\|\left( \int_v^tb(\theta_0,\bar X_r,\bar \mu_r)-b(\theta_0,\bar X_v,\bar\mu_v)\mathrm{d}r+\int_v^ub(\theta_0,\bar X_r,\bar \mu_r)-b(\theta_0,\bar X_v,\bar\mu_v)\mathrm{d}r\right)B_{su}^H\right\|_{L^p}\nonumber\\
    &\quad+\left\|\left( \int_v^sb(\theta_0,\bar X_r,\bar \mu_r)-b(\theta_0,\bar X_v,\bar\mu_v)\mathrm{d}r+\int_v^ub(\theta_0,\bar X_r,\bar \mu_r)-b(\theta_0,\bar X_v,\bar\mu_v)\mathrm{d}r\right)B_{ut}^H\right\|_{L^p}\nonumber\\
    &\leq c |t-v|^{1+H}|t-s|^H. \label{ineq: remainder sto sewing}
\end{align}
For the last two lines, we use first the Lipschitz assumption on $f$ then the regularity estimates on the drift and $B^H$. 

We now estimate $\hat A$ using the Mandelbrot--Van Ness
representation of fractional Brownian motion~\cite{MandelbrotVannessRepresentation}, i.e $B_s^{H}=\int_{-\infty}^sK(s,r)dW_r$ where the $W$ is a Brownian motion. We cut the integral in two parts:
\begin{align}
    Z_s^v:=\sigma \int_{-\infty}^v K(s,r)dW_r ,\qquad \tilde{B}_s=\sigma \int_v^s K(s,r)dW_r.\label{decomposition passe/futur fBm}
\end{align} 
So we have that $Z$ and $\tilde{B}$ are independent and $Z$ is $\mathcal{F}_v$ measurable. We denote $y_s=\psi_v+\sigma Z_s^v$ and conditionally on $\mathcal{F}_v$ by the Theorem 1.1 of~\cite{MalliavinNualart} we have the Wiener chaos expansion:
\[
f(\theta,\psi_s^v+\sigma B^H_s,\bar\mu_s)= f(\theta,y_s+\sigma \tilde{B}_s,\bar\mu_s)= a_0(s)+a_1(s)\tilde{B}_s+B_s^\perp, 
\]
where $B^\perp$ is the orthogonal of the subspace spanned by 1 and $(\tilde{B}_r)_{r\geq v}$, with:
\begin{align*}
    &a_0(s):=\E[f(\theta,y_s+\sigma \tilde{B}_s,\bar{\mu}_s)],\\
    &a_1(s):=\E[(\tilde{B}_s)^2]^{-1}\E[f(\theta,y_s+\sigma \tilde{B}_s,\bar{\mu}_s)\tilde{B}_s].
\end{align*}
Then by orthogonality we get:
\[
\E[\hat{A}_{st}|\mathcal{F}_v]=(a_0(s)+a_0(t)) Z^v_{st}+ \sigma \E[(a_1(s)\tilde{B}_s+a_1(t)\tilde{B}_t)\tilde{B}_{st}|\mathcal{F}_v].
\]
To check the second condition of the Stochastic Sewing we look at the delta of the increments and using~\eqref{delta of average}:
\[  
\E[ \delta\hat{A}_{sut}|\mathcal{F}_v]=D^0_{sut}+\sigma D^1_{sut}
\]
where
\[
D^0_{sut}:=(a_0(t)-a_0(u))Z^v_{su}-(a_0(u)-a_0(s))Z^v_{ut},
\]
\[
D^1_{sut}:=\E[(a_1(s)\tilde{B}_s+a_1(t)\tilde{B}_t)\tilde{B}_{st}|\mathcal{F}_v]-\E[(a_1(s)\tilde{B}_s+a_1(u)\tilde{B}_u)\tilde{B}_{su}|\mathcal{F}_v]-\E[(a_1(u)\tilde{B}_u+a_1(t)\tilde{B}_t)\tilde{B}_{ut}|\mathcal{F}_v].
\]
We will estimate both terms separately, we start by $D^0$. We need to estimate $a_0(u)-a_0(s)$ but $a_0$ is just the expectation of the function of a Gaussian random variable so we can rewrite it in the following way:
\[ 
a_0(s)=F(\theta,\psi_s^v+Z_s^v, (2H)^{-\frac{1}{2}}(s-v)^H,\bar{\mu}_s)
\]
where $F(\theta,m,\sigma,\nu)=\E[f(\theta,m+\sigma X,\nu)]$ and X is a centered Gaussian random variable. We prove the Lipschitz continuity in $m$ and $\sigma$ of $F$ by proving that its derivatives are bounded:
\begin{align*}
    \partial_m F(\theta,m,\sigma,\nu)&=\partial_m\frac{1}{\sqrt{2\pi}}\int_{\R} e^{-\frac{|x|^2}{2}}f(\theta,\sigma x+m,\nu)dx\\
    &=\frac{1}{\sqrt{2\pi}\sigma}\int_{\R} \partial_me^{-\frac{|y-m|^2}{2\sigma^2}}f(\theta,y,\nu)dx\\
    &=\frac{1}{\sqrt{2\pi}\sigma}\int_{\R} \frac{y-m}{\sigma^2}e^{-\frac{|y-m|^2}{2\sigma^2}}f(\theta,y,\nu)dx\\
    &=\frac{1}{\sqrt{2\pi}\sigma}\int_{\R} xe^{-\frac{|x|^2}{2}}(f(\theta,\sigma x+m,\nu)-f(\theta,m,\nu))dx
\end{align*}
where for the last line we used the fact that a centered Gaussian has expectation zero. By doing similar steps, we can obtain the following: 
\begin{align*}
    &\partial_\sigma F(\theta,m,\sigma,\nu)=\frac{-1}{\sqrt{2\pi}\sigma}\int_{\R} e^{-\frac{|x|^2}{2}}(f(\theta,\sigma x+m,\nu)-f(\theta,m,\nu))dx\\
    & \qquad \qquad \qquad \quad +\frac{1}{\sqrt{2\pi}\sigma}\int_{\R} x^2e^{-\frac{x^2}{2}}(f(\theta,\sigma x+m,\nu)-f(\theta,m,\nu))dx.
\end{align*}
Therefore:
\[ 
|\partial_mF(\theta,m,\sigma,\nu)|+|\partial_\sigma F(\theta,m,\sigma,\nu)| \leq c \|f\|_{Lip(x,\nu)}
\]
with the constant $c$ that does not depend on $\theta$. Using the fact that F is Lipschitz in $m$, $\sigma$ and $\nu$ and the fact that $\bar{\mu}$ is $H-\text{Hölder}$ continuous in $L^p$ we have:
\begin{align*}
    \|a_0(u)-a_0(s)\|_{L^p}&=\|F(\theta,Z_u^v+\psi_u^v,(2H)^{-\frac{1}{2}}(u-v)^H,\bar{\mu}_u)-F(\theta,Z_s^v+\psi_s^v,(2H)^{-\frac{1}{2}}(s-v)^H,\bar{\mu}_s)\|_{L^p}\\
    &\leq \|F(\theta,Z_u^v+\psi_u^v,(2H)^{-\frac{1}{2}}(u-v)^H,\bar{\mu}_u)-F(\theta,Z_u^v+\psi_u^v,(2H)^{-\frac{1}{2}}(u-v)^H,\bar{\mu}_s)\|_{L^p}\\
    &+\|F(\theta,Z_u^v+\psi_u^v,(2H)^{-\frac{1}{2}}(u-v)^H,\bar{\mu}_s)-F(\theta,Z_s^v+\psi_s^v,(2H)^{-\frac{1}{2}}(u-v)^H,\bar{\mu}_s)\|_{L^p}\\
    &+\|F(\theta,Z_s^v+\psi_s^v,(2H)^{-\frac{1}{2}}(u-v)^H,\bar{\mu}_s)-F(\theta,Z_s^v+\psi_s^v,(2H)^{-\frac{1}{2}}(s-v)^H,\bar{\mu}_s)\|_{L^p}\\
    &\leq c (|u-s|^H+ \|Z^v_{su}\|_{L^p}+ (2H)^{-\frac{1}{2}} ((u-v)^H-(s-v)^H)\\ 
    &\leq c (|u-s|^H+\|Z^v_{su}\|_{L^p}).
\end{align*}
Thus we get:
\begin{align*}
    \|D^0_{sut}\|_{L^p}\lesssim |t-s|^H( \|Z^v_{su}\|_{L^p}+\|Z_{ut}^v\|_{L^p})+\|Z^v_{su}\|_{L^p}\|Z^v_{ut}\|_{L^p}.
\end{align*}
We use the inequality (3.4) from~\cite{stochasticsewingupgrade} to estimate the regularity of $Z^v$ which states that, for all $H\in (0,1)$,
    \begin{equation}{\label{eq: estimate Z}}
  \| Z_{su}^v\|_{L^p}\lesssim (s-v)^{H-1}(u-s).      
    \end{equation}
Using this result and the fact that $v<s<u<t$, we get:
\begin{align}
\|D^0_{sut}\|_{L^p}&\lesssim  |s-v|^{-(1-H)}|t-s|^{1+H}+ |s-v|^{-(2-2H)}|t-s|^2\nonumber\\
&\lesssim |s-v|^{-(1-H)}|t-s|^{1+H}. \label{ineq: D0 sto sewing}
\end{align}
When $H\leq \frac13$, one can use the regularity of the derivative. To make it appear we start by rewriting $D^0$:
\begin{align*}
    D^0_{sut}&=(a_0(t)-a_0(u))Z^v_{su}-(a_0(u)-a_0(s))Z^v_{ut}\\
    &=(F(\theta,Z_t^v+\psi_t^v,(2H)^{-\frac{1}{2}}(t-v)^H,\bar{\mu}_t)-F(\theta,Z_t^v+\psi_t^v,(2H)^{-\frac{1}{2}}(t-v)^H,\bar{\mu}_u))Z^v_{su}\\
    &\quad-(F(\theta,Z_s^v+\psi_s^v,(2H)^{-\frac{1}{2}}(s-v)^H,\bar{\mu}_u)-F(\theta,Z_s^v+\psi_s^v,(2H)^{-\frac{1}{2}}(s-v)^H,\bar{\mu}_s))Z^v_{ut}\\
    &\quad+(F(\theta,Z_t^v+\psi_t^v,(2H)^{-\frac{1}{2}}(u-v)^H,\bar{\mu}_u)-F(\theta,Z_u^v+\psi_u^v,(2H)^{-\frac{1}{2}}(u-v)^H,\bar{\mu}_u)\\
    &\qquad\quad-\partial_mF(\theta,Z_u^v+\psi_u^v,(2H)^{-1/2}(u-v)^H,\bar\mu_u)Z^v_{ut})Z^v_{su}\\
    &\quad-(F(\theta,Z_u^v+\psi_u^v,(2H)^{-\frac{1}{2}}(u-v)^H,\bar{\mu}_u)-F(\theta,Z_s^v+\psi_s^v,(2H)^{-\frac{1}{2}}(u-v)^H,\bar{\mu}_u)\\
    &\qquad\quad-\partial_mF(\theta,Z_u^v+\psi_u^v,(2H)^{-1/2}(u-v)^H,\bar\mu_u)Z^v_{su})Z^v_{ut}\\
    &\quad + (F(\theta,Z_t^v+\psi_t^v,(2H)^{-\frac{1}{2}}(t-v)^H,\bar{\mu}_u)-F(\theta,Z_t^v+\psi_t^v,(2H)^{-\frac{1}{2}}(u-v)^H,\bar{\mu}_u))Z^v_{su}\\
    &\quad -(F(\theta,Z_s^v+\psi_s^v,(2H)^{-\frac{1}{2}}(u-v)^H,\bar{\mu}_u)-F(\theta,Z_s^v+\psi_s^v,(2H)^{-\frac{1}{2}}(u-v)^H,\bar{\mu}_u))Z^v_{ut}. 
\end{align*}
We have three different types of terms, we start by those differing only in the mean. We have that:
\begin{align*}
    \partial_{mm}^2F(\theta,m,\sigma,\nu)&=\partial_m\frac{1}{\sqrt{2\pi}\sigma}\int_{\R} xe^{-\frac{x^2}{2}}(f(\theta,\sigma x+mn\nu)-f(\theta,m,\nu))\mathrm{d}x\\
    &=\frac{1}{\sqrt{2\pi}\sigma}\int_{\R} xe^{-\frac{x^2}{2}}(\partial_mf(\theta,\sigma x+m,\nu)-\partial_mf(\theta,m,\nu))\mathrm{d}x.
\end{align*}
Hence, we obtain:
\begin{align}
    |F(\theta,Z_t^v+\psi_t^v,(2H)^{-\frac{1}{2}}(u-v)^H,&\bar{\mu}_u)-F(\theta,Z_u^v+\psi_u^v,(2H)^{-\frac{1}{2}}(u-v)^H,\bar{\mu}_u) \nonumber\\ 
    &\qquad\quad-\partial_mF(\theta,Z_u^v+\psi_u^v,(2H)^{-1/2}(u-v)^H,\bar\mu_u)Z^v_{ut}|\leq c|Z^v_{ut}|^2. \label{ineq: sto sewing F average}
\end{align}
Now, for the difference in variance, we notice that:
\begin{align*}
    \partial_\sigma F(\theta,m,\sigma,\nu)&=\frac{-1}{\sqrt{2\pi}\sigma}\int_{\R} e^{-\frac{|x|^2}{2}}(f(\theta,\sigma x+m,\nu)-f(\theta,m,\nu)-\partial_mf(\theta,m,\nu)\sigma x)dx\\
    &\quad+\frac{1}{\sqrt{2\pi}\sigma}\int_{\R} x^2e^{-\frac{x^2}{2}}(f(\theta,\sigma x+m,\nu)-f(\theta,m,\nu)-\partial_mf(\theta,m,\nu)\sigma x)dx.
\end{align*}
Therefore,
\begin{equation}\label{ineq: sto sewing diff sigma F}
|\partial_\sigma F(\theta,m,\sigma,\nu)|\leq \|\partial_xf\|_{Lip(x,\nu)} \sigma.
\end{equation}
Finally, we prove that the McKean-Vlasov law has more regularity than the process itself:
\begin{align}
    &|F(\theta,Z_t^v+\psi_t^v,(2H)^{-\frac{1}{2}}(t-v)^H,\bar{\mu}_t)-F(\theta,Z_t^v+\psi_t^v,(2H)^{-\frac{1}{2}}(t-v)^H,\bar{\mu}_u)|\nonumber\\
    &=|\E[f(\theta,y_t+\sigma\tilde B_t,\bar\mu_t)-f(\theta,y_t+\sigma\tilde B_t,\bar\mu_u)]|\nonumber\\
    &=|\E[\E'[\partial_\mu f(\theta,y_t+\sigma\tilde B_t,\tilde\mu_u)(\bar X_u')\bar X_{ut}']]|+O(|t-u|^{2H})\nonumber\\
    &=|\E[\E'[\partial_\mu f(\theta,y_t+\sigma\tilde B_t,\tilde\mu_u)(\bar X_u') B_{ut}']]|+O(|t-u|^{2H})\nonumber \\
    &\leq c |\E[\E'[\langle Df(\theta,y_t+\sigma \tilde B_t,\bar\mu_u)(\bar X'_u),\mathbbm{1}_{[u,t]}\rangle_\mathcal{H}]]+O(|t-u|^{2H})\nonumber\\
    &\leq c|\langle \mathbbm{1}_{[0,u]},\mathbbm{1}_{[u,t]}\rangle_\mathcal{H}+O(|t-s|^{2H})\nonumber\\
    &\leq c |t-u|^{2H}.\label{ineq: sto sewing law regularity}
\end{align}
By combining~\eqref{ineq: sto sewing F average},\eqref{ineq: sto sewing diff sigma F},~\eqref{ineq: sto sewing law regularity} and~\eqref{eq: estimate Z} we obtain:
\begin{equation}
    \|D^0_{sut}\|\lesssim |s-v|^{1-H}|t-s|^{1+2H} \label{ineq: sto sewing D0 H<1/3}
\end{equation}
since $t-s\leq s-v$. Now we can estimate $D^1$. We first rewrite it using~\eqref{delta of average}:
\begin{align*}
    D^1_{sut}=&(a_1(t)-a_1(u))\E[\tilde{B}_t\tilde{B}_{st}]+(a_1(s)-a_1(u))\E[\tilde{B}_s\tilde{B}_{st}]\\&
    -(a_1(s)-a_1(u))\E[\tilde{B}_s\tilde{B}_{su}]-(a_1(t)-a_1(u))\E[\tilde{B}_t\tilde{B}_{ut}].
\end{align*}
We focus on $(a_1(t)-a_1(u))\E[\tilde{B}_t\tilde{B}_{st}]$ since the computations for the other terms are the same.
We start by seeing that
\begin{align}\label{estimate B tilde}
|\E[\tilde{B}_s\tilde{B}_{st}]|\lesssim |t-s|^{2H},
\end{align}
a proof can be found in Lemma 3.2 of~\cite{stochasticsewingupgrade}. Thus, we need to bound $a_1(t)-a_1(u)$. As for $a_0(s)$, we can define 
\[a_1(s)=G(\theta, Y_s^v,(2H)^{-\frac{1}{2}}(s-v)^H,\bar{\mu}_s)\]
with
\[
G(\theta,m,\sigma,\nu)=\sigma^{-1}\E[f(\theta,m+\sigma X,\nu)X]=\E[\partial_xf(\theta,m+\sigma X,\nu)]
\]
where $X$ is a standard normal random variable, and the second equality follows from Stein's lemma.
Therefore, replacing $f$ with $\partial_x f$ and following the same steps as for $F$, we obtain
\[
|\partial_m G(\theta,m,\sigma,\nu)|+|\partial_\sigma G(\theta,m,\sigma,\nu)|\lesssim \|\partial_x f\|_{Lip(x,\nu)}.
\]
By the Lipschitz continuity of $G$, and following the same argument as for $a_0$, we obtain:
\begin{align*}
|a_1(t)-a_1(u)| &=\|G(\theta,Z_u^v+\psi_u^v,(2H)^{-\frac{1}{2}}(u-v)^H,\bar{\mu}_u)-G(\theta,Z_s^v+\psi_s^v,(2H)^{-\frac{1}{2}}(s-v)^H,\bar{\mu}_s)\|_{L^p}\\
    &\leq \|G(\theta,Z_u^v+\psi_u^v,(2H)^{-\frac{1}{2}}(u-v)^H,\bar{\mu}_u)-G(\theta,Z_u^v+\psi_u^v,(2H)^{-\frac{1}{2}}(u-v)^H,\bar{\mu}_s)\|_{L^p}\\
    &+\|G(\theta,Z_u^v+\psi_u^v,(2H)^{-\frac{1}{2}}(u-v)^H,\bar{\mu}_s)-G(\theta,Z_s^v+\psi_s^v,(2H)^{-\frac{1}{2}}(u-v)^H,\bar{\mu}_s)\|_{L^p}\\
    &+\|G(\theta,Z_s^v+\psi_s^v,(2H)^{-\frac{1}{2}}(u-v)^H,\bar{\mu}_s)-G(\theta,Z_s^v+\psi_s^v,(2H)^{-\frac{1}{2}}(s-v)^H,\bar{\mu}_s)\|_{L^p}\\
    &\lesssim |t-u|^H+|Z^v_{ut}|.
    \end{align*}
    Thus with~\eqref{eq: estimate Z}, we have:
    \[
    \|a_1(t)-a_1(u)\|_{L^p}\lesssim |t-u|^H +|s-v|^{-(1-H)}|t-s|.
    \]
    We get the final estimate for $D^1$ using~\eqref{estimate B tilde}:
    \begin{equation}\label{ineq:D1 sto sewing}
        \|D^1_{sut}\|_{L^p}\lesssim |t-s|^{3H}+|s-v|^{-(1-H)}|t-s|^{1+2H}.
    \end{equation}
    If $H\leq \frac13$, the only difference is that we use the extra regularity of the McKean-Vlasov law. By replacing $f$ by $\partial_xf$ in the computation of~\eqref{ineq: sto sewing law regularity}, we obtain:
    \[
    \|a_1(t)-a_1(u)\|\lesssim |t-u|^{2H}+\|Z^v_{ut}\|_{L^p}+|s-v|^{-(1-H)}|t-s|.
    \]
    Thus using~\eqref{eq: estimate Z} and~\eqref{estimate B tilde} we obtain:
\begin{equation}\label{ineq: sto sewing D1 H<1/3}
        \|D^1_{sut}\|_{L^p}\lesssim |t-s|^{4H}+|s-v|^{-(1-H)}|t-s|^{1+2H}.
    \end{equation}

For $H\in(\frac13,\frac12)$, by combining~\eqref{ineq: remainder sto sewing},\eqref{ineq: D0 sto sewing},\eqref{ineq:D1 sto sewing} we obtain:
\begin{align*}
\| \E[\delta A_{sut}|\mathcal{F}_v]\|_{L^p}&\lesssim |t-v|^{1+H}|t-s|^H+|s-v|^{-(1-H)}|t-s|^{1+H}+|t-s|^{3H}+|s-v|^{-(1-H)}|t-s|^{1+2H}\\
&\lesssim |t-v|^{1+H}|t-s|^H+|t-s|^{3H}+|s-v|^{-(1-H)}|t-s|^{1+H}.
\end{align*}
Although the bound is not written in exactly the form required
by the stochastic sewing lemma, we can treat its two terms
separately and apply the lemma with the following two choices
of exponents, satisfying condition~\eqref{cond sto sewing}:
\begin{align*}
    &\alpha=1-H,&\qquad \beta_1=1+H,\qquad&\gamma=1+H, \qquad \beta_2=H,\qquad \beta_3=2H;\\
    &\alpha=0,&\qquad \beta_1=3H,\,\quad\qquad&\gamma=1+H, \qquad \beta_2=H,\qquad \beta_3=2H.
\end{align*}
Among the conditions in~\eqref{cond sto sewing}, the most restrictive is $1+\alpha-\beta_1 < \alpha\frac{\gamma+\beta_2-1}{\gamma}$, which is clearly satisfied for both parameter choices of $(\alpha, \beta_1, \gamma, \beta_2)$ outlined above.
The desired convergence follows from the extended stochastic sewing theorem, Theorem 5.5 of~\cite{stochasticsewingupgrade}.

When $H\in(\frac14,\frac13]$, we combine~\eqref{ineq: remainder sto sewing},\eqref{ineq: sto sewing D0 H<1/3},\eqref{ineq: sto sewing D1 H<1/3} to obtain:
\begin{align*}
    \|\E[\delta \tilde A_{sut}|\mathcal{F}_v]\|_{L^p}\lesssim |t-v|^{1+H}|t-s|^H+|s-v|^{-(1-H)}|t-s|^{1+2H}+|t-s|^{4H}.
\end{align*}

As in the previous case, we have two choices of exponents:
\begin{align*}
    \alpha &= 1-H, & \beta_1 &= 1+2H, & \gamma &= 1+H, & \beta_2 &= H, & \beta_3 &= 2H; \\
    \alpha &= 0,   & \beta_1 &= 4H,   & \gamma &= 1+H, & \beta_2 &= H, & \beta_3 &= 2H.
\end{align*}
Once again, one can verify that these exponents satisfy condition~\eqref{cond sto sewing}, thereby establishing convergence for $H \in (\frac{1}{4}, \frac{1}{3}]$.

\end{proof}

We now prove the convergence results for the rough contrast. The first step is an analogue of the decomposition used in the Young case.

\begin{lemma}\label{diff contrast rough}
For $H\in(\frac14,\frac12)$, under Assumptions~\ref{Bdd moment},~\ref{Lip cond}, for every $N,n$ and
$\theta\in\Theta$,
\[
Q_n^N(\theta)-Q_n^N(\theta_0)=\widehat Q_n^N(\theta)-2\sigma\widetilde Q_n^N(\theta)+R_n^N(\theta),
\]
where
\[
\sup_{\theta\in\Theta} |R_n^N(\theta)|\xrightarrow[N,n\to\infty]{}0
\]
almost surely, and
\begin{align*}
\widehat Q_n^N(\theta)&=\frac{\alpha_n}{N}\sum_{i=1}^N\sum_{k=0}^{n-1}\left|\frac{\delta_{\theta,\theta_0}\left(b(\cdot,X_{t_k}^{i,N},\mu_{t_k}^N)+b(\cdot,X_{t_{k+1}}^{i,N},\mu_{t_{k+1}}^N)\right)}{2}\right|^2,\\
\widetilde Q_n^N(\theta)&=\frac1N\sum_{i=1}^N\sum_{k=0}^{n-1}\frac{\delta_{\theta,\theta_0}\left(b(\cdot,X_{t_k}^{i,N},\mu_{t_k}^N)+b(\cdot,X_{t_{k+1}}^{i,N},\mu_{t_{k+1}}^N)\right)}{2}B_{t_kt_{k+1}}^{i,H}.
\end{align*}
\end{lemma}

\begin{proof}
In this part, we use another approximation compared to the Young case of the trajectory of the particle in which the drift is approximated by the average of the left and right point:
\begin{equation}\label{Taylor rough particle}
    X^{i,N}_{t_kt_{k+1}}=\frac{\alpha_n}{2}(b(\theta_0,X^{i,N}_{t_k},\mu^N_{t_k})+b(\theta_0,X^{i,N}_{t_{k+1}},\mu^N_{t_{k+1}}))+\sigma B^{i,H}_{t_kt_{k+1}}+\tilde r^{i,N}_k,
\end{equation}
where the remainder is defined as follows
\begin{equation}\label{remainder average approximation}
    \tilde r^{i,N}_k:=\int_{t_k}^{t_{k+1}}b(\theta_0,X^{i,N}_s,\mu^N_s)-\frac{b(\theta_0,X^{i,N}_{t_k},\mu^N_{t_k})+b(\theta_0,X^{i,N}_{t_{k+1}},\mu^N_{t_{k+1}})}{2}\mathrm{d}s.
\end{equation}
We compute the difference by first expanding the square and then using the previous expansion:
\begin{align*}
    &N\alpha_n(Q_n^N(\theta)-Q_n^N(\theta_0)) =\sum_{i=1}^N\sum_{k=0}^{n-1} (X^{i,N}_{t_kt_{k+1}}-\frac{\alpha_n}{2}b(\theta,X^{i,N}_{t_k},\mu^N_{t_k})+b(\theta,X^{i,N}_{t_{k+1}},\mu^N_{t_{k+1}}))^2  \\
    &-(X^{i,N}_{t_kt_{k+1}}-\frac{\alpha_n}{2}b(\theta_0,X^{i,N}_{t_k},\mu^N_{t_k})+b(\theta_0,X^{i,N}_{t_{k+1}},\mu^N_{t_{k+1}}))^2\\
    &=\sum_{i=1}^N\sum_{k=0}^{n-1} \frac{\alpha_n^2}{4}|b(\theta,X^{i,N}_{t_k},\mu^N_{t_k})+b(\theta,X^{i,N}_{t_{k+1}},\mu^N_{t_{k+1}})|^2 -\sum_{i=1}^N\sum_{k=0}^{n-1}\frac{\alpha_n^2}{4}|b(\theta_0,X^{i,N}_{t_k},\mu^N_{t_k})+b(\theta_0,X^{i,N}_{t_{k+1}},\mu^N_{t_{k+1}})|^2\\
    &-\sum_{i=1}^N\sum_{k=0}^{n-1}\alpha_n\delta_{\theta,\theta_0}(b(.,X^{i,N}_{t_k},\mu^N_{t_k})+b(.,X^{i,N}_{t_{k+1}},\mu^N_{t_{k+1}}))X^{i,N}_{t_kt_{k+1}}.
\end{align*}
We use the expansion~\eqref{Taylor rough particle}:
\begin{align*}
   & N\alpha_n(Q_n^N(\theta)-Q_n^N(\theta_0))\\ &=\sum_{i=1}^N\sum_{k=0}^{n-1}\frac{\alpha_n^2}{4}|b(\theta,X^{i,N}_{t_k},\mu^N_{t_k})+b(\theta,X^{i,N}_{t_{k+1}},\mu^N_{t_{k+1}})|^2+ \frac{\alpha_n^2}{4}|b(\theta_0,X^{i,N}_{t_k},\mu^N_{t_k})+b(\theta_0,X^{i,N}_{t_{k+1}},\mu^N_{t_{k+1}})|^2 \\
   &-2\alpha_n^2\frac{(b(\theta,X^{i,N}_{t_k},\mu^N_{t_k})+b(\theta,X^{i,N}_{t_{k+1}},\mu^N_{t_{k+1}}))
    (b(\theta,X^{i,N}_{t_k},\mu^N_{t_k})+b(\theta,X^{i,N}_{t_{k+1}},\mu^N_{t_{k+1}}))}{2}\\
    &-\sum_{i=1}^N\sum_{k=0}^{n-1}\alpha_n^2(b(\theta,X^{i,N}_{t_k},\mu^N_{t_k})+b(\theta,X^{i,N}_{t_{k+1}},\mu^N_{t_{k+1}}))(b(\theta_0,X^{i,N}_{t_k},\mu^N_{t_k})+b(\theta_0,X^{i,N}_{t_{k+1}},\mu^N_{t_{k+1}})) \\
    &-2\sigma\sum_{i=1}^N\sum_{k=0}^{n-1}\alpha_nB^{i,H}_{t_kt_{k+1}}\frac{\delta_{\theta,\theta_0}b(\cdot,X^{i,N}_{t_{k+1}},\mu^N_{t_{k+1}})+b(\cdot,X^{i,N}_{t_k},\mu^N_{t_k})}{2}\\
    & \underbrace{-2\sum_{i=1}^N\sum_{k=0}^{n-1}\alpha_n \tilde r^{i,N}_k\frac{\delta_{\theta,\theta_0}b(\cdot,X^{i,N}_{t_{k+1}},\mu^N_{t_{k+1}})+b(\cdot,X^{i,N}_{t_k},\mu^N_{t_k})}{2}}_{:=N\alpha_nR^N_n(\theta)}=N\alpha_n(\hat Q_n^N(\theta) -2\sigma \tilde Q_n^N(\theta)+R_n^N(\theta)).
    \end{align*} 

    The remainder can be controlled by the same estimates used
to prove \eqref{estimate remainder}, applied to the symmetric drift approximation.
This yields that $R_n^N$ goes to $0$ almost surely uniformly in $\theta$.
\end{proof}

\begin{lemma}\label{CV somme Riemann middle approx}
Let $f,g:\Theta\times\R\times\mathcal{P}_2(\R)\to\R$ satisfy the Lipschitz and polynomial-growth conditions of Lemma~\ref{CV somme Riemann classique}. Then, for every $H\in(\frac14,\frac12)$,
\begin{align*}
&\sup_{\theta\in\Theta}\Bigg|\frac{\alpha_n}{N}\sum_{i=1}^N\sum_{k=0}^{n-1}\frac{f(\theta,X_{t_k}^{i,N},\mu_{t_k}^N)+f(\theta,X_{t_{k+1}}^{i,N},\mu_{t_{k+1}}^N)}{2}\frac{g(\theta,X_{t_k}^{i,N},\mu_{t_k}^N)+g(\theta,X_{t_{k+1}}^{i,N},\mu_{t_{k+1}}^N)}{2}\\
&\hspace{2cm}-\E\left[\int_0^T f(\theta,\bar X_s^1,\bar\mu_s)g(\theta,\bar X_s^1,\bar\mu_s)\,\mathrm{d}s\right]\Bigg|\xrightarrow[N,n\to\infty]{}0
\end{align*}
almost surely. In particular, under Assumptions~\ref{Bdd moment} and~\ref{Lip cond}
\[
\sup_{\theta\in\Theta}\left|\widehat Q_n^N(\theta)-\E\left[\int_0^T|\delta_{\theta,\theta_0}b(\cdot,\bar X_s,\bar\mu_s)|^2\,\mathrm{d}s\right]\right|\xrightarrow[N,n\to\infty]{}0
\]
almost surely.
\end{lemma}

\begin{proof}
The proof follows the same decomposition as in Lemma~\ref{CV somme Riemann classique}. The additional endpoint terms introduced by the symmetric approximation are controlled using the Hölder regularity of the particle system.
\end{proof}

We next establish the convergence of the stochastic integral term.

\begin{lemma}\label{CV somme Riemann fBm, rough case}
Let $H\in(\frac14,\frac12)$ and $f:\Theta\times\R\times\mathcal{P}_2(\R)\to\R$ satisfy the same  assumptions as in Lemma \ref{lemma Strat integral}. Then
\[
\sup_{\theta\in\Theta} \left| \frac1N\sum_{i=1}^N\sum_{k=0}^{n-1}\frac{f(\theta,X_{t_k}^{i,N},\mu_{t_k}^N)+f(\theta,X_{t_{k+1}}^{i,N},\mu_{t_{k+1}}^N)}{2}B_{t_kt_{k+1}}^{i,H}-\E\left[\int_0^Tf(\theta,\bar X_s^1,\bar\mu_s)\circ\mathrm{d}B_s^{1,H}\right]\right|\to0\]
almost surely. Consequently, under Assumptions~\ref{Bdd moment},~\ref{Lip cond} and $b$ satisfy the same assumptions as in Lemma~\ref{lemma Strat integral}
\[
\sup_{\theta\in\Theta}\left|\widetilde Q_n^N(\theta)-\E\left[\int_0^T \delta_{\theta,\theta_0} b(\cdot,\bar X_s^1,\bar\mu_s) \circ\mathrm{d}B_s^{1,H} \right]\right|\to0
\]
almost surely.
\end{lemma}

\begin{proof}
For this case, we are going to change the order of the arguments meaning compared to the proof of Lemma~\ref{CV somme Riemann fBm, smooth case}. We will first use the propagation of chaos then the law of large numbers and finally the convergence of the partition. We choose this order because we want to have the limit measure which is non random when we do the convergence of the discretizations towards the integral. Thus we keep the definition of $L^2(\theta)$ as in the proof of Lemma~\ref{lemma diff contrast fct Young case} and we start splitting:
    \begin{align*}
        \sup_{\theta \in \Theta}|&\frac{1}{N}\sum_{i=1}^N \sum_{k=0}^{n-1}\frac{f(\theta,X^{i,N}_{t_k},\mu_{t_k}^N)+f(\theta,X^{i,N}_{t_{k+1}},\mu^N_{t_{k+1}})}{2}B^{i,H}_{t_kt_{k+1}}-L^2(\theta)|\\
        &\leq  \sup_{\theta \in \Theta}|\frac{1}{N}\sum_{i=1}^N \sum_{k=0}^{n-1}\left(\frac{f(\theta,X^{i,N}_{t_k},\mu_{t_k}^N)+f(\theta,X^{i,N}_{t_{k+1}},\mu^N_{t_{k+1}})}{2}-\frac{f(\theta,\bar{X}^i_{t_k},\bar{\mu}_{t_k})+f(\theta,\bar{X}^i_{t_{k+1}},\bar{\mu}_{t_{k+1}})}{2}\right)B^{i,H}_{t_kt_{k+1}}|\\
        &+ \sup_{\theta \in \Theta}|\frac{1}{N}\sum_{i=1}^N \sum_{k=0}^{n-1}\frac{f(\theta,\bar{X}^i_{t_k},\bar{\mu}_{t_k})+f(\theta,\bar{X}^i_{t_{k+1}},\bar{\mu}_{t_{k+1}})}{2}B^{i,H}_{t_kt_{k+1}}-L^2(\theta)|.
    \end{align*}
For the first term, we use the propagation of chaos as in Theorem~\ref{prop of chaos}, we take $\varepsilon>0$
 \begin{align*}
     \E[&\sup_{\theta \in \Theta}|\frac{1}{N}\sum_{i=1}^N \sum_{k=0}^{n-1}\left(\frac{f(\theta,X^{i,N}_{t_k},\mu_{t_k}^N)+f(\theta,X^{i,N}_{t_{k+1}},\mu^N_{t_{k+1}})}{2}-\frac{f(\theta,\bar{X}^i_{t_k},\bar{\mu}_{t_k})+f(\theta,\bar{X}^i_{t_{k+1}},\bar{\mu}_{t_{k+1}})}{2}\right)B^{i,H}_{t_kt_{k+1}}|^p]\\
     &\leq \E[\sup_{\theta \in \Theta}|\sum_{k=0}^{n-1}\left(\frac{f(\theta,X^{1,N}_{t_k},\mu_{t_k}^N)+f(\theta,X^{1,N}_{t_{k+1}},\mu^N_{t_{k+1}})}{2}-\frac{f(\theta,\bar{X}^1_{t_k},\bar{\mu}_{t_k})+f(\theta,\bar{X}^1_{t_{k+1}},\bar{\mu}_{t_{k+1}})}{2}\right)B^{1,H}_{t_kt_{k+1}}|^p]\\
     &\leq \E[ (\sum_{k=0}^{n-1} \| f(\theta, X^{1,N},\mu^N)-f(\theta,\bar{X}^1,\bar{\mu}^1)\|_{H-\varepsilon,[0,t_{k+1}]}^p \|B^{1,H}\|_{H-\varepsilon,[t_kt_{k+1}]}\alpha_n^{H-\varepsilon} t_{k+1}^{H-\varepsilon})^p ]\\
     &\leq c \frac{1}{N^{\gamma p}} (\alpha_n\sum_{k=0}^{n-1}t_{k+1})^{p(H-\varepsilon)}.
 \end{align*}
We used the fact that the particles are identically distributed and then the regularity estimates on the solutions~\ref{prop solution}. Regarding the second part, we proceed in two steps: first, we prove the convergence for a fixed $\theta$, then we will use Ascoli-Arzela to make the convergence uniform in $\theta$. To start, we are going to use Lemma~\ref{Array LLN}, which is an array Law of large Numbers. 

We fix $\theta \in \Theta$ and we define $Y_n^i=\sum_{k=0}^{n-1}\frac{f(\theta,\bar{X}^i_{t_k},\bar{\mu}_{t_k})+f(\theta,\bar{X}^i_{t_{k+1}},\bar{\mu}_{t_{k+1}})}{2}B^{i,H}_{t_kt_{k+1}}$. Hence, using Lemma~\ref{lemma Strat integral}, we know that there exists $\varepsilon>0$ such that for all $p\geq2$:
\[
\|Y_n^i-\int_0^Tf(\theta,\bar X^i_s,\bar \mu_s)\circ \mathrm{d}B^{i,H}\|_{L^p(\Omega)}\lesssim \alpha_n^\varepsilon.
\]
We check the assumptions of the previous lemma, and hence we can conclude that for each $\theta$, we have the wanted convergence. \\
To make the convergence uniform, we use Ascoli-Arzela, hence we need to show the almost sure equicontinuity. To prove it, we use the assumption on the derivative of $f$ and the fact that the bound on the convergence of the Stratonovich integration does not depend on $\theta$. Let $\theta,\tilde\theta \in\Theta$, we do a Taylor expansion in $\theta$:
\begin{align*}
    |\frac{1}{N}\sum_{i=1}^N Y_n^i(\theta)-\frac{1}{N}\sum_{i=1}^N Y_n^i(\tilde\theta)|&\leq |\frac{1}{N}\sum_{i=1}^N \sum_{k=0}^{n-1}\frac{\delta_{\theta,\tilde\theta}f(.,\bar X^i_{t_k},\bar\mu_{t_k})+f(.,\bar X^i_{t_{k+1}},\bar\mu_{t_{k+1}})}{2}B^{i,H}_{t_kt_{k+1}}|\\
    &\leq |\theta-\tilde\theta| |\frac{1}{N}\sum_{i=1}^N \sum_{k=0}^{n-1}\frac{f(\xi,\bar X^i_{t_k},\bar\mu_{t_k})+f(\xi,\bar X^i_{t_{k+1}},\bar\mu_{t_{k+1}})}{2}B^{i,H}_{t_kt_{k+1}}|,
\end{align*}
with $\xi\in[\theta,\tilde\theta]$. Then, by the law of large numbers, the sum is finite almost surely and with a bound that does not depend on $\theta$. Thus, it is equicontinuous. Therefore, the convergence is uniform and it concludes the proof.
\end{proof}

\begin{proof}[Proof of Theorem~\ref{CV contrast rough case}]
The first statement follows by combining Lemmas \ref{diff contrast rough}, \ref{CV somme Riemann middle approx} and \ref{CV somme Riemann fBm, rough case}. For the second statement, at $\theta=\theta_0$ we have
\begin{align*}
Q_n^N(\theta_0)&=\frac{\sigma^2}{N\alpha_n}\sum_{i=1}^N\sum_{k=0}^{n-1}\left(|B_{t_kt_{k+1}}^{i,H}|^2-\alpha_n^{2H}\right)+\frac{2}{N\alpha_n}\sum_{i=1}^N\sum_{k=0}^{n-1}\widetilde r_k^{i,N}B_{t_kt_{k+1}}^{i,H}+\frac1{N\alpha_n}\sum_{i=1}^N\sum_{k=0}^{n-1}|\widetilde r_k^{i,N}|^2,
\end{align*}
where we recall that $\tilde r^{i,N}_k$ is as in~\eqref{remainder average approximation} and we also introduce the notation $\tilde r^i_k$ for its McKean-Vlasov counterparts. The first and third terms are controlled as in the Young case. The second term requires a more precise analysis because $H<\frac12$. For the remaining term we use the same type of strategy as in the convergence of contrast function in Proposition~\ref{CV contrast in theta 0} meaning we use the propagation of chaos and the law of large numbers to replace the interacting particles by the mean field limit:
    \begin{align}
       & |\frac{1}{N\alpha_n}\sum_{i=1}^N\sum_{k=0}^{n-1}\tilde r^{i,N}_k B^{i,H}_{t_kt_{k+1}}|\leq  |\frac{1}{N\alpha_n}\sum_{i=1}^N\sum_{k=0}^{n-1}\bar r^{i,N}_kB^{i,H}_{t_kt_{k+1}}-\frac{1}{N\alpha_n}\sum_{i=1}^N\sum_{k=0}^{n-1}\bar r^{i}_k B^{i,H}_{t_kt_{k+1}}|\nonumber\\
        &+| \frac{1}{N\alpha_n}\sum_{i=1}^N\sum_{k=0}^{n-1}\bar r^{i}_k B^{i,H}_{t_kt_{k+1}}-\frac{1}{\alpha_n} \E[\sum_{k=0}^{n-1}\tilde r^1_k B^H_{t_kt_{k+1}}]|+|\frac{1}{\alpha_n} \E[\sum_{k=0}^{n-1}\bar r^1_k B^H_{t_kt_{k+1}}]|.\label{decomposition product remainder fbm}
    \end{align}
    We study each term separately. For the first one, we use the propagation of chaos:
    \begin{align}
        \E[|\frac{1}{N\alpha_n}\sum_{i=1}^N\sum_{k=0}^{n-1}\tilde r^{i,N}_kB^{i,H}_{t_kt_{k+1}}-\frac{1}{N\alpha_n}\sum_{i=1}^N\sum_{k=0}^{n-1}\tilde r^{i}_k B^{i,H}_{t_kt_{k+1}}|^p]&\leq \alpha_n^{1-2p}\sum_{k=0}^{n-1}\E[|B_{t_kt_{k+1}}|^p|\bar r^i_k-\tilde r^{i,N}_k|^p] \nonumber \\ 
        &\leq \alpha_n^{1-2p}\sum_{k=0}^{n-1}\E[|B_{t_kt_{k+1}}|^{2p}]^{\frac{1}{2}}\E[|\bar r^i_k-\tilde r^{i,N}_k|^{2p}]^{\frac{1}{2}}.\label{estimate diff between product fbm and remainder}
    \end{align}
    We study the difference between $\bar r^i_k$ and $\tilde r^{i,N}_k$:
    \begin{align*}
        \E[|\bar r^i_k-\tilde r^{i,N}_k|^p]&=\E[|\int_{t_k}^{t_{k+1}}b(\theta_0,\bar X_s^i,\bar\mu_s)-\frac{b(\theta_0,\bar X^{i}_{t_k},\bar\mu_{t_k})+b(\theta_0,\bar X^{i}_{t_{k+1}},\bar\mu_{t_{k+1}})}{2}\\ 
        &-(b(\theta_0,  X_s^{i,N},\mu^N_s)-\frac{b(\theta_0,X^{i,N}_{t_k},\mu^N_{t_k})+b(\theta_0,X^{i,N}_{t_{k+1}},\mu^N_{t_{k+1}})}{2})\mathrm{d}s|^p]\\
        &\leq c \alpha_n^{p(1+H-\varepsilon)}\E[\|b(\theta_0,\bar X^i,\bar \mu)-b(\theta_0,X^{i,N},\mu^N)\|_{H-\varepsilon,[t_kt_{k+1}]}^p] \leq c\frac{1}{N^{\gamma p}}\alpha_n^{p(1+H-\varepsilon)}.
    \end{align*}
    We plug this estimation in~\eqref{estimate diff between product fbm and remainder}:
    \begin{align}
        \E[|\frac{1}{N\alpha_n}\sum_{i=1}^N\sum_{k=0}^{n-1}\tilde r^{i,N}_kB^{i,H}_{t_kt_{k+1}}-\frac{1}{N\alpha_n}\sum_{i=1}^N\sum_{k=0}^{n-1}\tilde r^{i}_k B^{i,H}_{t_kt_{k+1}}|^p]&\leq c \alpha_n^{-2p}\alpha_n^{pH}\alpha_n^{p(1+H-\varepsilon)}N^{-\gamma p}\nonumber\\
        &\leq c \left(\frac{n^{2-4H-2\varepsilon}}{N^{2\gamma}}\right)^{p/2}. \label{estimate precise sum remainder 1}
    \end{align}
 By choosing correctly $\gamma$ and $\varepsilon$, and using Assumption~\ref{CV rate step size 1}, we can get the convergence towards 0.\\ 
    For the second term of~\eqref{decomposition product remainder fbm}, we use the Rosenthal inequality given in Proposition~\ref{Rosenthal inequality}:
    \begin{align}
      &  \E[| \frac{1}{N\alpha_n}\sum_{i=1}^N\sum_{k=0}^{n-1}\tilde r^{i}_k B^{i,H}_{t_kt_{k+1}}-\frac{1}{\alpha_n} \E[\sum_{k=0}^{n-1}\tilde r^1_k B^H_{t_kt_{k+1}}]|] \\
      &\leq \frac{c}{N^p}\left( \sum_{i=1}^N\E[|\frac{1}{\alpha_n}\sum_{k=0}^{n-1}B_{t_kt_{k+1}}^{i,H}\tilde r^i_k-\E[B_{t_kt_{k+1}}^{i,H}\tilde r^i_k]|^2]\right)^{p/2}
        +\frac{c}{N^p}\sum_{i=1}^N\E[|\frac{1}{\alpha_n}\sum_{k=0}^{n-1}B_{t_kt_{k+1}}^{i,H}\tilde r^i_k-\E[B_{t_kt_{k+1}}^{i,H}\tilde r^i_k]|^p]\nonumber\\
        &\leq c \left(\frac{n^{2-4H}}{N}\right)^{p/2}+c\frac{n^{p(1-2H)}}{N^{p-1}}. \label{estimate precise sum remainder 2}
    \end{align}
   This follows from the regularity estimates for $B^H$ and $\bar{r}^{i}$. Invoking Assumption~\ref{CV rate step size 1} once more, we conclude that this term converges to zero. The final term in~\eqref{decomposition product remainder fbm} requires a more refined analysis. To this end, we begin by rewriting it using the definition of $\bar{r}^i_k$:
    \begin{align*}
        & \E[\frac{1}{\alpha_n}\sum_{k=0}^{n-1}B^{H}_{t_kt_{k+1}}\tilde r^i_k]=\E[\frac{1}{\alpha_n}\sum_{k=0}^{n-1}B^H_{t_kt_{k+1}}\int_{t_k}^{t_{k+1}}b(\theta_0,\bar X_s,\bar\mu_s)-\frac{1}{2}(b(\theta_0,\bar X_{t_k},\bar\mu_{t_k})+b(\theta_0,\bar X_{t_{k+1}},\bar\mu_{t_{ k+1}}))\mathrm{d}s]\\
        &=\E[\frac{1}{\alpha_n}\sum_{k=0}^{n-1}B^H_{t_kt_{k+1}}\int_{t_k}^{t_{k+1}}b(\theta_0,\bar X_s,\bar\mu_s)-b(\theta_0,\bar X_{t_k},\bar\mu_{t_k})\mathrm{d}s -\frac{\alpha_n}{2}(b(\theta_0,\bar X_{t_{k+1}},\bar\mu_{t_{k+1}})-b(\theta_0,\bar X_{t_k},\bar\mu_{t_k})]\\
        &=\E[\frac{1}{\alpha_n}\sum_{k=0}^{n-1}\partial_xb(\theta_0,\bar X_{t_k},\bar\mu_{t_k})B^H_{t_kt_{k+1}}\int_{t_k}^{t_{k+1}}\bar X_{t_ks}\mathrm{d}s-\frac{\alpha_n}{2}\bar X_{t_kt_{k+1}}] \\
        & \qquad +\E[\frac{1}{\alpha_n}\sum_{k=0}^{n-1}B^H_{t_kt_{k+1}}\E[\partial_\mu b(\theta_0,\bar X_{t_k},\bar\mu_{t_k})(\bar X_{t_k}')\int_{t_k}^{t_{k+1}}\bar X_{t_ks}\mathrm{d}s-\frac{\alpha_n}{2}\bar X_{t_kt_{k+1}}]]+O(\alpha_n^{3H-1})\\
        &=\E[\frac{1}{\alpha_n}\sum_{k=0}^{n-1}\partial_xb(\theta_0,\bar X_{t_k},\bar\mu_{t_k})B^H_{t_kt_{k+1}}\int_{t_k}^{t_{k+1}}\bar X_{t_ks}\mathrm{d}s-\frac{\alpha_n}{2}\bar X_{t_kt_{k+1}}]+O(\alpha_n^{3H-1}).
    \end{align*}
   Applying a Taylor expansion in space and measure, we observe that the measure derivative vanishes, as it is taken under an expectation with respect to a centered fBm. By bounding the Taylor remainder for the limiting particle system, we conclude that the noise is the sole term not absorbed into the remainder.
    \begin{align*}
        \E[\frac{1}{\alpha_n}\sum_{k=0}^{n-1} B^H_{t_kt_{k+1}}\tilde r^i_k]&=\E[\frac{\sigma}{\alpha_n}\sum_{k=0}^{n-1}  \partial_xb(\theta_0,\bar X_{t_k},\bar{\mu}_{t_k})B^H_{t_kt_{k+1}}\left(\int_{t_k}^{t_{k+1}} B^H_{t_ks}\mathrm{d}s-\frac{\alpha_n}{2}B^H_{t_kt_{k+1}}\right)]+O(\alpha_n^{3H-1})\\
        &=\frac{\sigma}{\alpha_n}\sum_{k=0}^{n-1} \text{Cov}(\partial_xb(\theta_0,\bar X_{t_k},\bar{\mu}_{t_k}),B^H_{t_kt_{k+1}}\left(\int_{t_k}^{t_{k+1}} B^H_{t_ks}\mathrm{d}s-\frac{\alpha_n}{2}B^H_{t_kt_{k+1}}\right))+O(\alpha_n^{3H-1}).
    \end{align*}
   The expectation can be recast as a covariance since, for all $k\geq 0$,
\begin{equation}
    \E\left[B^H_{t_kt_{k+1}}\left(\int_{t_k}^{t_{k+1}} B^H_{t_ks}\,\mathrm{d}s-\frac{\alpha_n}{2}B^H_{t_kt_{k+1}}\right)\right]=0.
\end{equation}
Following the proof of Lemma~\ref{lemma Strat integral}, we rely on the Mandelbrot--Van Ness representation of the fBm to decompose it into two independent components. Adopting the notation from~\eqref{decomposition passe/futur fBm} and setting $Z^s_{st}=Z_{st}$, we observe that $\tilde{B}_{st} \perp\!\!\!\perp \mathcal{F}_{t_k}$. Utilizing this independence, we obtain for all $k\geq 0$:
    \begin{align*}
        \E[\tilde r^i_k B^H_{t_kt_{k+1}}]&=\text{Cov}(\partial_xb(\theta_0,\bar X_{t_k},\bar\mu_{t_k}),Z_{t_kt_{k+1}}+\tilde B_{t_kt_{k+1}}\int_{t_k}^{t_{k+1}}Z_{t_ks}+\tilde B_{t_ks}\mathrm{d}s +\frac{\alpha_n}{2}(Z_{t_kt_{k+1}}+\tilde B_{t_kt_{k+1}}))\\
        &=\int_{t_k}^{t_{k+1}}\text{Cov}(\partial_xb(\theta_0,\bar X_{t_k},\bar\mu_{t_k}),Z_{t_kt_{k+1}}(Z_{t_ks}-\frac{1}{2}Z_{t_kt_{k+1}}))\mathrm{d}s.
    \end{align*}
Since the process $Z$ is a stochastic integral, we can apply the associated product rule (see Proposition 1.1.3 in~\cite{MalliavinNualart}):
\begin{align*}
    Z_{t_kt_{k+1}}&\left(Z_{t_ks}-\frac{1}{2}Z_{t_kt_{k+1}}\right)\\
    &=\int_{-\infty}^{t_k} \bigl(K(t_{k+1},r)-K(t_k,r)\bigr)\mathrm{d}W_r \int_{-\infty}^{t_k} \left(K(s,r)-\frac{1}{2}\bigl(K(t_{k+1},r)+K(t_k,r)\bigr)\right)\mathrm{d}W_r\\
    &=\int_{-\infty}^{t_k}\int_{-\infty}^{t_k} \bigl(K(t_{k+1},\cdot)-K(t_k,\cdot)\bigr)\left(K(s,\cdot)-\frac{1}{2}\bigl(K(t_{k+1},\cdot)+K(t_k,\cdot)\bigr)\right)(r,q)\,\mathrm{d}W_r\,\mathrm{d}W_q\\
    &\quad+\left\langle K(t_{k+1},\cdot)-K(t_k,\cdot),K(s,\cdot)-\frac{1}{2}\bigl(K(t_{k+1},\cdot)+K(t_k,\cdot)\bigr)\right\rangle_\mathcal{H}.
\end{align*}
We omit the explicit expression for the scalar product since, being a deterministic quantity, it vanishes when computing the covariance. Furthermore, utilizing the fact that Skorokhod integrals have zero expectation, we can express the entire term as an expectation. To proceed with integration by parts, we denote the first- and second-order Malliavin derivatives with respect to $W$ by $D^W$ and $D^{2,W}$, respectively~\cite{MalliavinNualart}. This yields:
\begin{align*}
    \E[\tilde r^i_k B^H_{t_kt_{k+1}}]&=\int_{t_k}^{t_{k+1}}\E[\partial_xb(\theta_0,\bar X_{t_k},\bar\mu_{t_k})\int_{-\infty}^{t_k}\int_{-\infty}^{t_k}(K(t_{k+1},\cdot)-K(t_k,\cdot)) (K(s,\cdot)\\
    &\qquad-\frac{1}{2}(K(t_{k+1},\cdot)+K(t_k,\cdot)))(r,q)dW_rdW_q]\mathrm{d}s\\
    &=\int_{t_k}^{t_{k+1}}\int_{-\infty}^{t_k}\int_{-\infty}^{t_k} \E[D^{2,W}_{r,q}\partial_xb(\theta_0,\bar X_{t_k},\bar\mu_{t_k})(K(t_{k+1},\cdot)-K(t_k,\cdot)) \\
    &\qquad(K(s,\cdot)-\frac{1}{2}(K(t_{k+1},\cdot)+K(t_k,\cdot)))(r,q)]\mathrm{d}r\mathrm{d}q\mathrm{d}s\\
    &=\int_{t_k}^{t_{k+1}}\int_{-\infty}^{t_k}\int_{-\infty}^{t_k} \E[\partial_x^3b(\theta_0,\bar X_{t_k},\bar\mu_{t_k})D^{2,W}_{r,q}\bar X_{t_k}](K(t_{k+1},\cdot)-K(t_k,\cdot)) \\
    &\qquad(K(s,\cdot)-\frac{1}{2}(K(t_{k+1},\cdot)+K(t_k,\cdot)))(r,q)\mathrm{d}r\mathrm{d}q\mathrm{d}s\\
    &\leq \int_{t_k}^{t_{k+1}}\int_{-\infty}^{t_k}\int_{-\infty}^{t_k} \E[\partial_x^3b(\theta_0,\bar X_{t_k},\bar\mu_{t_k})|r-t_k|^{H-\frac{1}{2}}||q-t_k|^{H-\frac{1}{2}}M](K(t_{k+1},\cdot)-K(t_k,\cdot))\\
    &\qquad (K(s,\cdot)-\frac{1}{2}(K(t_{k+1},\cdot)+K(t_k,\cdot)))(r,q)\mathrm{d}r\mathrm{d}q\mathrm{d}s\\
    &\leq c \int_{t_k}^{t_{k+1}}\int_{-\infty}^{t_k}\int_{-\infty}^{t_k}|r-t_k|^{H-\frac{1}{2}}|q-t_k|^{H-\frac{1}{2}}(K(t_{k+1},\cdot)-K(t_k,\cdot))\\ 
    &\qquad(K(s,\cdot)-\frac{1}{2}(K(t_{k+1},\cdot)+K(t_k,\cdot)))(r,q)\mathrm{d}r\mathrm{d}q\mathrm{d}s 
\end{align*}
with $M$ given in Lemma \ref{lemma bound DW}. We obtain a purely deterministic integral, and the change of variable $x=\frac{t_k-r}{\alpha_n}, y=\frac{t_k-q}{\alpha_n}$ yields:
\begin{align*}
    &\int_{-\infty}^{t_k}\int_{-\infty}^{t_k}|r-t_k|^{H-\frac{1}{2}}|q-t_k|^{H-\frac{1}{2}}(K(t_{k+1},\cdot)-K(t_k,\cdot))(K(s,\cdot)-\frac{1}{2}(K(t_{k+1},\cdot)+K(t_k,\cdot)))(r,q)\mathrm{d}r\mathrm{d}q\\
    &\leq \int_0^\infty \int_0^\infty \alpha_n^{H-\frac{1}{2}}x^{H-\frac{1}{2}}\alpha_n^{H-\frac{1}{2}}y^{H-\frac{1}{2}}(\alpha_n^{H-\frac{1}{2}}|x+1|^{H-\frac{1}{2}}-\alpha_n^{H-\frac{1}{2}}x^{H-\frac{1}{2}}) (\alpha_n^{H-\frac{1}{2}}|y+1|^{H-\frac{1}{2}}-\alpha_n^{H-\frac{1}{2}}y^{H-\frac{1}{2}})\frac{\alpha_n^2}{2}\mathrm{d}x\mathrm{d}y\\
    &= \frac{1}{2}\alpha_n^{4H}I^2,
\end{align*}
where $I=\int_0^\infty |x+1|^{H-\frac{1}{2}}-x^{H-\frac{1}{2}}<\infty$. By putting everything together, we get:
\begin{align}
|\E[\frac{1}{\alpha_n}\sum_{k=0}^{n-1}\tilde r _k^1 B^{1,H}_{t_kt_{k+1}}]|\leq c\alpha_n^{3H-1}.\label{estimate precise sum remainder 3}
\end{align}
Thus, this term vanishes for $H\in(1/3,1/2)$.

When $H\in(1/4,1/3]$, the previous bound is no longer sufficient to ensure convergence since $3H-1\leq 0$. We therefore require a more delicate estimate. Specifically, we expand $b(\theta_0,\bar X_s,\bar \mu_s)$ to second order via Taylor's formula:
\begin{align*}
    \E\left[\frac{1}{\alpha_n}\sum_{k=0}^{n-1} B^H_{t_kt_{k+1}}\tilde r^i_k\right]
&=\E\left[\frac{\sigma}{\alpha_n}\sum_{k=0}^{n-1} \partial_xb(\theta_0,\bar X_{t_k},\bar{\mu}_{t_k})B^H_{t_kt_{k+1}}\left(\int_{t_k}^{t_{k+1}} B^H_{t_ks}\,\mathrm{d}s-\frac{\alpha_n}{2}B^H_{t_kt_{k+1}}\right)\right]\\
&+\frac{\sigma^2}{2\alpha_n}\sum_{k=0}^{n-1}\E\left[\partial^2_{xx}b(\theta_0,\bar X_{t_k},\bar{\mu}_{t_k})B^H_{t_kt_{k+1}}\left(\int_{t_k}^{t_{k+1}} |B^H_{t_ks}|^2\,\mathrm{d}s-\frac{\alpha_n}{2}|B^H_{t_kt_{k+1}}|^2\right)\right]+{O}(\alpha_n^{4H-1}).
\end{align*}
In the Taylor expansion, we omit the drift term, which is absorbed into the remainder due to its higher regularity. Furthermore, as in the first expansion, the measure derivatives vanish in expectation. Focusing on the second-order term, we decompose the product of the fBm increments into Wiener chaoses:
\begin{align*}
    B_{t_kt_{k+1}}^H|B_{t_ks}|^2&=\mathrm{I}_1(\mathbbm{1}_{[t_kt_{k+1}]})\mathrm{I}_1(\mathbbm{1}_{[t_k,s]})^2\\
    &=\mathrm{I}_3(\mathbbm{1}_{[t_kt_{k+1}]}\otimes\mathbbm{1}_{[t_k,s]}^{\otimes 2})+2\E[B_{t_ks}^HB_{t_kt_{k+1}}^H]\mathrm{I}_1(\mathbbm{1}_{[t_k,s]})+|s-t_k|^{2H}\mathrm{I}_1(\mathbbm{1}_{[t_kt_{k+1}]}),
\end{align*}
where $\mathrm{I}_p$ denotes the multiple Skorokhod integral of order $p$. This yields a decomposition into the first and third Wiener chaoses, leading to the following rewriting of the expectation:
\begin{align*}
    \E&\left[\partial^2_{xx}b(\theta_0,\bar X_{t_k},\bar{\mu}_{t_k})B^H_{t_kt_{k+1}}\left(\int_{t_k}^{t_{k+1}} |B^H_{t_ks}|^2\,\mathrm{d}s-\frac{\alpha_n}{2}|B^H_{t_kt_{k+1}}|^2\right)\right]\\
    &=\E\left[\partial^2_{xx}b(\theta_0,\bar X_{t_k},\bar\mu_{t_k})\int_{t_k}^{t_{k+1}}\left(2\E[B_{t_ks}^HB_{t_kt_{k+1}}^H]\mathrm{I}_1(\mathbbm{1}_{[t_k,s]})+\left(|s-t_k|^{2H}-\frac{3}{2}\alpha_n^{2H}\right)\mathrm{I}_1(\mathbbm{1}_{[t_kt_{k+1}]})\right)\mathrm{d}s\right]\\
    &\quad+\E\left[\partial^2_{xx}b(\theta_0,\bar X_{t_k},\bar\mu_{t_k})\int_{t_k}^{t_{k+1}}\left(\mathrm{I}_3(\mathbbm{1}_{[t_kt_{k+1}]}\otimes\mathbbm{1}_{[t_k,s]}^{\otimes 2})-\frac{1}{2}\mathrm{I}_3(\mathbbm{1}_{[t_kt_{k+1}]}^{\otimes 3})\right)\mathrm{d}s\right].
\end{align*}
We apply Malliavin integration by parts to show that this expectation is suitably small, starting with the terms of the first Wiener chaos:
\begin{align*}
    \left|\E\left[\partial^2_{xx}b(\theta_0,\bar X_{t_k},\bar\mu_{t_k})\mathrm{I}_1(\mathbbm{1}_{[t_k,s]})\right]\right|
    &=\left|\E\left[\langle D\partial^2_{xx}b(\theta_0,\bar X_{t_k},\bar\mu_{t_k}),\mathbbm{1}_{[t_k,s]}\rangle_\mathcal{H}\right]\right|\\
    &\leq c \langle\mathbbm{1}_{[0,t_k]},\mathbbm{1}_{[t_k,s]}\rangle_\mathcal{H}\leq c |s-t_k|^{2H},
\end{align*}
where we used the fact that the Malliavin derivative is supported on $[0,t_k]$. Consequently, we obtain the following bound for the first chaos terms:
\begin{align*}
    \Bigg|\E\Bigg[&\partial^2_{xx}b(\theta_0,\bar X_{t_k},\bar\mu_{t_k})\int_{t_k}^{t_{k+1}}\bigg(2\E[B_{t_ks}^HB_{t_kt_{k+1}}^H]\mathrm{I}_1(\mathbbm{1}_{[t_k,s]})+\left(|s-t_k|^{2H}-\frac{3}{2}\alpha_n^{2H}\right)\mathrm{I}_1(\mathbbm{1}_{[t_kt_{k+1}]})\bigg)\,\mathrm{d}s\Bigg]\Bigg|\\
    &\leq c \int_{t_k}^{t_{k+1}} \left(|s-t_{k}|^{4H}+\frac{3}{4}\alpha_n^{4H}\right)\mathrm{d}s\leq c \alpha_n^{4H+1}.
\end{align*}
Next, we bound the terms belonging to the third Wiener chaos using a similar strategy:
\begin{align*}
  &  \left|\E\left[\partial^2_{xx}b(\theta_0,\bar X_{t_k},\bar\mu_{t_k})\mathrm{I}_3(\mathbbm{1}_{[t_kt_{k+1}]}\otimes\mathbbm{1}_{[t_k,s]}^{\otimes 2})\right]\right|
    =\left|\E\left[\langle D^3\partial^2_{xx}b(\theta_0,\bar X_{t_k},\bar\mu_{t_k}),\mathbbm{1}_{[t_kt_{k+1}]}\otimes\mathbbm{1}_{[t_k,s]}^{\otimes 2}\rangle_{\mathcal{H}^3}\right]\right|\\
    &\leq \left|\langle\mathbbm{1}^{\otimes 3}_{[0,t_k]},\mathbbm{1}_{[t_kt_{k+1}]}\otimes\mathbbm{1}_{[t_k,s]}^{\otimes 2}\rangle_{\mathcal{H}^3}\right|\leq c \alpha_n^{2H}|s-t_k|^{4H}.
\end{align*}
Substituting this estimate back into the expectation yields:
\begin{align*}
    \Bigg|\E\Bigg[&\partial^2_{xx}b(\theta_0,\bar X_{t_k},\bar\mu_{t_k})\int_{t_k}^{t_{k+1}}\left(\mathrm{I}_3(\mathbbm{1}_{[t_kt_{k+1}]}\otimes\mathbbm{1}_{[t_k,s]}^{\otimes 2})-\frac{1}{2}\mathrm{I}_3(\mathbbm{1}_{[t_kt_{k+1}]}^{\otimes 3})\right)\mathrm{d}s\Bigg]\Bigg|\\
    &\leq c\int_{t_k}^{t_{k+1}}\left(\alpha_n^{2H}|s-t_k|^{4H}+\frac{1}{2}\alpha_n^{6H}\right)\mathrm{d}s\leq c \alpha_n^{6H+1}.
\end{align*}
Finally, combining all the previous bounds, we obtain:
\begin{equation}\label{estimate precise sum remainder 4}
    \left|\E\left[\frac{1}{\alpha_n}\sum_{k=0}^{n-1}\tilde{r}^1_k B^{1,H}_{t_kt_{k+1}}\right]\right|\leq c\alpha_n^{4H-1}.
\end{equation}

\end{proof}

\subsection{Proof of the CLT}\label{section proofs CLT}

Since the estimator does not admit an explicit expression, we exploit
the fact that it is a minimizer of the squared contrast. The first-order
condition gives
\[
\partial_\theta(Q_n^N)^2(\theta_n^N)=0.
\]
A Taylor expansion around $\theta_{0,j}$ for $j\in\{1,...,p\}$ therefore yields
\begin{align}
-\partial_{\theta}(Q_n^N)^2(\theta_0)
&=
(\theta_{n,j}^N-\theta_{0,j})
\int_0^1
\partial^2_{\theta\theta}(Q_n^N)^2
\bigl(
s\theta_0+(1-s)\theta_n^N
\bigr)
\,\mathrm{d}s.
\label{Taylor gradient constrast function}
\end{align}

We first establish the triangular-array CLT used for the contrast
function.

\begin{proof}[Proof of Theorem~\ref{Linderberg CLT}]
 The proof is inspired by~\cite{proofclt} and relies on the following operator: for X a random variable we define:
        \begin{align*}
            T_X: &C_b^{unif}\rightarrow \R \\ 
            &f \rightarrow (x\mapsto \E[f(x+X)]
        \end{align*}
        where $C_b^{unif}$ is the space of bounded uniformly continuous real functions. 
        One can easily check that the operator is linear, continuous and that if $X_1,X_2$ are independent random variables that $T_{X_1+X_2}=T_{X_1}T_{X_2}$ and they commute.
        We are interested in such operator because it gives us a criterion for the convergence in distribution:
        \begin{lemma}
            Let $(X_n)_{n\geq 1}$ be a sequence of real random variables, if for all f such that it admits two derivatives which belongs to $C_b^{unif}$ we have for some random variable X:
            \[ \lim_{n\rightarrow \infty}\|T_{X_n}f-T_Xf\|_\infty=0\]
            then $(X_n)$ converges in distribution to X. 
        \end{lemma}
        The proof of the Lemma is in~\cite{proofclt}. We start by denoting $S_n^N=\frac{1}{\sqrt{N}}\sum_{i=1}^N Y_n^i$ and $Z\sim \mathcal{N}(0,1)$ using the fact that for $n$ fixed the $(Y_n^i)$ are i.i.d we have
        \begin{align*}
            &T_{S_n^N}=T_{\frac{Y_n^1}{\sqrt{N}}}...T_{\frac{Y_n^N}{\sqrt{N}}}=T_{\frac{Y_n^1}{\sqrt{N}}}^N,\\
            &T_{\sigma Z}=T^N_{\frac{\sigma Z}{\sqrt{N}}}.
        \end{align*}
        Using the fact that the operators are a contraction and that they commute with each other by independence we obtain: 
        \[\|T_{S_n^N}f-T_{\sigma Z}f\|_\infty \leq N\| T_{\frac{Y_n^i}{\sqrt{N}}}f-T_{\frac{\sigma Z}{\sqrt{N}}}f\|_\infty.\]
        Thus, we are going to focus on the right term and show that it converges to 0. We take $f \in C^2$ with both derivatives in $C^{unif}_b$. We start by a Taylor expansion of $f$, there exits $\eta$ between $x$ and $x+y$ such that:
        \[f(x+y)=f(x) +yf'(x)+\frac{y^2}{2}f^{}(x)+\frac{y^2}{2}(f''(\eta)-f''(x)).\]
        We take an $\varepsilon>0 $ since the second order derivative is uniformly continuous we have for all $\varepsilon>0$, there exists a $\delta>0$ such that for all $|x-y|\leq \delta$ we have $[f''(x)-f''(y)|\leq \varepsilon$.
        By applying the Taylor formula to $T_{\frac{Y^1_n}{\sqrt{N}}}f$ for all $x$ and using the fact that $Y_n^1$ is centered we obtain:
        \begin{align*}
            T_{\frac{Y^1_n}{\sqrt{N}}}f(x)&=f(x)+\frac{f'(x)}{\sqrt{N}}\E[Y_n^1]+\frac{\sigma_n^2}{2N}f''(x)+\frac{1}{2}\E[(\frac{Y_n^1}{\sqrt{N}})^2(f''(\eta)-f(x))]\\
            &=f(x)+\frac{\lambda_n^2}{N}f''(x)+\frac{1}{2}\E[(\frac{Y_n^1}{\sqrt{N}})^2(f''(\eta)-f''(x))\mathbf{1}_{|\eta-x|\leq \delta}]+\frac{1}{2}\E[(\frac{Y_n^1}{\sqrt{N}})^2(f''(\eta)-f''(x))\mathbf{1}_{|\eta-x|\geq \delta}].
        \end{align*}
        Similarly we get: 
        \begin{align*}
            T_{\frac{\sigma Z}{\sqrt{N}}}f(x)=f(x)+\frac{\lambda_n^2}{2N}f''(x)+\frac{1}{2N}\E[Z^2(f''(\eta)-f''(x))\mathbf{1}_{|\eta-x|\leq \delta}]+\frac{1}{2N}\E[Z^2(f''(\eta)-f''(x))\mathbf{1}_{|Z|\geq \sqrt{N}\delta}].
        \end{align*}
        We consider the difference and using the fact that $f^"$ is bounded and uniformly continuous, we get that:
        \begin{align*}
            |T_{\frac{Y_n^1}{\sqrt{N}}}f(x)-T_{\frac{\sigma Z}{\sqrt{N}}}f(x)|&\leq \frac{1}{2N}( \|f''\|_\infty |\lambda_n^2-\lambda^2|+\varepsilon(\lambda_n^2+\lambda^2)\\
            &+\|f''\|_\infty\E[|Y_n^1|\mathbf{1}_{|Y_n^1|^2\geq \sqrt{N}\delta}]+\|f''\|_\infty\E[|Z|^2\mathbf{1}_{|Z|\geq \sqrt{N}\delta}]).
        \end{align*}
        Hence, using the Lindeberg condition and the fact that the variance converges, we get that 
        \[\lim_{N,n\rightarrow\infty} N\|T_{\frac{Y_n^1}{\sqrt{N}}}f-T_{\frac{\sigma Z}{\sqrt{N}}}f\|_\infty=0.\]
\end{proof}

We now apply the previous result to the leading term in the expansion
of $Q_n^N(\theta_0)$.

\begin{proof}[Proof of Theorem~\ref{CLt contrast function}]
        We recall the computation in the proof of Proposition~\ref{CV contrast in theta 0}:
        \[
        Q_n^N(\theta_0)=\frac{\sigma^2}{N\alpha_n}\sum_{i=1}^N\sum_{k=0}^{n-1}\left(|B_{t_kt_{k+1}}^{i,H}|^{2}-\alpha_n^{2H}\right)+\frac{1}{N\alpha_n}\sum_{i=1}^N\sum_{k=0}^{n-1} \left(2\sigma B_{t_kt_{k+1}}^{i,H}r^{i,N}_k+|r^{i,N}_k|^2\right).\]
        First, we prove that the last two terms go to zero under the CLT normalization. We have different estimate depending on the case $H>\frac{1}{2}$ and $H\in(\frac{1}{4},\frac{1}{2})$. We start by the Young case and we use the estimate~\ref{estimate remainder}:
        \begin{align*}
            \E[\left|\frac{\sqrt{N\gamma_n^2}}{N\alpha_n}\sum_{i=1}^N \sum_{k=0}^{n-1} \left(2\sigma B_{t_kt_{k+1}}^{i,H}r^{i,N}_k+|r^{i,N}_k|^2\right)\right|^p]&\leq c (N\gamma^2_n)^{p/2} (\alpha_n^{p(2H-1)}+\alpha_n^{2pH})\leq c (N\alpha_n)^{p/2}.
        \end{align*}
        So using the Assumption~\ref{CV rate step size 2} we have the convergence to zero.\\
        For $H\in(1/3,1/2)$, we use the estimates~\eqref{estimate remainder},~\eqref{estimate precise sum remainder 1},~\eqref{estimate precise sum remainder 2},\eqref{estimate precise sum remainder 3} we have for all $\gamma\in(0,\frac{1}{2})$:
        \begin{align*}
            \E[\left|\frac{\sqrt{N\gamma_n^2}}{N\alpha_n}\sum_{i=1}^N \sum_{k=0}^{n-1} \left( 2\sigma B_{t_kt_{k+1}}^{i,H}r^{i,N}_k+|r^{i,N}_k|^2\right)\right|^p]&\leq c (Nn^{4H-3})^{p/2}(\frac{n ^{p(1-2H)}}{N^{p\gamma}}+\frac{n^{p(1-2H)}}{N^{p-1}}+\alpha_n^{p(3H-1)}+\alpha_n^{2pH})\\
            &\leq c (\alpha_nN^{1-2\gamma})^{p/2}+(\alpha_n^{1+2H}N)^{p/2}+(\alpha_n^3 N)^{p/2}.
        \end{align*}
        Thus using Assumption~\ref{CV rate step size 3} we have the convergence to zero.
        Finally, when $H\in(1/4,1/3]$, with the bounds~\eqref{estimate precise sum remainder 1},\eqref{estimate precise sum remainder 2} and~\eqref{estimate precise sum remainder 4}, we have for all $\gamma \in(0,\frac12)$:
        \begin{align*}
            \E[\left|\frac{\sqrt{N\gamma_n^2}}{N\alpha_n}\sum_{i=1}^N \sum_{k=0}^{n-1} \left( 2\sigma B_{t_kt_{k+1}}^{i,H}r^{i,N}_k+|r^{i,N}_k|^2\right)\right|^p]&\leq c (Nn^{4H-3})^{p/2}(\frac{n ^{p(1-2H)}}{N^{p\gamma}}+\frac{n^{p(1-2H)}}{N^{p-1}}+\alpha_n^{p(4H-1)}+\alpha_n^{2pH})\\
            &\leq c (\alpha_nN^{1-2\gamma})^{p/2}+(\alpha_n^{1+4H}N)^{p/2}+(\alpha_n^3 N)^{p/2}.
        \end{align*}
        Hence, using Assumption~\ref{CV rate step size 4}, we have the convergence towards zero.
        
We now focus on the first term. We apply the result just established to the sequence of random variables
$$Y_n^i := \frac{\sigma^2 \gamma_n}{\alpha_n} \sum_{k=0}^{n-1} \left(|B_{t_k t_{k+1}}^{i,H}|^2 - \alpha_n^{2H}\right).$$
These variables satisfy the i.i.d.\ assumption, as they are driven by independent fractional Brownian motions. Furthermore, the convergence of their variance is guaranteed by Proposition~\ref{L2_Estimate_fBm}, with the scaling sequence $\gamma_n$ appropriately chosen to ensure a non-degenerate limit.

To conclude, it remains to verify the Lindeberg condition. Since the sequence belongs to a fixed Wiener chaos, Gaussian hypercontractivity ensures the equivalence of moments, meaning that all higher-order moments are uniformly bounded in $n$. Let us fix a moment $2p > 2$, and let $p'$ denote its Hölder conjugate such that $1/p + 1/p' = 1$. By applying Hölder's inequality followed by Markov's inequality, we obtain for any $\delta > 0$:
\begin{align*}
    \E\Big[|Y_n^1|^2 \mathbf{1}_{|Y_n^1| \geq \sqrt{N}\delta}\Big] 
    &\leq \E\big[|Y_n^1|^{2p}\big]^{1/p} \P\big(|Y_n^1| \geq \sqrt{N} \delta\big)^{1/p'} \leq \sup_{k \in \N} \E\big[|Y_k^1|^{2p}\big]^{1/p} \left( \frac{\lambda_n}{\delta^2 N} \right)^{1/p'},
\end{align*}
where $\lambda_n := \E[|Y_n^1|^2]$ converges by Proposition~\ref{L2_Estimate_fBm}. As $N, n \rightarrow \infty$, the right-hand side converges to zero, confirming that the Lindeberg condition holds. Consequently, applying Theorem~\ref{Linderberg CLT} yields the desired convergence.
\end{proof}

\begin{proof}[Proof of Proposition~\ref{assymptotic first and second order derivative}]
We present the proof in the Young regime; the rough case follows by
the same argument, with the corresponding symmetric approximation. First,
\[
\partial^2_{\theta\theta}(Q_n^N)^2(\theta)=2\bigl(\partial_{\theta} Q_n^N(\theta)\bigr)^2+2Q_n^N(\theta)\partial^2_{\theta\theta}Q_n^N(\theta).
\]
The uniform convergence of $Q_n^N$ was established in the consistency proof. Differentiating the contrast gives, 
\begin{align*}
\partial_{\theta}Q_n^N(\theta)&=-\frac{2}{N}\sum_{i=1}^N\sum_{k=0}^{n-1}\partial_{\theta} b(\theta,X_{t_k}^{i,N},\mu_{t_k}^N)\left(X_{t_kt_{k+1}}^{i,N} -\alpha_nb(\theta,X_{t_k}^{i,N},\mu_{t_k}^N)\right).
\end{align*}
Substituting the expansion of the particle increments yields
\begin{align*}
\partial_{\theta}Q_n^N(\theta)&=\frac{2\alpha_n}{N}\sum_{i=1}^N\sum_{k=0}^{n-1}\partial_{\theta}b(\theta,X_{t_k}^{i,N},\mu_{t_k}^N) \left(b(\theta,X_{t_k}^{i,N},\mu_{t_k}^N)-b(\theta_0,X_{t_k}^{i,N},\mu_{t_k}^N)\right)\\
&\quad-\frac{2\sigma}{N}\sum_{j=1}^N\sum_{k=0}^{n-1}\partial_{\theta}b(\theta,X_{t_k}^{i,N},\mu_{t_k}^N) B_{t_kt_{k+1}}^{i,H}- \frac{2}{N}\sum_{j=1}^N\sum_{k=0}^{n-1}\partial_{\theta}b(\theta,X_{t_k}^{i,N},\mu_{t_k}^N)r_k^{i,N}.
\end{align*}

Similarly,
\begin{align*}
\partial_{\theta\theta}^2Q_n^N(\theta)&=\frac{2\alpha_n}{N}\sum_{i=1}^N\sum_{k=0}^{n-1}\left( \partial_{\theta}b(\theta,X_{t_k}^{i,N},\mu_{t_k}^N)\right)^2\\
&\quad+\frac{2\alpha_n}{N}\sum_{i=1}^N\sum_{k=0}^{n-1}\partial_{\theta}^2 b(\theta,X_{t_k}^{i,N},\mu_{t_k}^N)\left(b(\theta,X_{t_k}^{i,N},\mu_{t_k}^N)-b(\theta_0,X_{t_k}^{i,N},\mu_{t_k}^N)\right)\\
&\quad -\frac{2\sigma}{N}\sum_{i=1}^N\sum_{k=0}^{n-1}\partial_{\theta}^2b(\theta,X_{t_k}^{i,N},\mu_{t_k}^N) B_{t_kt_{k+1}}^{i,H}-\frac{2}{N}\sum_{i=1}^N\sum_{k=0}^{n-1}\partial_{\theta}^2b(\theta,X_{t_k}^{i,N},\mu_{t_k}^N)r_k^{i,N}.
\end{align*}

The terms without $r_k^{i,N}$ converge by Lemmas~\ref{CV somme Riemann classique} and~\ref{CV somme Riemann fBm, smooth case}, while the terms involving$ r_k^{j,N}$ vanish by~\eqref{estimate remainder}. This yields the claimed convergence of the first- and second-order derivatives.
\end{proof}

\begin{proof}[Proof of Theorem~\ref{CLT}]
Combining Equation~\eqref{Taylor gradient constrast function} with Theorem~\ref{CLt contrast function} and Proposition~\ref{assymptotic first and second order derivative} yields
$$
\sqrt{N\gamma_n^2} (\theta_{n}^N-\theta_{0}) \xrightarrow[N,n\to\infty]{\mathcal{D}} \mathcal{N} \left(0, \frac{c_H^2}{4} \left(\E\left[\int_0^T \partial_{\theta}b(\theta_0, \bar{X}_s, \bar{\mu}_s) \circ \mathrm{d}B_s^H\right]\right)^{-2}\right).
$$

\end{proof}

\subsection{Proof for the generalization with unknown diffusion and regularity}
\label{section proofs Unknown H and sigma}

Most of the arguments in this section follow the proofs for the case in which $H$ and $\sigma$ are  known. We therefore focus on the additional estimates required to handle the plug-in estimators.

We first study the estimators of $H$ and $\sigma$ based on a single trajectory.

\begin{proof}[Proof of Proposition~\ref{estimator one particle sigma and H}]
We prove that the empirical quadratic variations computed from the interacting trajectories has the same asymptotic as the associated to only $B^H$.
    Consider the following expansion:
    \begin{align*}
        S_nX^{i,N}(a^m)&=\frac{1}{n-3m}\sum_{k=3m}^{n-1}|V^{a^m}X^{i,N}(\frac{k}{n})|^2\\
        &=\frac{1}{n-3m}\sum_{k=3m}^{n-1}|V^{a^m}I^{i,N}(\frac{k}{n})|^2+\frac{\sigma^2}{n-3m}\sum_{k=ml}^{n-1}|V^{a^m}B^{i,H}(\frac{k}{n})|^2\\
        &+\frac{2 \sigma}{n-3m}\sum_{k=ml}^{n-1}V^{a^m}I^{i,N}(\frac{k}{n})V^aB^{i,H}(\frac{k}{n}),
    \end{align*}
    where $I^{i,N}(t)=X^i_0+\int_0^tb(\theta_0,X^{i,N}_s,\mu_s^N)ds$.
    Note that, using a Taylor expansion we get for all $\varepsilon>0$:
    \begin{align*}
    V^{a^m}I^{i,N}(\frac{k}{n})&=\sum_{q=0}^{3m}a^m_qI^{i,N}(\frac{(k-q)T}{n})\\
    &=\sum_{q=0}^{3m}a^m_q(I^{i,N}(\frac{kT}{n})-\frac{qT}{n}b(\theta_0,X^{i,N}_{t_k},\mu^N_{t_k})+\int_{\frac{(k-q)T}{n}}^\frac{kT}{n}b(\theta_0,X^{i,N}_{t_k},\mu^N_{t_k})-b(\theta_0,X^{i,N}_s,\mu^N_s)\mathrm{d}s)\\
    &=I^{i,N}(\frac{kT}{n})\sum_{q=0}^{3m}a^m_q-b(\theta_0,X^{i,N}_{t_k},\mu^N_{t_k})\sum_{q=0}^{3m}a^m_q \frac{qT}{n}+\int_{\frac{(k-q)T}{n}}^\frac{kT}{n}b(\theta_0,X^{i,N}_{t_k},\mu^N_{t_k})-b(\theta_0,X^{i,N}_s,\mu^N_s)\mathrm{d}s\\
    &=\sum_{q=0}^{3m}a^m_q\int_{\frac{(k-q)T}{n}}^\frac{kT}{n}b(\theta_0,X^{i,N}_{t_k},\mu^N_{t_k})-b(\theta_0,X^{i,N}_s,\mu^N_s)\mathrm{d}s =O(n^{-1-H+\varepsilon}),
    \end{align*}
    the linear term disappears since $a$ is a filter of order 3. Using the regularity of $B^H$ we have: $V^aB^{i,H}(\frac{i}{n})=O(n^{-H+\varepsilon})$. Hence, we obtain:
    \begin{equation}\label{eq: empirical second moment particle}
        S_nX^{i,N}(a^m) = \sigma^2 S_nB^{i,H}(a^m)+O(n^{-1-2H+\varepsilon}).
    \end{equation}
   We follow the proof of Proposition 4 in~\cite{estimationparameterfBm} and since the two empirical means are asymptotically equivalent and we apply continuous functions to obtain the estimators, we have the same asymptotics. 
\end{proof}

\begin{proof}[Proof of Theorem~\ref{average estimator of sigma and H}]

Although the particle-based estimators are dependent, we first
show that they admit an asymptotic representation in terms of
independent functionals of the driving fractional Brownian
motions, with a negligible remainder. We can then apply
Theorem~\ref{Linderberg CLT} to the independent leading terms. For each particle, the estimators are constructed as linear functionals of the logarithm of the empirical second moments. Using the expansion in~\eqref{eq: empirical second moment particle}, we obtain
\begin{align*}
    \log(S_n X^{i,N}(a^m)) &= \log\bigl(\sigma^2S_n B^{i,H}(a^m) + S_nX^{i,N}(a^m)-\sigma^2S_nB^{i,H}(a^m)\bigr) \\
    &= \log(\sigma^2S_n B^{i,H}(a^m)) + \frac{S_nX^{i,N}(a^m)-\sigma^2S_nB^{i,H}(a^m)}{\sigma^2S_n B^{i,H}(a^m)}+O\Big(\Big(\frac{S_nX^{i,N}(a^m)-\sigma^2S_nB^{i,H}(a^m)}{\sigma^2S_n B^{i,H}(a^m)}\Big)^2\Big).
\end{align*}
By the self-similarity of the fBm and the fact that $a$ is a filter, we have that
\[
\E[S_nB^{i,H}(a^m)]=\frac{1}{n-3m}\sum_{k=3m}^{n-1}\E[|V^{a^m}B^{i,H}(\frac {kT}{n})|^2]=-\frac{\alpha_n^{2H}}{2}\sum_{q,r=0}^{3m} a_q^ma_r^m |q-r|^{2H}=\alpha_n^{2H}\pi_H^{a^m}(0).
\]
Hence we can rewrite $S_n B^{i,H}(a^m)$ in the following way:
\begin{align*}
S_nB^{i,H}(a^m)&=\alpha_n^{2H}\pi_H^{a^m}(0)\left(1+\left(\frac{1}{n-3m}\sum_{k=3m}^{n-1}\frac{|V^{a^m}B^{i,H}(\frac{kT}{n})|^2}{\alpha_n^{2H}\pi_H^{a^m}(0)}-1\right)\right)\\
&=\alpha_n^{2H}\pi_H^{a^m}(0)\left(1+\left(\frac{1}{n-3m}\sum_{k=3m}^{n-1}\frac{|V^{a^m}B^{i,H}(\frac{kT}{n})|^2}{\E[|V^{a^m}B^{i,H}(\frac{kT}{n})|^2]}-1\right)\right).
\end{align*}
Thus using Proposition 3 of~\cite{estimationparameterfBm} we obtain:
\[
\alpha_n^{-2H}S_nB^{i,H}(a^m)=\pi_{H}^{a^m}(0)+O_\P(\sqrt{\alpha_n}).
\]

We can denote the linear operations applied to the logarithm by $f$. Summing over all $N$ particles, the aggregated term takes the form
\[
    \frac{1}{N}\sum_{i=1}^N f\bigl(\log(\sigma^2S_n B^{i,H}(a^m))\bigr) + O_\P(\alpha_n).
\]
The $\varepsilon$-loss in the previous pathwise estimate comes
from using the sample-path Hölder regularity of fractional
Brownian motion. Here, we instead work with $L^p$ estimates,
for which the increments satisfy the exact scaling
$\|B_t^H-B_s^H\|_{L^p}\lesssim |t-s|^H$. The corresponding
moment estimates therefore retain the exponent $H$, without
an $\varepsilon$-loss, yielding an $O_{\mathbb P}(\alpha_n)$
remainder in the logarithmic expansion. The dominant contribution arises from the fractional Brownian motions. Since these driving processes are independent, we achieve the standard $\sqrt{N}$ convergence rate. Furthermore, for $H>\frac{1}{2}$ under Assumption~\ref{CV rate step size 2}, the remainder term properly converges under the $\sqrt{Nn}$ normalization.

When $H<\frac12$, since we cannot use Assumption~\ref{CV rate step size 2}, we will have to detail the remainder. First, we have:
\[
S_nX^{i,N}(a^m)-\sigma^2S_nB^{i,H}\left(a^m\right)=\frac{2 \sigma}{n-3m}\sum_{k=3m}^{n-1}V^{a^m}I^{i,N}(\frac{kT}{n})V^{a^m}B^{i,H}(\frac{kT}{n})+O_\P(\alpha_n^{2H+2}).
\]
We focus on the first term.  As in the proof of the CLT  we start by using the propagation of chaos to obtain the same expression but with independent particles:
\begin{align*}
 &\E[\sqrt{Nn}|\frac{2}{n-3m}\sum_{k=3m}^{n-1}V^{a^m}I^{i,N}(\frac{kT}{n})V^{a^m}B^{i,H}(\frac{kT}{n})-\frac{2}{n-3m}\sum_{k=3m}^{n-1}V^{a^m}\bar I^i(\frac{kT}{n})V^{a^m}B^{i,H}(\frac{kT}{n})|^2]\\
 & \lesssim \sqrt{N^{1-2\gamma}\alpha_n^{1+4H}},   
\end{align*}
with $\bar I^i$ being the counterparts of $I^{i,N}$ for the McKean-Vlasov particles. We defined
\[
Z^i_n=\frac{1   }{n-3m}\sum_{k=3m}^{n-1}V^{a^m}\bar I^i(\frac{kT}{n})V^{a^m}B^{i,H}(\frac{kT}{n}).
\]
Moreover, we have by Rosenthal inequality~\eqref{Rosenthal inequality}:
\begin{align*}
    \E[|\frac{\sqrt{Nn}}{N}\sum_{i=1}^N\frac{Z^i_n}{S_nB^{i,H}(a^m)}-\E[\frac{Z^i_n}{S_nB^{i,H}(a^m)}]|^p]&\lesssim \alpha_n^{\frac p2}+\frac{\alpha_n^{\frac{p}{2}}}{N^{\frac p2-1}}.
\end{align*}
Thus, the only term that remains is the expectation. To prove that it goes fast enough to zero, we perform a Taylor expansion:
\begin{align*}
    Z^i_n&=\frac{1}{n-3m}\sum_{k=3m}^{n-1}V^{a^m}\bar I^i(\frac{kT}{n})V^{a^m}B^{i,H}(\frac{kT}{n})\\
    &=\frac{1}{n-3m}\sum_{k=3m}^{n-1}\sum_{q=0}^{3m}a_q^m\int_{\frac{(k-q)T}{n}}^\frac{kT}{n}b(\theta_0,\bar X^{i}_{t_k},\bar\mu_{t_k})-b(\theta_0,\bar X^{i}_s,\bar \mu_s)\mathrm{d}sV^{a^m}B^{i,H}(\frac{kT}{n})\\
    &=\frac{\sigma}{n-3m}\sum_{k=3m}^{n-1}\partial_xb(\theta_0,\bar X_{t_k},\bar\mu_{t_k})\sum_{q=0}^{3m}a^m_q\int_{\frac{(k-q)T}{n}}^{\frac{kT}{n}}B^{i,H}_{t_k s}\mathrm{d}sV^{a^m}B^{i,H}(\frac{kT}{n})\\
    &+\frac{\sigma}{n-3m}\sum_{k=3m}^{n-1}\sum_{q=0}^{3m} a_q^m \int_{\frac{(k-q)T}{n}}^{\frac{kT}{n}}\E'[\partial_\mu b(\theta_0,\bar X_{t_k}^i,\bar\mu_{t_k})(\bar X'_{t_ks})\times (B^{i,H}_{t_k})']\mathrm{d}sV^{a^m}B^{i,H}(\frac{kT}{n})+o(\alpha_n^{1+3H})
\end{align*}
with $\bar X'$ and $(B^{H})'$ being independent copies of $\bar X^i$ and $B^{i,H}$. We start  with the term involving the spatial derivative, we want to decompose the product of the fBm in the Wiener chaos, we first compute the following expectation:
\begin{align*}
    &\E[\sum_{q=0}^{3m}a_q^m\int_{\frac{(k-q)T}{n}}^{\frac{kT}{n}}B^{i,H}_{t_k s}\mathrm{d}sV^{a^m}B^{i,H}(\frac{kT}{n})]=\alpha_n\E[\sum_{q=0}^{3m}a_q^m\int_{0}^{q}B^{i,H}_{t_kt_k-\frac{(k-u)T}{n}}\mathrm{d}uV^{a^m}B^{i,H}(\frac{kT}{n})] \\
    &=\alpha_n\sum_{q=0}^{3m}\sum_{r=0}^{3m}a_q^m a_r^m\int_0^q \E[(B^{i,H}_{t_k-\frac{uT}{n}}-B^{i,H}_{t_k})B^{i,H}_{t_k-\frac{rT}{n}}]\mathrm{d}u\\
    &=-\frac{\alpha_n^{1+2H}}{2}\sum_{q=0}^{3m}\sum_{r=0}^{3m}a_q^m a_r^m\int_0^q|r-u|^{2H}\mathrm{d}u 
    \\
    &=-\frac{\alpha_n^{1+2H}}{2}\int_0^{3m}\sum_{q=0}^{3m}\sum_{r=0}^{3m}a_q^m a_r^m |r-3m+v|^{2H}(1-\mathbbm{1}_{v\leq 3m-q})\mathrm{d}v \\
    &=\frac{\alpha_n^{1+2H}}{2}\int_0^{3m}\sum_{q=0}^{3m}\sum_{r=0}^{3m}a_q^m a_r^m |r-3m+v|^{2H}\mathbbm{1}_{v\leq 3m-q}\mathrm{d}v\\
    &=\frac{\alpha_n^{1+2H}}{2}\int_0^{3m}\sum_{q'=0}^{3m}\sum_{r'=0}^{3m}a_{3m-q'}^m a_{3m-r'}^m |r'-v|^{2H}\mathbbm{1}_{v\leq q'}\mathrm{d}u\\
    &=\frac{\alpha_n^{1+2H}}{2}\int_0^{3m}\sum_{q'=0}^{3m}\sum_{r'=0}^{3m}a_{q'}^m a_{r'}^m |r'-v|^{2H}\mathbbm{1}_{v\leq q'}\mathrm{d}u
\end{align*}
where in the last equality we used that $a^m_q=-a^m_{3m-q}$. Thus we have proven:
\[
 \E[\sum_{q=0}^{3m}a_q^m\int_{\frac{(k-q)T}{n}}^{\frac{kT}{n}}B^{i,H}_{t_k s}\mathrm{d}sV^{a^m}B^{i,H}(\frac{kT}{n})]=0.
\]
We want to decompose this product onto the Wiener chaos, thus we denote:
\begin{align*}
    &\sum_{q=0}^{3m}a_q^m\int_{\frac{(k-q)T}{n}}^{\frac{kT}{n}}B^{i,H}_{t_k s}\mathrm{d}s=-I^1(\underbrace{\sum_{q=0}^{3m}a_q^m\int^{\frac{kT}{n}}_{\frac{(k-q)T}{n}}\mathbbm1_{[s,t_k]}\mathrm{d}s}_{:=\psi_n^k}),\\
    &V^{a^m}B^{i,H}(\frac{kT}{n})=\sum_{q=0}^{3m}a_q^mB^{i,H}_{t_k-\frac{qT}{n}}=I^1(\underbrace{\sum_{q=0}^{3m}a_q^m \mathbbm1_{[0,t_k-\frac{qT}{n}]}}_{:=\phi_n^k}).
\end{align*}
Now we use the product formula:
\begin{align*}
    \sum_{q=0}^{3m}a_q^m\int_{\frac{(k-q)T}{n}}^{\frac{kT}{n}}B^{i,H}_{t_k s}\mathrm{d}sV^{a^m}B^{i,H}(\frac{kT}{n})=I^2(\psi_n^k\otimes\phi_n^k).
\end{align*}
There is only one term since the product has zero expectation. Hence, by doing integration by parts we obtain
\begin{align*}
    |\E[\frac{\sigma}{n-3m}\sum_{k=3m}^{n-1}&\frac{\partial_xb(\theta_0,\bar X_{t_k},\bar\mu_{t_k})}{S_nB^{i,H}(a^m)} \sum_{q=0}^{3m}a^m_q\int_{\frac{(k-q)T}{n}}^{\frac{kT}{n}}B^{i,H}_{t_k s}\mathrm{d}sV^{a^m}B^{i,H}(\frac{kT}{n})]|\\
    &=|\frac{\sigma}{n-3m}\sum_{k=3m}^{n-1}\E[\frac{\partial_xb(\theta_0,\bar X_{t_k},\bar\mu_{t_k})}{S_nB^{i,H}(a^m)}I^2(\psi_n^k\otimes \phi_n^k)]| \lesssim\alpha_n^{1+2H}. 
\end{align*}
It now remains to handle the term with the Lions derivative. Applying integration by parts once again, we obtain:
\begin{align*}
    |\E[\frac{\sigma}{n-3m}\sum_{k=3m}^{n-1}\E'[&\frac{\partial_\mu b(\theta_0,\bar X_{t_k}^i,\bar\mu_{t_k})(\bar X'_{t_k})}{S_nB^{i,H}(a^m)}\times \sum_{q=0}^{3m} a_q^m \int_{\frac{(k-q)T}{n}}^{\frac{kT}{n}}(B^{i,H}_{t_k})'\mathrm{d}s]V^{a^m}B^{i,H}(\frac{kT}{n})]|\\
    &=|\frac{\sigma}{n-3m}\sum_{k=3m}^{n-1}\E[\E'[\frac{\partial_\mu b(\theta_0,\bar X_{t_k}^i,\bar\mu_{t_k})(\bar X'_{t_k})}{S_nB^{i,H}(a^m)}I^1(\psi_n^k)](I^1)'(\phi_n^k)]|
    \lesssim \alpha_n^{1+2H}
\end{align*}
with $(I^1)'$ being the Skorokhod integral given the noise $B'$. Combining the previous estimates with the Taylor remainder
bound, we obtain:
\[
|\sqrt{Nn}\E[\frac{Z_n^i}{S_nB^{i,H}(a^m)}]|\lesssim\sqrt{Nn}(\alpha_n^{1+2H}+\alpha_n^{1+H})\lesssim\sqrt{N\alpha_n^{1+2H}}.
\]
Thus for $H\in(\frac13,\frac12)$ it goes to zero using the Assumption~\ref{CV rate step size 3}. Hence we have the following expansion:
\[
\sqrt{Nn}\frac{1}{N}\sum_{i=1}^N f(\log(\sigma^2S_nB^{i,H}(a^m)))+O_\P(\sqrt{N^{1-2\gamma}\alpha_n}+\alpha_n^{\frac 12}+\alpha_n^{\frac 12}N^{-\frac 14}+\sqrt{N\alpha_n^{1+2H}}).
\]
For $H\in(\frac14,\frac13]$, we are going to extend the Taylor expansion in the computation of $Z_n^i$:
\begin{align*}
Z_n^i&=\frac{\sigma}{n-3m}\sum_{k=3m}^{n-1}\partial_xb(\theta_0,\bar X_{t_k},\bar\mu_{t_k})\sum_{q=0}^{3m}a^m_q\int_{\frac{(k-q)T}{n}}^{\frac{kT}{n}}B^{i,H}_{t_k s}\mathrm{d}sV^{a^m}B^{i,H}(\frac{kT}{n})\\
&+\frac{\sigma}{n-3m}\sum_{k=3m}^{n-1}\sum_{q=0}^{3m}a^m_q\int_{\frac{(k-q)T}{n}}^{\frac{kT}{n}}\E'[\partial_\mu b(\theta_0,\bar X_{t_k}^i,\bar\mu_{t_k})(\bar X'_{t_k})\times (B^{i,H}_{t_k})']\mathrm{d}sV^{a^m}B^{i,H}(\frac{kT}{n})\\
&+\frac{\sigma^2}{2(n-3m)}\sum_{k=3m}^{n-1}\partial^2_xb(\theta_0,\bar X_{t_k},\bar\mu_{t_k})\sum_{q=0}^{3m}a^m_q\int_{\frac{(k-q)T}{n}}^{\frac{kT}{n}}|B^{i,H}_{t_k s}|^2\mathrm{d}sV^{a^m}B^{i,H}(\frac{kT}{n})\\
&+\frac{\sigma^2}{2(n-3m)}\sum_{k=3m}^{n-1}\sum_{q=0}^{3m}a^m_q\int_{\frac{(k-q)T}{n}}^{\frac{kT}{n}}\E'[\partial_y\partial_\mu b(\theta_0,\bar X_{t_k}^i,\bar\mu_{t_k})(\bar X'_{t_k})|(B^{i,H}_{t_ks})'|^2]\mathrm{d}sV^{a^m}B^{i,H}(\frac{kT}{n})\\
&+\frac{\sigma^2}{2(n-3m)}\sum_{k=3m}^{n-1}\sum_{q=0}^{3m}a^m_q\int_{\frac{(k-q)T}{n}}^{\frac{kT}{n}}\E'[\E^{''}[\partial^2_\mu b(\theta_0,\bar X_{t_k}^i,\bar\mu_{t_k})(\bar X'_{t_k},\bar X^{''}_{t_k}) (B^{i,H}_{t_ks})'(B^{i,H}_{t_ks})^{''}]]\mathrm{d}sV^{a^m}B^{i,H}(\frac{kT}{n})\\
&+o(\alpha_n^{1+4H}).
\end{align*}
We have already bounded the first two terms. For the third term (involving the second-order spatial derivative), we follow the same strategy as in the proof of the convergence of $Q_n^N(\theta_0)$, while for the last two terms we apply integration by parts to obtain
\[
|\sqrt{Nn}\E[\frac{Z_n^i}{S_nB^{i,H}(a^m)}]|\lesssim \sqrt{N\alpha_n^{1+4H}}.
\]
Hence, it converges towards zero because of Assumption~\ref{CV rate step size 4}. Thus, we obtain the following expression:
\[
\sqrt{Nn}\frac{1}{N}\sum_{i=1}^N f(\log(\sigma^2S_nB^{i,H}(a^m)))+O_\P(\sqrt{N^{1-2\gamma}\alpha_n}+\alpha_n^{\frac 12}+\alpha_n^{\frac 12}N^{-\frac 14}+\sqrt{N\alpha_n^{1+4H}}).
\]
We have proven that for all $H>\frac14$, under the related assumptions that the dominant term is the same as in the pure fBm case, which yields the same asymptotics for the estimators.

\end{proof}

We next compare the contrast with the plug-in version obtained by replacing $H$ and $\sigma$ with their estimators.

\begin{proposition}\label{contrast function sigma,H}
For $H\in(\frac14,1)$, $H\neq\frac12$, under Assumptions~\ref{Bdd moment},~\ref{Lip cond}, and~\ref{CV rate step size 1},
\[
\sup_{\theta\in\Theta}\left|Q_n^N(\theta)-\widetilde Q_n^N(\theta)\right|\xrightarrow[N,n\to\infty]{}0
\]
almost surely.
\end{proposition}

\begin{proof}[Proof of Proposition~\ref{contrast function sigma,H}]
In both the rough and the Young cases, the only difference between the contrast functions is the correction term, and therefore we get:
    \begin{align*}
        \sup_{\theta \in \Theta} |Q_n^N(\theta)-\tilde Q_n^N(\theta)|&=\left|\sigma^2\alpha_n^{2H-2}-(\sigma_n^N)^2 \alpha_n^{2H^N_n-2}\right|\\
        &\leq \sigma^2\alpha_n^{-2}|\alpha_n^{2H}-\alpha_n^{2H^N_n}|+\alpha_n^{2H^N_n-2}|\sigma_n^N+\sigma||\sigma_n^N-\sigma|\\
        &\leq c (\alpha_n^{2(H_n^N\wedge H)-2}\log(n)|H_n^N-H|+\alpha_n^{2H_n^N-2}|\sigma_n^N-\sigma|)\\
        &\leq c \frac{\alpha_n^{2(H_n^N\wedge H)-2}\log(n)}{\sqrt{Nn}}=c\left(\frac{n^{3-4(H\wedge H_n^N)}\log(n)^2}{N}\right)^{\frac{1}{2}}
    \end{align*}
    which goes to zero under Assumption~\ref{CV rate step size 1}.
\end{proof}

\begin{proof}[Proof of Theorem~\ref{consistency joint estimation}]
The proof follows the same argument as the proof of Theorem~\ref{Consistency Young case}. Proposition~\ref{contrast function sigma,H} shows that the plug-in contrast is asymptotically equivalent to the contrast with known $H$ and $\sigma$. We can therefore apply Proposition~\ref{CV contrast in theta 0} and Theorem~\ref{Diff CV} in the Young regime, and Theorem~\ref{CV contrast rough case} in the rough regime.
\end{proof}

It remains to prove the central limit theorem stated in Theorem~\ref{CLT joint estimation}.

\begin{proof}[Proof of Theorem~\ref{CLT joint estimation}]

The overall strategy remains the same as for the consistency result established in Theorem~\ref{consistency joint estimation}: we study the difference between the contrast functions evaluated with the true parameters $(H, \sigma)$ and the estimated parameters. Here, however, we consider the CLT normalization. Following the same computations as in the proof of Proposition~\ref{contrast function sigma,H}, we obtain the following expansion:
\begin{align*}
    & \sqrt{N\gamma_n^2}\tilde Q_n^N(\theta_0) = \sqrt{N\gamma_n^2}Q_n^N(\theta_0) + \sqrt{N\gamma_n^2}\bigl(\tilde Q_n^N(\theta_0) - Q_n^N(\theta_0)\bigr) \\
    &= \sqrt{N\gamma_n^2}Q_n^N(\theta_0) - \sqrt{N\gamma_n^2} \, 2T\sigma^2 \alpha_n^{2H-2} \biggl( \log(\alpha_n)(H_n^N-H) + \frac{1}{\sigma}(\sigma_n^N-\sigma) \biggr)  + \sqrt{N\gamma_n^2} \, O_{\P}\biggl( \frac{\log(n)^2}{Nn} \biggr).
\end{align*}

The asymptotic behavior now depends on the value of $H$. If $H \geq \frac{3}{4}$, we have $\gamma_n\alpha_n^{2H-2} = O(\sqrt{n})$, meaning the dominant contribution arises purely from the first term. Conversely, when $H < \frac{3}{4}$, the asymptotic fluctuations of $H_n^N$ and $\sigma_n^N$ must be explicitly taken into account. Indeed, utilizing the algebraic relations between the estimators detailed in the proof of Proposition 4 in~\cite{estimationparameterfBm}, we deduce that
\[
    \sqrt{Nn} \left( \log(\alpha_n)(H_n^N-H) + \frac{1}{\sigma}(\sigma_n^N-\sigma) \right) \xrightarrow{\mathcal{D}} \mathcal{N}\Bigl(0, \bigl(\log^2(T) + 2^{2H+1}\log(2)(4-2^{2H})^{-1}\bigr)\Sigma^2\Bigr).
\]
The leading terms in the linearized expansions belong to the
second Wiener chaos. Their marginal Gaussian convergence,
together with the convergence of their covariance, yields the
joint convergence. Furthermore, exploiting the self-similarity of the fractional Brownian motion increments, one can show that the covariance between the two terms converges to a constant $\rho_H$ given by
\[
    \rho_H = \frac{1}{\|A\|^2} \sum_{m=1}^2 W_{H,m} \sum_{r\in\mathbb{N}} \E\bigl[ B^H_1 V^{a^m} B^H(1+r) \bigr],
\]
where $W_{H,m} = 2A'_m - 2^{2H+1}\log(2)(4-2^{2H})^{-1}A_m$, following the notation in~\cite{estimationparameterfBm}. Consequently, we obtain the weak convergence:
\[
    \sqrt{N\gamma_n^2}\tilde Q_n^N(\theta_0) \xrightarrow[N,n\rightarrow+\infty]{\mathcal{D}} \mathcal{N} \Bigl(0, \underbrace{c_H^2 + T^2\sigma^4\bigl(\log^2(T) + 2^{2H+1}\log(2)(4-2^{2H})^{-1}\bigr)\Sigma^2 + 2\rho_H}_{:= \Lambda^2_H} \Bigr).
\]

Finally, proceeding exactly as in the proof of Theorem~\ref{CLT}, we establish the CLT for the estimator by applying the Taylor expansion in~\eqref{Taylor constrast function}.

\end{proof}

\appendix
\section{Appendix}
\subsection{Auxiliary results}\label{section preliminary results}
The following lemma is a direct consequence of the Borel--Cantelli lemma and allows us to convert sufficiently fast $L^p$ convergence into almost sure convergence.

\begin{lemma}\label{fast Lp implies as}
Let $p_0\in\N^*$. Let $(Y_n)_{n\in\N}$ be a sequence of random variables such that, for every $p\geq p_0$, $\E[|Y_n|^p] \leq c_p n^{-\alpha p}.$
Then, for every $\varepsilon>0$, there exists a random variable $\eta_\varepsilon\in L^p$ for every $p\geq p_0$ such that $|Y_n|\leq \eta_\varepsilon n^{-\alpha+\varepsilon}$
almost surely.
\end{lemma}
The proof is a direct application of the Borel-Cantelli lemma; see~\cite{kloeden2007pathwise}. We will also use the following array version of the law of large numbers.

\begin{lemma}\label{Array LLN}
Let $(Y_n^i)_{n,i\geq1}$ be real-valued random variables such that, for every $n$, $(Y_n^i)_{i\geq1}$ are i.i.d. Suppose that there exists an i.i.d. sequence $(Y^i)_{i\geq1}$ and $p>2$ such that $\sum_{n\geq0}\E[|Y_n^i-Y^i|^p]<+\infty.$ Then,
\[
\frac1N\sum_{i=1}^NY_n^i\xrightarrow[N,n\to\infty]{}\E[Y^1]
\]
almost surely.
\end{lemma}
The proof is again a consequence of the Borel-Cantelli lemma. We next record an interpolation inequality that is used repeatedly to pass from uniform estimates to Hölder estimates.

\begin{lemma}\label{interpolation inequality}
Let $0<\alpha<\beta$ and $f\in C^\beta([0,T])$. Then $f\in C^\alpha([0,T])$ and $\|f\|_\alpha\leq 2\|f\|_\infty^{1-\frac{\alpha}{\beta}}\|f\|_\beta^{\frac{\alpha}{\beta}}.$
\end{lemma}

\begin{proof}[Proof of Lemma~\ref{interpolation inequality}]
Set
\[
h=\left(\frac{\|f\|_\infty}{\|f\|_\beta}\right)^{1/\beta}.
\]
If $|x-y|\geq h$, then
\[
|f(x)-f(y)|\leq 2\|f\|_\infty\leq 2\|f\|_\infty^{1-\frac{\alpha}{\beta}}
\|f\|_\beta^{\frac{\alpha}{\beta}}
|x-y|^\alpha.
\]
If $|x-y|\leq h$, then
\[
|f(x)-f(y)|\leq \|f\|_\beta|x-y|^\beta
\leq \|f\|_\infty^{1-\frac{\alpha}{\beta}} \|f\|_\beta^{\frac{\alpha}{\beta}} |x-y|^\alpha.
\]
Combining the two estimates yields the result.
\end{proof}

We will also use the following standard estimate for Young integrals. Being a standard result we omit its proof.

\begin{lemma}\label{estimates Young integral}
Let $f\in C^\alpha([a,b])$ and $g\in C^\beta([a,b])$ with $\alpha+\beta>1$. Then there exists a constant $c_{\alpha,\beta}>0$ such that
\[
\left| \int_a^b (f(t)-f(a))\,\mathrm{d}g(t)\right|\leq c_{\alpha,\beta}\|f\|_\alpha\|g\|_\beta|b-a|^{\alpha+\beta}.
\]
\end{lemma}

We finally record the second-moment asymptotics for increments of fractional Brownian motion.

\begin{lemma}\label{L2_Estimate_fBm}
Let $B^H$ be a fractional Brownian motion with Hurst parameter $H\in(0,1)$. Then, for $H<\frac34$,
\[
\lim_{N,n\to\infty} \frac{N}{\alpha_n^{4H-3}} \E\left[ \left|\frac{1}{N\alpha_n}\sum_{i=1}^N\sum_{k=0}^{n-1}\left(|B_{t_kt_{k+1}}^{i,H}|^2-\alpha_n^{2H}\right)\right|^2\right]= c_H.
\]
For $H=\frac34$,
\[
\lim_{N,n\to\infty}\frac{N}{\log(n)}\E\left[\left|\frac{1}{N\alpha_n}\sum_{i=1}^N\sum_{k=0}^{n-1}\left(|B_{t_kt_{k+1}}^{i,H}|^2-\alpha_n^{2H}\right)\right|^2\right]=c_H.
\]
For $\frac34<H<1$,
\[
\lim_{N,n\to\infty}N\E\left[\left|\frac{1}{N\alpha_n}\sum_{i=1}^N\sum_{k=0}^{n-1}\left(|B_{t_kt_{k+1}}^{i,H}|^2-\alpha_n^{2H}\right)\right|^2\right]=c_H.
\]
\end{lemma}

\begin{proof}
Using the independence of the fractional Brownian motions,
\begin{align*}
\E\left[\left|\frac{1}{N\alpha_n}\sum_{i=1}^N\sum_{k=0}^{n-1}\left(|B_{t_kt_{k+1}}^{i,H}|^2-\alpha_n^{2H}\right)\right|^2\right]&=\frac{1}{N\alpha_n^2}\E\left[\left|\sum_{k=0}^{n-1}\left(|B_{t_kt_{k+1}}^H|^2-\alpha_n^{2H}\right)\right|^2\right]\\
&=\frac{\alpha_n^{4H-2}}{N}\E\left[\left|\sum_{k=0}^{n-1}\left(|B_{kk+1}^H|^2-\alpha_n^{2H}\right)\right|^2\right],
\end{align*}
where we used the self-similarity of fractional Brownian motion. Define $\gamma(i)=|i+1|^{2H}+|i-1|^{2H}-2|i|^{2H}.$
A direct computation gives
\begin{align*}
& \E\left[\left|\sum_{k=0}^{n-1}|B_{k+1}^H-B_k^H|^2-n\right|^2\right]=\operatorname{Var}\left(\sum_{k=0}^{n-1}|B_{k+1}^H-B_k^H|^2\right)\\
&=2\sum_{k=0}^{n-1}\sum_{\ell=0}^{n-1} \operatorname{Cov}\left(B_{k+1}^H-B_k^H,B_{\ell+1}^H-B_\ell^H\right)^2 =2\sum_{k=0}^{n-1}\sum_{\ell=0}^{n-1}\gamma(k-\ell)^2.
\end{align*}
Equivalently,
\[
\E\left[\left|\sum_{k=0}^{n-1}\left(|B_{k+1}^H-B_k^H|^2-1\right)\right|^2\right]=2\sum_{i=0}^{n-1}(n-i)\gamma(i)^2=:S_n.
\]
A Taylor expansion gives $\gamma(i)=2H(2H-1)i^{2H-2}+o(i^{2H-2}),$
and hence $\gamma(i)^2=4H^2(2H-1)^2i^{4H-4}+o(i^{4H-4}).$
When $H<\frac34$, the series
$\sum_i\gamma(i)^2$ converges, and therefore $S_n\sim c_H n.$
When $H=\frac34$, $S_n\sim c_H n\log n.$
Finally, when $\frac34<H<1$, $S_n \sim c_H n^{4H-2}$.
Substituting these asymptotics into the previous rescaling yields the three claims.
\end{proof}

We also need uniform moment and regularity estimates for the particle and McKean-Vlasov systems. Observe that the first two points can also be found in~\cite{amorino2025fractional}, we prove them again here as our model is no longer linear in the parameter.

\begin{proposition}\label{prop solution}
Under Assumptions~\ref{Bdd moment} and~\ref{Lip cond}, for every $N\in\N$, $i\in\{1,\ldots,N\}$, $p\geq2$, and $\varepsilon>0$ such that $H-\varepsilon>0$, there exists a constant $c>0$, independent of $N$, such that
\begin{enumerate}
    \item $\sup_{t\in[0,T]}
(    \E[|X_t^{i,N}|^p]+\E[|\bar X_t^i|^p]+\E[W_2^p(\mu_t^N,\delta_0)]+     \E[W_2^p(\bar\mu_t,\delta_0)])\leq c$;
    \item $\E[|X_t^{i,N}-X_s^{i,N}|^p]+\E\left[\left|W_2(\mu_t^N,\delta_0)-W_2(\mu_s^N,\delta_0)    \right|^p\right]\leq c|t-s|^{pH};$    
    \item $\E[\|X^{i,N}\|_{H-\varepsilon,[0,T]}^q]+\E[\|W_2(\mu^N,\delta_0)\|_{H-\varepsilon,[0,T]}^q]\leq c;$
    
    \item $ \E[|\bar X_t^i-\bar X_s^i|^p]+\E\left[\left|W_2(\bar\mu_t,\delta_0)-W_2(\bar\mu_s,\delta_0)    \right|^p\right]\leq c|t-s|^{pH};$

    \item $\E[\|\bar X^i\|_{H-\varepsilon,[0,T]}^q]+\E[\|W_2(\bar\mu,\delta_0)\|_{H-\varepsilon,[0,T]}^q]\leq c.$
\end{enumerate}
\end{proposition}

\begin{proof}
The proof relies on Grönwall's lemma and the moment estimates for fractional Brownian motion. We give the argument for the interacting particle system; the McKean-Vlasov case follows similarly. By the Lipschitz condition on the drift,
\begin{align*}
\E\left[ \sup_{t\leq T}|X_t^{i,N}|^p \right]&\leq c\E\left[|X_0^i|^p+\int_0^T\sup_{r\leq s}|b(\theta_0,X_r^{i,N},\mu_r^N)|^p\,\mathrm{d}s+\sigma^p\sup_{t\leq T}|B_t^{i,H}|^p\right]\\
&\leq c\left(1+\E[|X_0^i|^p]+\E\left[\sup_{t\leq T}|B_t^{i,H}|^p\right]\right)+c\int_0^T\E\left[\sup_{r\leq s}|X_r^{i,N}|^p\right]\,\mathrm{d}s,
\end{align*}
where we have used the fact that the Wasserstein distance from $\delta_0$ is controlled by the empirical moments of the particles. Grönwall's lemma yields the first estimate.
For the temporal regularity,
\begin{align*}
\E[|X_t^{i,N}-X_s^{i,N}|^p]&\leq c\E\left[\left|\int_s^tb(\theta_0,X_r^{i,N},\mu_r^N)\,\mathrm{d}r\right|^p\right]+c\E[|B_t^{i,H}-B_s^{i,H}|^p]\\
&\leq c|t-s|^p+c|t-s|^{pH}\leq c|t-s|^{pH}.
\end{align*}
The corresponding Wasserstein estimate follows from the definition of the Wasserstein distance.

The Hölder estimates then follow from Kolmogorov's continuity criterion, together with the interpolation inequality in Lemma~\ref{interpolation inequality}.
\end{proof}

Finally, we use the following form of Rosenthal's inequality.

\begin{proposition}\label{Rosenthal inequality}
Let $(Y_i)_{i\geq1}$ be a sequence of centered independent random variables. Then, for every $p>2$, there exists $C>0$ such that
\[
\E\left[\left|\sum_{i=1}^N Y_i\right|^p\right]\leq C\left(\left(\sum_{i=1}^N\E\left[|Y_i|^2\right]\right)^{p/2}+\sum_{i=1}^N\E[|Y_i|^p]\right).
\]
\end{proposition}

The proof can be found in \cite{Rosenthalinequality}.

 \vspace{0.5 cm}

\subsection{Elements of Malliavin Calculus for Fractional Brownian Motion}\label{app:malliavin_fbm}
We first recall the Mandelbrot--Van Ness representation of fractional Brownian motion (see~\cite{MandelbrotVannessRepresentation}). A process $(B_t^H)_{t\geq 0}$ is a centered fractional Brownian motion with Hurst parameter $H\in(0,1)$ if there exists a two-sided Brownian motion $(W_t)_{t\in\R}$ such that
\[
B_t^H=\int_{\R} K(t,s)\,\mathrm{d}W_s,
\]
where
\[
K(t,s)=c_H[(t-s)_+^{H-\frac12}-(-s)_+^{H-\frac12}].
\]

Next we introduce the derivative operator and its adjoint, the divergence operator. Let $\mathcal{E}$ be the set of step functions on $[0,T]$. The canonical Hilbert space $\mathcal{H}$ associated with $B^H$ is the closure of $\mathcal{E}$ with respect to the scalar product generated by linearly extending
\begin{equation*}
\langle \mathbf{1}_{[0,t]}, \mathbf{1}_{[0,s]} \rangle_{\mathcal{H}} = R_H(t,s).
\end{equation*}
When $H > \frac{1}{2}$, for any $f, g \in |\mathcal{H}|$ (where $|\mathcal{H}|$ denotes the Banach space of measurable functions such that $\|f\|_{|\mathcal{H}|}^2 := \int_0^T \int_0^T |f(t)||f(s)|\phi(t,s)\,\mathrm{d}s\,\mathrm{d}t < \infty$), the scalar product takes the explicit integral form:
\begin{equation*}
\langle f, g \rangle_{\mathcal{H}} = \int_0^T \int_0^T f(t) g(s) \phi(t,s) \, \mathrm{d}s \, \mathrm{d}t, \quad \text{where } \phi(t,s) := H(2H-1)|t-s|^{2H-2}.
\end{equation*}
The mapping $\mathbf{1}_{[0,t]} \mapsto B^H_t$ extends to a linear isometry $h \mapsto B^H(h)$ between $\mathcal{H}$ and the Gaussian subspace of $L^2(\Omega)$ generated by $B^H$, yielding an isonormal Gaussian process $\{B^H(h) : h \in \mathcal{H}\}$.

Let $\mathcal{S}$ be the class of smooth and cylindrical random variables of the form $F = f(B^H(h_1), \dots, B^H(h_n))$, where $f \in C_p^\infty(\mathbb{R}^n)$ and $h_i \in \mathcal{H}$. The Malliavin derivative of $F$ is the $\mathcal{H}$-valued random variable given by:
\begin{equation*}
DF = \sum_{i=1}^n \frac{\partial f}{\partial x_i}(B^H(h_1), \dots, B^H(h_n)) h_i.
\end{equation*}
We denote by $\mathbb{D}^{1,2}$ the closure of $\mathcal{S}$ with respect to the norm $\|F\|_{1,2}^2 = \mathbb{E}[|F|^2] + \mathbb{E}[\|DF\|_{\mathcal{H}}^2]$.

The divergence operator $\delta$, or Skorokhod integral with respect to $B^H$, denoted $\delta(u) = \int_0^T u_t \, \delta B^H_t$, is the adjoint of $D$. Its domain $\mathrm{Dom}(\delta)$ consists of all $u \in L^2(\Omega; \mathcal{H})$ such that $|\mathbb{E}[\langle DF, u \rangle_{\mathcal{H}}]| \le c_u \|F\|_{L^2(\Omega)}$ for all $F \in \mathbb{D}^{1,2}$.

\begin{lemma}[Integration by Parts (IBP) and Duality]
For any $F \in \mathbb{D}^{1,2}$ and $u \in \mathrm{Dom}(\delta)$, the duality relation holds:
\begin{equation*}
\mathbb{E}\big[ F \delta(u) \big] = \mathbb{E}\big[ \langle DF, u \rangle_{\mathcal{H}} \big].
\end{equation*}
Moreover, if $F \in \mathbb{D}^{1,2}$, $u \in \mathrm{Dom}(\delta)$, and $Fu \in L^2(\Omega; \mathcal{H})$ such that $F \delta(u) - \langle DF, u \rangle_{\mathcal{H}} \in L^2(\Omega)$, then $Fu \in \mathrm{Dom}(\delta)$ and the integration by parts formula holds:
\begin{equation*}
\delta(Fu) = F \delta(u) - \langle DF, u \rangle_{\mathcal{H}}.
\end{equation*}
\end{lemma}

The space $L^2(\Omega)$ decomposes into the orthogonal direct sum of Wiener chaoses $L^2(\Omega) = \bigoplus_{p=0}^{\infty} \mathcal{H}_p$. For $p \ge 1$, let $I_p : \mathcal{H}^{\odot p} \to \mathcal{H}_p$ denote the multiple Wiener-Itô integral of order $p$ with respect to $B^H$, where $\mathcal{H}^{\odot p}$ is the symmetric tensor power equipped with the norm $\sqrt{p!} \, \|\cdot\|_{\mathcal{H}^{\otimes p}}$. By convention, $I_0(c) = c$ for $c \in \mathbb{R}$.

\begin{lemma}[Product Formula]
Let $f \in \mathcal{H}^{\odot p}$ and $g \in \mathcal{H}^{\odot q}$. The product of two multiple integrals expands as:
\begin{equation*}
I_p(f) I_q(g) = \sum_{r=0}^{p \wedge q} r! \binom{p}{r} \binom{q}{r} I_{p+q-2r}\big(f \tilde{\otimes}_r g\big),
\end{equation*}
where $f \otimes_r g \in \mathcal{H}^{\otimes (p+q-2r)}$ is the contraction of order $r$ defined in $\mathcal{H}^{\otimes r}$, and $f \tilde{\otimes}_r g$ denotes its symmetrization. For $H > \frac{1}{2}$ and $f \in |\mathcal{H}|^{\odot p}$, $g \in |\mathcal{H}|^{\odot q}$, the $r$-th contraction is given explicitly by:
\begin{align*}
(f \otimes_r g)(t_1, \dots, t_{p+q-2r}) &= \int_{[0,T]^{2r}} f(t_1, \dots, t_{p-r}, u_1, \dots, u_r) \\
&\quad \times g(t_{p-r+1}, \dots, t_{p+q-2r}, v_1, \dots, v_r) \prod_{i=1}^r \phi(u_i, v_i) \, \mathrm{d}u_i \, \mathrm{d}v_i.
\end{align*}
\end{lemma}

\begin{lemma}[Young-Skorokhod Conversion Formula]
Assume $H \in (1/2, 1)$ and let $u = (u_t)_{t \in [0,T]}$ be a stochastic process in $\mathbb{D}^{1,2}(|\mathcal{H}|)$ such that
\begin{equation*}
\mathbb{E} \left[ \int_0^T \int_0^T |D_s u_t| \, \phi(t,s) \, \mathrm{d}s \, \mathrm{d}t \right] < \infty \quad \text{and} \quad \int_0^T \int_0^T |D_s u_t| \, \phi(t,s) \, \mathrm{d}s \, \mathrm{d}t < \infty \text{ a.s.}.
\end{equation*}
Then $u \in \mathrm{Dom}(\delta)$, the Young integral $\int_0^T u_t \circ \mathrm{d}B^H_t$ exists, and we have the conversion formula:
\begin{equation*}
\int_0^T u_t \circ \mathrm{d}B^H_t = \int_0^T u_t \, \delta B^H_t + \int_0^T \int_0^T D_s u_t \, \phi(t,s) \, \mathrm{d}s \, \mathrm{d}t,
\end{equation*}
where $\phi(t,s) = H(2H-1)|t-s|^{2H-2}$.
\end{lemma}

A direct consequence of this lemma is the identity
\begin{align}\label{eq: identity-means}
  \E \left[\int_0^T u_t \circ \mathrm{d}B^H_t\right] = \E \left[\int_0^T \int_0^T D_s u_t \, \phi(t,s) \, \mathrm{d}s \, \mathrm{d}t\right],  
\end{align}
since the Skorokhod integral has zero mean by definition.

For complete proofs and further details, a standard reference is~\cite{MalliavinNualart}.

\begin{lemma}\label{lemma malliavin derivative McKean particle} Let us denote by $D^i$ the Malliavin derivative with respect to $B^{H,i}$. 
    Under Assumptions~\ref{Bdd moment},\ref{Lip cond},\ref{strong assumption on derivative} applied to $b$, we have:
    \begin{enumerate}
        \item for all $t,s\in[0,T]$ and $i\in\N$: $D^i_s\bar X^i_t=(\sigma e^{\int_s^t\partial_xb(\theta_0,\bar X_s,\bar\mu_s)ds})\mathbbm{1}_{s<t}$,
        \item for all $t,s\in[0,T]$and $i,j \in \N, i\neq j$: $D_s^j\bar X^i_t=0$.
    \end{enumerate}
\end{lemma}
The proof is given in~\cite{amorino2025fractional}.

\begin{lemma}\label{lemma bound DW}[Singular Estimate for $D^W X_t$ when $H < 1/2$]
Let $H \in (0, 1/2)$, we denote $W^i$ the underlying Brownian motion of $B^{H,i}$. There exists an almost surely finite random variable $M \in \bigcap_{p \ge 1} L^p(\Omega)$ such that for all $0 \le s < t \le T$ and each component $i \geq1$:
    \begin{equation*}
    \big| D_s^{W^i}\bar X_t^i \big| \le M \, (t - s)^{H - \frac{1}{2}}.
    \end{equation*}
\end{lemma}
The proof is given in Section 5.3 of~\cite{MalliavinNualart}.

\subsection{Lions Derivatives on the Space of Probability Measures}\label{app: Lions derivative}
 We consider a rich, atomless probability space $(\Omega, \mathcal{F}, \mathbb{P})$. Following Carmona and Delarue~\cite{car18}, any functional $U : \mathcal{P}_2(\mathbb{R}^d) \to \mathbb{R}$ can be lifted to a functional $\widetilde{U} : L^2(\Omega; \mathbb{R}^d) \to \mathbb{R}$ defined by:
\begin{equation*}
\widetilde{U}(X) := U(\mathcal{L}(X)), \quad X \in L^2(\Omega; \mathbb{R}^d),
\end{equation*}
where $\mathcal{L}(X) = X_\sharp \mathbb{P}$ denotes the law of $X$.

\begin{lemma}[Definition of the L-Derivative on $\mathcal{P}_2(\mathbb{R}^d)$]
The functional $U$ is said to be differentiable at $\mu \in \mathcal{P}_2(\mathbb{R}^d)$ in the Lions sense (or $L$-differentiable) if its lift $\widetilde{U}$ is Fréchet differentiable at some $X \in L^2(\Omega; \mathbb{R}^d)$ such that $\mathcal{L}(X) = \mu$.

\begin{enumerate}
    \item By Riesz's representation theorem, there exists a unique element $D\widetilde{U}(X) \in L^2(\Omega; \mathbb{R}^d)$ such that for any $Y \in L^2(\Omega; \mathbb{R}^d)$:
    \begin{equation*}
    \lim_{\varepsilon \to 0} \frac{\widetilde{U}(X + \varepsilon Y) - \widetilde{U}(X)}{\varepsilon} = \mathbb{E}\big[ \langle D\widetilde{U}(X), Y \rangle \big].
    \end{equation*}
    
    \item There exists a Borel measurable function $\partial_\mu U(\mu) : \mathbb{R}^d \to \mathbb{R}^d$ belonging to $L^2(\mathbb{R}^d, \mu; \mathbb{R}^d)$, uniquely defined $\mu$-almost everywhere and independent of the choice of representative $X$, such that:
    \begin{equation*}
    D\widetilde{U}(X) = \partial_\mu U(\mu)(X) \quad \text{a.s.}.
    \end{equation*}
    The mapping $x \mapsto \partial_\mu U(\mu)(x)$ is called the Lions derivative (or $L$-derivative) of $U$ at $\mu$.
    \item We assume that $(\mu,x)\rightarrow \partial_\mu U(\mu)(x)$ is jointly Lipschitz, then:
    \begin{equation}\label{expansion measure Lions}
    U(\mu)=U(\mu_0)+\E[\partial_\mu U(\mu)(X)(X-X_0)]+o(\|X-X_0\|_2^2).
    \end{equation}
\end{enumerate}
\end{lemma}


\end{document}